\documentclass[10pt]{article}
\usepackage{amsfonts,amssymb,amsmath,a4wide,xypic,remreset,amscd}

\def\proof {\noindent{\sc{Proof. }}}
\def\qed {\mbox{}\hfill {\small \fbox{}} \\}  
\newcommand{\Zm}{\ensuremath{\mathbb{Z}}}
\newcommand{\Om}{\ensuremath{\mathbb{O}}}
\newcommand{\Sm}{\ensuremath{\mathbb{S}}}

\newcommand{\Rm}{\ensuremath{\mathbb{R}}}
\newcommand{\Vm}{\ensuremath{\mathbb{V}}}
\newcommand{\Tm}{\ensuremath{\mathbb{T}}}
\newcommand{\Pm}{\ensuremath{\mathbb{P}}}
\newcommand{\Mm}{\ensuremath{\mathbb{M}}}

\newcommand{\Nm}{\ensuremath{\mathbb{N}}}
\newcommand{\Em}{\ensuremath{\mathbb{E}}}

\newcommand{\Um}{\ensuremath{\mathbb{U}}}

\newcommand{\Dm}{\ensuremath{\mathbb{D}}}
\def\lto{\longrightarrow}
\def\lmto{\longmapsto}

\def\leq{\leqslant}
\def\geq{\geqslant}
\newcommand{\fleche}{\ensuremath{\vdash}}
\newcommand{\ffleche}{\mbox{\ensuremath{
\dashv \vdash}}}

\newcommand{\mX}{\ensuremath{\mathcal{X}}}

\newcommand{\mI}{\ensuremath{\mathcal{I}}}
\newcommand{\mK}{\ensuremath{\mathcal{K}}}

\newcommand{\mG}{\ensuremath{\mathcal{G}}}

\newcommand{\mN}{\ensuremath{\mathcal{N}}}
\newcommand{\mM}{\ensuremath{\mathcal{M}}}
\newcommand{\mS}{\ensuremath{\mathcal{S}}}
\newcommand{\mB}{\ensuremath{\mathcal{E}}}
\newcommand{\mA}{\ensuremath{\mathcal{A}}}

\newcommand{\mH}{\ensuremath{\mathcal{H}}}

\newcommand{\vs}{\vspace{.2cm}}

\makeatletter
\renewcommand{\section}{\@startsection
{section}
{1}
{0mm}
{-1.2\baselineskip}
{\baselineskip}
{\center \scshape}}

\renewcommand{\subsection}{\@startsection
{subsection}
{2}
{0mm}
{\baselineskip}
{0.5\baselineskip}
{\normalfont \normalsize \bfseries}}

\renewcommand{\subsubsection}{\@startsection
{subsubsection}
{3}
{0mm}
{-\baselineskip}
{0mm}
{\normalfont \normalsize \bfseries}}

\@removefromreset{subsection}{section}
\makeatother

\renewcommand{\thesubsection}{\arabic{subsection}}

\author{Patrick Bernard}

\title{}

\begin{document}
\begin{center}
\begin{scshape}
\begin{Large}
The Dynamics of Pseudographs in Convex Hamiltonian  Systems
\footnote{version 3, october 2007; 
version 1 submitted in october 2004}
\vspace{1cm}\\
\end{Large}
\begin {large}
Patrick Bernard
\end{large}
\end{scshape}
\end{center}

JAMS, 21 (2008) 615-669.

\begin{small}
\textsc{abstract. }
We study the evolution, under convex Hamiltonian flows on cotangent bundles
of compact manifolds, of certain distinguished subsets of the phase space.
These subsets are generalizations of Lagrangian graphs,
we call them pseudographs.
They  emerge in a natural way from
Fathi's weak KAM theory.
By this method, we find various orbits which connect 
prescribed regions of the phase space.
Our study is inspired by works of John Mather.
As an application, we obtain the existence of diffusion 
in a large class of a priori unstable systems
and provide a solution to the large gap problem.
We hope that our method will have applications to more examples.
\\

\textsc{R\'esum\'e. }
Nous \'etudions  l'\'evolution, par le flot d'un Hamiltonien convexe
sur une vari\'et\'e compacte, de certains ensembles de l'espace
des phases.
Nous appelons pseudographes ces ensembles,
qui sont des g\'en\'eralisations de graphes
Lagrangiens apparaissant de mani\`ere naturelle dans la th\'eorie
KAM faible de Fathi.
Par cette m\'ethode, nous trouvons  diverses orbites
qui joignent des domaines donn\'es de l'espace des phases.
Notre \'etude s'inspire de travaux de John Mather.
Nous obtenons l'existence de diffusion dans une large classe
de syst\`emes \`a priori instables comme application de cette 
m\'ethode, qui permet de r\'esoudre le probleme de l'\'ecart 
entre les tores invariants.  
Nous esp\'erons que la m\'ethode s'appliquera \`a d'autres exemples.

\tableofcontents
\vspace{2cm}

\newpage
\end{small}

\section{Introduction}
In all this paper, $M$ denotes a 
connected compact manifold without boundary, of dimension $d$,
and $TM$ and $T^*M$ are its tangent and cotangent bundle.
We shall consider the periodic   time-dependent Hamiltonian system
generated by a function $H:\Rm\times T^*M \lto \Rm$,
and denote by $\phi_s^t$ the flow from time $s$ to time $t$.

\subsubsection{}\label{question}
In order to motivate our discussion, we  begin with 
a precise question:
Given two Lagrangian manifolds $\mG$ and $\mG'$ in the cotangent
bundle, which are graphs over the base $M$, does there exist
a trajectory which connects $\mG$ and $\mG'$, or in other words does there
exist times $s<t$ such that the Lagrangian manifold 
$\phi_s^t(\mG)$ intersects $\mG'$?

\subsubsection{}
This question formulates some well known problems.
As an example, let us suppose that $M=\Tm^d$, 
and identify the cotangent bundle 
$T^*\Tm^d$ with $\Tm^d \times \Rm^d$.
Let us consider the  Hamiltonian $H_0=h(p)$,
where  $h:\Rm^d\lto \Rm$ is a real function.
Such Hamiltonians will be called fully integrable 
in the sequel.
It is known that they leave invariant  
the tori $\Tm_p:=\Tm^d\times\{p\}$, for $p\in \Rm^d$.
As a consequence, the answer to the previous question
is obviously negative for $\mG=\Tm_p$ and $\mG'=\Tm_{p'}$,
when $p\neq p'$.
What happens for Hamiltonians $H$ which are close to $H_0$?
For example, it is known that the solar system can be described by a fully
integrable Hamiltonian $H_0$ if the interactions between planets
are neglected. 
In this example, the variables $p\in \Rm^d$
encode the parameters of the elliptic trajectories of the planets.
It is well known that these parameters would not change in time
if the interaction between planets did not exist.
Understanding for which values of $p$ and $p'$ 
the question \ref{question}
has a positive answer with 
$\mG=\Tm_p$ and $\mG'=\Tm_{p'}$,
amounts to understand to what extent the elliptic trajectories
will deform under the influence of mutual interactions.
In other words, it amounts to understand the secular dynamics,
and the stability of the solar system.
We will not treat these specific examples in the present
papers, although they are parts of our motivations.
See \cite{Arnold} and \cite{MS} for beautiful and deep
examples of 
perturbations of fully integrable systems.

\subsubsection{}
Question \ref{question}
is especially interesting when the Lagrangian manifolds $\mG$
and $\mG'$ have different Liouville classes
(which corresponds to the case $p\neq p'$ in the discussion
above).
In this case, we have a problem of 
non exact Lagrangian intersection, and it seems that  
the powerful tools developed to deal with 
exact intersections provide no interesting insight.
In order to study this problems,
we  make strong assumptions
on the Hamiltonian $H$, namely that it is convex, super-linear,
periodically time-dependent, and complete, see details in 
\ref{HypothesesH}.
We will define an equivalence relation, called forcing relation
and denoted by 
 $\ffleche$ on the set
$H^1(M,\Rm)$ of cohomology classes of Lagrangian graphs,
in such a way that, if $c\ffleche c'$,
(we will say that $c$ and $c'$ \textit{force} each over)
 then the answer to question
 \ref{question} is positive for each Lagrangian graphs
$\mG$ of cohomology $c$ and $\mG'$ of cohomology $c'$.
The definition of this equivalence relation is one of the major
ideas of the present paper.
The key point in considering an equivalence relation
is that local informations on the equivalence classes 
can be put together to obtain global information.
On the other hand, most of the mechanisms known so far to study questions
related to  \ref{question},
the theorem of Birkhoff for twist maps,
the geometric construction of Arnold,
as well as the variational construction of Mather,
 can be expressed in this unified setting
as local informations on the forcing classes
(the classes of equivalence of the relation $\ffleche$).
Our main goal in the present paper will be to detail this fact
and to study the local properties of the forcing classes.

\subsubsection{}\label{applicationintro}
In order to demonstrate the usefulness of our theory, let us provide 
an example.
Proofs and more general statements are given in section \ref{GAE}.
We take $M=\Tm\times \Tm^{d-1}$,
and denote by $(q,p)=(q_1,q_2,p_1,p_2)$ the points of $T^*M$,
where $q_1\in \Tm$, $q_2\in \Tm^{d-1}$,
 $p_1\in \Rm$, $p_2\in \Rm^{d-1}$.
We consider the time-periodic Hamiltonian
$$
H(t,q,p)
=
H_1(t,q_1,p_1)+
|p_2|^2-
V(q_2)F(t,q)
$$
and we assume that the conditions of 
convexity, super-linearity and completeness are satisfied.
In addition, we assume that $F:\Tm\times \Tm^d\lto \Rm$
takes positive  values, and that 
$V:\Tm^{d-1}\lto \Rm$ takes positive values except at a single
point, say $0$, where its takes the value $0$.
The manifold $\Tm\times\Rm:=\{q_2=0, p_2=0\}$
is then invariant under the Hamiltonian Flow. 
The restricted flow is generated by the restricted 
Hamiltonian $H_1$.
Under these hypotheses, it is not hard to prove (we will do it) that 
each rotational invariant circle of the restricted 
dynamics $H_1$ admits a homoclinic orbit.
We make two additional non-degeneracy assumptions:

(H1) The Hamiltonian $H_1$ is generic in the sense that 
its irrotational invariant circles of rational rotation number 
are completely periodic.
(We allow periodic circles in order to include the case where 
$H_1$ is integrable).

(H2) We assume a non-degeneracy hypothesis on the
set of action minimizing  homoclinic
orbits to  the invariant circles of $H_1$.
This hypothesis is detailed in section \ref{GAE}, 
it should be seen as analogous to the classical hypothesis
of transversality of the stable and unstable manifolds
in the construction of Arnold.
Although we expect in the future to prove that 
this condition is generic in some sense, we do not discuss
any genericity issue here.

Under these hypotheses, our abstract results imply the following.
\vs\\
\textsc{Theorem }
\begin{itshape}
If $P$ and $P'$ are given real numbers, there exists
a Hamiltonian trajectory $(q(t),p(t))$
and an integer $t\in\Nm$
such that $p_1(0)=P$ and $p_1(t)=P'$.
\end{itshape}

\subsubsection{}
The systems described in example \ref{applicationintro}
are \textit{a priori} unstable according
to the terminology in use in the world of Arnold's diffusion.
This is due to the presence of the distinguished
invariant manifold $\{q_2=p_2=0\}$, which in many situations
is normally hyperbolic.
It appears clearly in the fundamental paper of  Arnold, \cite{Arnold}
that the presence of such a hyperbolic invariant manifold
intersecting $\mG$ and $\mG'$ greatly favors a positive answer 
to question \ref{question}.
\textit{A priori} unstable systems have been widely studied
because they appear naturally in the perturbation of 
completely integrable systems, and are easier to deal with.

In the work of Arnold, it is also assumed that the restriction
of the dynamics 
to the hyperbolic manifold  is integrable, say $H_1=|p_1|^2$ in our example.
This means that this invariant manifold is foliated by invariant
tori which he called whiskered tori because of the presence of 
hyperbolicity.
These whiskered tori are the building blocks of Arnold's construction,
so that this second hypothesis was  very important.
The main point in our application is that we do not make this assumption.
We only assume that the restricted dynamics is generic,
in a clearly specified sense.

In the context of perturbations of fully integrable
systems, the restriction of the flow to the hyperbolic manifold is
close to integrable, and KAM theory implies the existence
of many whiskered tori. However, when computing precisely
the various quantities that appear in Arnold's construction,
one observes that there does not exist  enough tori in general.
More precisely, the gap between tori is too big, this
is the Large Gap problem, see for example \cite{Lochak} for a 
more precise explanation.

Overcoming this problem has long been considered as a major challenge.
While the classical approaches based on refinements on 
the scheme of Arnold were worked out in that direction,
new variational methods were introduced, 
by John Mather in \cite{MatherFourier}.
It is also worth mentioning the work of Bessi, \cite{Bessi},
where the results sketched by Arnold are proved using variational methods.
This paper contains one on the first relevant achievements
of variational methods in these kind of questions, and it has been 
very influential.
However, these variational methods were facing the same kind of 
difficulties as classical methods.
In several special instances, the Large Gap problem  can be bypassed 
 because for specific reasons there exist more whiskered
tori.
This remark has been exploited to obtain many non-trivial results
from Arnold's construction or variational methods.
 For example, orbits of unbounded 
speed where built in \cite{BT} using the scheme of Arnold.
A similar result had previously  been obtained
by  John Mather,
\cite{Matherconjecture}, using  variational methods,
see also \cite{Kal}.
Other works exploit the same remark in different directions, see 
for example \cite{BBB}, which elaborates on \cite{Bessi},
 and many other texts.

Solutions 
to the Large Gap problem
have recently  been given by Delshams, de la Llave and Seara,
see \cite{DLS}, and by Treschev, see 
\cite{Tre}
using elaborations on Arnold's method. The details in  these works
are far from simple.
Cheng and Yan have  also proposed a solution 
 using elaborations on the variational methods
initiated by Mather, see \cite{Cheng}, and  \cite{Cheng2},  as well 
as Z. Xia, see \cite{Xia1} and \cite{Xia}.
Compared to these papers, the spirit of our work is different.
We present mechanisms of instability which are more
general, but more abstract.
We present some examples 
for illustration, and in order 
to give the reader a hint of how the abstract mechanisms can be used,
but we do not try at that point to describe the more general applications. Neither do we discusss the genericity of our 
hypotheses.

The influence of John Mather's published and unpublished works
on  the developpement of variational approaches 
could not be overestimated. He has announced in 
\cite{Matherannonce} very deep results on the perturbation of fully
integrable systems in dimension $2$,
and given indications on proofs in various talks and lectures.
I hope that the  tools
developed in the present  paper will contribute to clarify
and extend  
these results.

\subsubsection{}
Let us now enter more precisely into matter.
Given a Lipschitz function $u:M\lto \Rm$
and a closed smooth form $\eta$ on $M$, we consider
the subset $\mG_{\eta,u}$ of $T^*M$ defined by
$$
\mG_{\eta,u}
=\big\{
(x,\eta_x+du_x), x\in M \text{ such that }  du_x \text{ exists} \big\}.
$$
We call the subset $\mG\subset T^*M$ an 
\textit{overlapping pseudograph} if there exists a closed smooth form
$\eta$ and a semi-concave function $u$ such that 
$\mG=\mG_{\eta,u}$. See Appendix \ref{semiconcave}
for the definition of semi-concave functions.
Each pseudograph $\mG$ has a well defined cohomology
$c(\mG)\in H^1(M,\Rm)$, see \ref{norm}, which is just the 
De Rham cohomology $[\eta]$ of the closed form appearing in the
definition of $\mG$.
We denote by $\Pm$ the set of overlapping pseudographs.
If $M=\Tm$ is a circle, 
then overlapping pseudographs are graphs of functions 
which have only  discontinuities with 
downward jumps, or in other words functions which
can be locally written as the sum
of a continuous and a decreasing function.
Such sets were introduced in \cite{KO}, where they
are used in very elegant proofs of many known
properties of Twist maps. 
In higher dimension, overlapping pseudographs naturally
arise from Fathi's approach of Mather theory.

\subsubsection{}
We define the forcing relation $\ffleche$ on $H^1(M,\Rm)$
as follows:
We say that
 $c$ and $c'$ force each other  (in short $c\ffleche c'$)
  if there exists an integer $N\in \Nm$
such that, for each pseudograph $\mG$ of cohomology $c$ (resp. $c'$),
there exists a pseudograph $\mG'$ of cohomology $c'$ (resp. $c$)
such that
$$
\overline{\mG'}\subset 
\bigcup_{1\leq i\leq N} \phi_0^i(\mG),
$$
where $\overline{\mG'}$ is the closure of $\mG'$ in $T^*M$.
This definition is certainly one of the most important novelties 
in the present paper.
Note that,
if $c\ffleche c'$, if $\mG$ is a Lagrangian graph of 
cohomology $c$, and if $\mG'$ is a Lagrangian graph of cohomology
$c'$, then  there exists a Hamiltonian orbit
which connects $\mG$ and $\mG'$.
As a consequence, understanding the equivalence classes of this
relation is a useful tool in the study of our motivating question.
Our main goal in the present paper will be to find sufficient
conditions for two classes to be equivalent.
It turns out that, although the definition seems very strong,
the existence of non-trivial forcing classes can be proved in many
 interesting situations, as example \ref{applicationintro}. 
In fact, many of the known constructions  of orbits connecting
prescribed regions of phase space 
(Birkhoff's theory of twist maps, Mather's 
variational construction of connecting
orbits, Arnold's geometric construction of diffusion)
can be adapted to this framework, and rephrased as the existence of 
large forcing classes.

\subsubsection{}
We shall define, following Fathi, an operator $\Phi:\Pm\lto \Pm$ in \ref{Phi},
with the following fundamental properties:
$$
\overline{\Phi(\mG)}\subset \phi(\mG),
$$
where $\phi:=\phi_0^1$ is the time-one map of the Hamiltonian flow,
and 
$c(\Phi(\mG))=c(\mG)$.
Fathi's weak KAM theorem,
\cite{Fathi1} states that, for each $c\in H^1(M,\Rm)$,
the operator $\Phi$ has fixed points of cohomology $c$. We 
call $\Vm_c$ the set of these fixed points, see section 
\ref{Aubrymather} for details.
The fixed points $\mG$ satisfy
$$
\overline{\mG}\subset \phi(\mG),
$$
and give rise to compact invariant sets 
$$
\tilde \mI(\mG):=\bigcap _{i\in \Nm}\phi^{-i}(\bar \mG).
$$
This provides a new way, due to Albert Fathi,
to define various invariant sets previously introduced  by Mather in 
\cite{Mather} and \cite{MatherFourier}.
 
\subsubsection{}
More precisely, to each cohomology $c\in H^1(M,\Rm)$
we associate the non-empty compact invariant sets 
$$
\tilde \mM(c)\subset \tilde \mA(c)\subset \tilde \mN(c),
$$
where 
$$
\tilde \mA(c):= \bigcap _{\mG\in \Vm_c}\tilde \mI(\mG)
\text{  and  }
\tilde \mN(c):= \bigcup _{\mG\in \Vm_c}\tilde \mI(\mG),
$$
are respectively called the Aubry set and the Ma\~n\'e set, 
and $\tilde \mM(c)$, called the Mather set, is the union
of the supports of the invariant measures of the action of $\phi$ on 
$\tilde \mA(c)$ (or equivalently on $\tilde \mN(c)$),
see \ref{ensembles} for more details.
A standing notation will be to denote by $\tilde \mX$
subsets of $T^*M$, and by $\mX$ their projection on $M$.

Beyond answering question \ref{question},
understanding the forcing relation $\ffleche$ has many dynamical consequences:

\subsubsection{}\label{orbites}
\textsc{Proposition.  }
\begin{itshape}
\begin{itemize}
\item[$(i)$]
Let $\mG$ and $\mG'$ be two Lagrangian graphs of cohomologies 
$c$ and $c'$. If $c\ffleche c'$ then there exists a time $t\in \Nm$
such that $\phi_0^t(\mG)$ intersects $\mG'$.
\item[$(ii)$]
If $c\ffleche c'$,
there exists a heteroclinic trajectory of the Hamiltonian flow
between $\tilde \mA(c)$ and $\tilde \mA(c')$.
\item[$(iii)$]
Let $c_i,i\in \Zm,$ be a sequence of cohomology classes
all of which force the others.
Fix, for each $i$ a neighborhood $U_i$ of $\tilde \mM(c_i)$
in $T^*M$.
There exists a trajectory of the Hamiltonian flow which visits
in turn all the sets $U_i$. In addition, 
if the sequence stabilizes to $c-$ on the left, or (and) to $c+$ on the right,
the trajectory can be assumed negatively asymptotic to   $\mA(c-)$
or (and) positively asymptotic to  $\mA(c+)$.\vs
\end{itemize}
\end{itshape}
The proof is given in section \ref{relation}.
Let us now state our main results which, as announced above,
describe the local structure of the forcing classes.

\subsubsection{}
For each $\mG\in \Vm$, we define
the subspace $R(\mG)$ of $H^1(M,\Rm)$
as the set of cohomology classes of smooth closed
one-forms whose support is disjoint from $\mI(\mG)$.
For each cohomology class $c\in H^1(M,\Rm)$,
we define the subspace $R(c)$ as 
$$
R(c)=\bigcap _{\mG\in \Vm_c}R(\mG)
\subset H^1(M,\Rm).\vs
$$
The following Theorem reformulates and extends  results of John Mather,
see \cite{MatherFourier} and also \cite{Fourier} and \cite{Cheng}.
It is proved in section \ref{mathersection}.
\vs\\
\textsc{Theorem. }\label{mather}
\begin{itshape}
For each $c_0\in H^1(M,\Rm)$, there exists
 a positive $\epsilon$ such that
the following holds:
Each class  $c\in H^1(M,\Rm)$
such that $c-c_0\in R(c_0)$ and
$\|c-c_0\|\leq \epsilon $
satisfies
$
c_0\ffleche c.
$\vs
\end{itshape}

In order to illustrate this result, it is useful to consider the case 
of twist maps $M=\Tm$. In this case, the reader should check that 
$R(c)=\Rm$ or $0$, and that 
 $R(c)=0$ if and only if there exists a rotational  invariant circle of
cohomology $c$, see (few) more details in section \ref{TM}.
The above result then roughly says that rotational invariant circles
are the only obstructions to the evolution of action variables, and recovers
the theory of Birkhoff.
We will  explained in section \ref{GAE}  how this
result allows to overcome the possible absence of invariant circles in example
\ref{applicationintro}.
\subsubsection{}\label{weakfinite}
There is a natural partition of the Aubry set $\tilde \mA(c)$
into compact invariant subsets $\tilde \mS$ called the static classes,
see section \ref{sectionstatic}. A generalized version of the following 
theorem is proved in section \ref{finitesection}.
\vs\\
\textsc{Theorem. }
\begin{itshape}
Assume that there exists only finitely many static classes 
in $\tilde \mA(c)$, and that the set 
$\tilde \mN(c)-\tilde \mA(c)$ is not empty and contains
finitely many orbits.
Then the cohomology $c$ is in the interior of its forcing class.\vs
\end{itshape}

This result may be seen as a reinterpretation in our langage
of the geometric construction of Arnold.
We explain in section \ref{GAE} how it allows to take 
into account the possible presence of invariant circles of $H_1$
in example \ref{applicationintro}.
We mention that, for a generic Lagrangian, all the Aubry sets
$\mA(c), c\in H^1(M,\Rm)$ have finitely many static classes,
see \cite{BG}.

\subsubsection{}
Let us now present the content of the paper.
The whole paper heavily relies on the notion of semi-concave 
function and of equi-semi-concave sets of functions. 
These notions are presented in Appendix  \ref{semiconcave}.
In Appendix \ref{uniform}, we prove some background
results, essentially due to Mather and Fathi, about the 
properties of the Action.

\textbf{Mather-Fathi Theory. \\}
This first part is a survey of the theory of Mather, 
Ma\~n\'e and Fathi of globally minimizing orbits,
from a point of view very close to the one of Fathi.
This survey is presented not only for the convenience of the reader,
but also because we need various  variations on existing
results,and also we need to recast the theory in our framework.
In section \ref{calculus}, we present the context, detail 
the standing hypotheses, and recall some known results 
of the calculus of variations which will be of constant use
(proofs are given in Appendix \ref{uniform}).
Pseudographs are defined and their basic properties studied in
Section \ref{pseudographs}.
In Section \ref{Aubrymather}, we use these pseudographs to
present Fathi's point of view on Aubry-Mather theory.
The theory is continued in section \ref{sectionstatic},
where we explain Ma\~n\'e's decomposition in static classes
of the Aubry set, and the construction of 
homoclinic orbits, due to
Fathi \cite{Fathi3}, Contreras and Paternain
\cite{CP}, (see also \cite{ETDS})
 which will play a central role in section 
\ref{finitesection}.

\textbf{Abstract mechanism. \\}
This part contains the main novelties of our paper.
In section \ref{relation}
we define the forcing relation $\ffleche$
and  explain how various orbits can
be built once this  relation  is understood.
We  prove Proposition \ref{orbites}.
We then introduce and study evolution operators on $\Pm$,
which are elaborations around the Lax-Oleinik operator,
in section \ref{operators}.
Section \ref{covering} is a parenthesis where we study
the action of taking finite Galois covering, which will
be essential for applications. 
The idea of taking finite Galois coverings comes from
 Fathi \cite{Fathi3}.
In section \ref{mathersection}, we prove and comment 
Theorem \ref{mather}.
In section 
\ref{finitesection}
we study the case where there exist only finitely many static classes.
We generalize and prove Theorem \ref{weakfinite}.

\textbf{Applications. \\}
In this short part, we detail 
some straightforward applications of the results
obtained before. We hope that it is possible to obtain 
much more applications by applying our results
with Hamiltonian methods such as normal form
theory, but this aspect is not discussed here.
Section 
\ref{TM} briefly mentions the application to twist maps.
Section 
\ref{GAE}
details \ref{applicationintro} above.
%
%
%
%
%
%
%
%
\section{Mather-Fathi Theory}
This part is an overview of the theory of Mather, Ma\~n\'e
and Fathi of globally minimizing orbits, which is oriented
towards our future needs.
We introduce our main objects.
Our point of view is close to the one of Fathi.
Most of the material  exposed here is a  small deformation of 
results in  \cite{Mather},
\cite{Fathibook}, \cite{Mane},  \cite{CDI}, or  \cite{CP}.

\subsection{Calculus of variations}\label{calculus}
\subsubsection{}\label{HypothesesH}
We shall consider $C^2$ Hamiltonian functions 
$H:\Rm\times T^*M\lto \Rm$.
We will denote by $P=(x,p)$ the points of $T^*M$.
The Cotangent bundle is endowed with its standard symplectic structure.
We denote by $X(t,P)$ or $X(t,x,p)$
the Hamiltonian vector-field of $H$, which is a time-dependent 
vector-field on $T^*M$.
We fix once and for all a Riemannian metric
$g$  on $M$, and use it to define
norms of tangent vectors and tangent covectors of $M$.
We will denote this norm  indifferently by $|P|$ or by $|p|$
when $P=(x,p)\in T^*_xM$.
We assume  the following standard set of hypotheses.
\begin{enumerate}
\item \textsc{Periodicity. } $H(t+1,P)=H(t,P)$ 
for each $(t,P)\in \Rm\times T^*M.$
\item
\textsc{Completeness. } The Hamiltonian vector-field $X$ 
generates a complete flow of diffeomorphisms  on $T^*M$. We denote by 
$\phi_s^t:T^*M\lto T^*M$ the flow from time $s$ to time $t$, 
and by $\phi$ the flow $\phi_0^1$.
\item \textsc{Convexity. } For each $(t,x)\in \Rm\times M$, the
function $p\lto H(t,x,p)$ is convex on $T_x^*M$,
with positive definite Hessian. Shortly,
$\partial^2_pH >0$.
\item \textsc{Superlinearity. }
For each $(t,x)\in \Rm\times M$, the function $p\lmto H(t,x,p)$
is super-linear, which means that 
$
\lim_{|p|\lto\infty} H(t,x,p)/|p|=\infty.$
\end{enumerate}

\subsubsection{}\label{HypothesesL}
We associate to the Hamiltonian $H$ a Lagrangian function
$L:\Rm\times TM\lto \Rm$
defined by 
$$L(t,x,v)=\sup _{p\in T^*_x M} p(v)-H(t,x,p).
$$
The fiberwise differential $\partial_p H$ of $H$ can be seen as a mapping
$$\partial_pH:\Rm\times T^*M\lto\Rm\times TM,$$
this mapping is a diffeomorphism, whose inverse is given by 
$$\partial_vL:\Rm\times TM\lto\Rm\times T^*M.$$
We have the relations
$L(t,x,v)
=\partial_vL(t,x,v)(v)-H(t,x,\partial_vL(t,x,v))$
and
$H(t,x,p)
=\partial_pH(t,x,p)(p)-L(t,x,\partial_pH(t,x,p)).$
The Lagrangian $L$ satisfies the following properties, which follow from 
the analogous properties of $H$:
\begin{enumerate}
\item \textsc{Periodicity. } $L(t+1,V)=L(t,V)$ 
for each $(t,V)\in \Rm\times TM.$
\item \textsc{Convexity. } For each $(t,x)\in \Rm\times M$, the 
function $v\lmto L(t,x,v)$ is a convex function on $T_xM$,
with positive definite Hessian. Shortly, $\partial^2_vL>0$.
\item \textsc{Superlinearity. }
For each $(t,x)\in \Rm\times M$, the function $v\lmto L(t,x,v)$
is super-linear on $T_xM$.
\end{enumerate}
See Appendix \ref{uniform}  for  comments related to these hypotheses.
The hypotheses listed above are very suitable to use the calculus
of variations.
\subsubsection{}\label{mini}
Let us fix two times $s>t$ in $\Rm$ and two points $x$ and $y$ in $M$.
Let $\Sigma(t,x;s,y)$ be the set of absolutely continuous
curves $\gamma:[t,s]\lto M$ such that
$\gamma(t)=x$ and 
$\gamma(s)=y$.
As usual, we define the action of the curve $\gamma$
as 
$
A(\gamma)=
\int _t^s L(\sigma,\gamma(\sigma),\dot\gamma(\sigma))\, d\sigma$.
It is known that, for each $C$, the set of curves $\gamma$
in $\Sigma(t,x;s,y)$
which satisfy $A(\gamma)\leq C$ is compact for the topology of uniform 
convergence. As a consequence, there exist curves minimizing the action.
Let us define 
the value 
$$
A(t,x;s,y)=
\min_{\gamma\in \Sigma( t,x;s,y)}
\int _t^s L(\sigma,\gamma(\sigma),\dot\gamma(\sigma))\, d\sigma,
$$
and let $\Sigma_m (t,x;s,y)$ be the set of curves in $\Sigma$
reaching the above minimum.
The set  $\Sigma_m (t,x;s,y)$ is not empty, and it is compact for the 
topology of uniform convergence.
Each curve  $\gamma(\sigma)\in \Sigma_m$
is $C^2$ and satisfies 
the Euler-Lagrange equations.
Setting 
$$
p(\sigma )=
\partial_v L(\sigma,\gamma(\sigma),\dot \gamma(\sigma)),
$$
which is equivalent to
$$
\dot \gamma(\sigma)=\partial_p H
(\sigma,\gamma(\sigma),p(\sigma)),
$$
these equations are 
$$
\dot p(\sigma) =\partial_x L(\sigma,\gamma(\sigma),\dot \gamma(\sigma))
=-\partial_x H(\sigma,\gamma(\sigma),p(\sigma))
$$
Hence the curve $(\gamma(\sigma),p(\sigma))$
is a trajectory of the Hamiltonian flow.

\subsubsection{}\label{surdiff}
For each minimizing curve $\gamma\in \Sigma_m (t,x;s,y)$,
we have 
$$
-p(t)=-\partial_v L(t, x, \dot \gamma(t))
\in \partial^+_xA(t,x;s,y),
$$
where 
$ \partial^+_xA(t,x;s,y)$
denotes the set of proximal super-differentials of 
$q\lmto A(t,q;s,y)$ at point $q=x$, see Appendix A.
We also have
$$
p(s)=\partial_v L(s, y, \dot \gamma(s))
\in
\partial^+_yA(t,x;s,y).
$$
For each $t'>t$, the set of functions
$(x,y)\lmto A(t,x;s,y),s\geq t'$
is  equi-semi-concave on $M\times M$,
hence equi-Lipschitz, see Appendix A. 
In addition, the three following properties are equivalent:
\begin{itemize}
\item[$(i)$] The set  $\Sigma_m (t,x;s,y)$ contains only one point.
\item[$(ii)$] The function $A(t,.;s,y)$ is differentiable at $x$.
\item[$(iii)$] The function  $A(t,x;s,.)$ is differentiable at $y$.
\end{itemize}
If these equivalent properties hold, and if $\gamma(\sigma)$ is the
unique  curve of $\Sigma_m (t,x;s,y)$,
then setting $p(\sigma)=\partial_v L(\sigma, y, \dot \gamma(\sigma))$,
we have 
$$
p(t)=-\partial_x A(t,x;s,y)
\text{ and }
p(s)=\partial_y A(t,x;s,y).
$$

\subsubsection{}
Let $\eta$ be a smooth  one-form.
We will see the form $\eta$ sometimes as 
a section of the cotangent bundle $\eta:M\lto T^* M$
and sometimes as a 
fiberwise linear 
function on the tangent bundle
$\eta:TM\lto \Rm$.
If the form $\eta$ is closed
then the diffeomorphism $\phi_{\eta}:(x,p)\lmto(x,p+\eta_x)$
of $ T^*M$ is symplectic.
The Hamiltonian  
$$H_\eta(t,x,p)=H\circ\phi_{\eta}(t,x,p)=H(t,x,p+\eta_x)
$$
satisfies our hypotheses.
The associated Lagrangian is 
$(L-\eta)(t,x,v)=L(t,x,v)-\eta_x(v),$
where $\eta$ is considered as a function on $TM$.
The following diagram commutes for each $t$.
$$
\xymatrix{
& {T^*M} \ar@{->}[dr]^{H}\\
TM \ar@{->}[ur]^{\partial_v L}
\ar@{->}[dr]_{\partial_v (L-\eta)}
& & \Rm\\
& T^*M \ar@{->}[uu]_{\phi_{\eta}}
\ar@{->}[ur]_{H_{\eta}} }
$$

\subsubsection{}
We will also consider the modified action 
$$
A_{\eta}(t,x;s,y)=
\inf_{\gamma\in \Sigma( t,x;s,y)}
\int _t^s L(\sigma,\gamma(\sigma),\dot\gamma(\sigma))
-\eta_{\gamma(\sigma)}(\dot\gamma(\sigma))\, d\sigma,
$$
which of course satisfies all the properties of  \ref{surdiff}, 
with the modified expressions
$$
\eta_x-p(t)
\in \partial^+_xA_{\eta}(t,x;s,y) \text{ and }
p(s)-\eta_y
\in
\partial^+_y A_{\eta}(t,x;s,y)
$$
when $(\gamma(\sigma),p(\sigma)):[s,t]\lto T^*M$
is the Hamiltonian trajectory associated to a curve
$\gamma\in \Sigma( t,x;s,y)$ minimizing $A_{\eta}$.

\subsubsection{}\label{etac}
Let $\Omega$ be the set of closed smooth forms on $M$.
It is useful
to fix once and for all a linear section $S$
of the projection $\Omega \lto H^1(M,\Rm)$.
In other words, $S$ is a linear mapping from $H^1(M,\Rm)$
to $\Omega$ such that $[S(c)]=c$.
We shall abuse notations and denote by $c$ the 
form $S(c)$, in such a way that the symbol  $c$ 
denotes either a cohomology class or a standard form 
representing this cohomology class.

\subsubsection{}
The following consequence of 
Appendix \ref{uniform} will be useful.
See appendix \ref{semiconcave} for the definition of equi-semi-concave.
\vs\\
\textsc{Proposition. }
\begin{itshape}
If $C$ is a bounded subset of $H^1(M,\Rm)$,
and $\epsilon$ is a positive number, 
the functions
$
A_{c}(s,.;t,.), c\in C, t\geq s+\epsilon
$
are equi-semi-concave on $M\times M$.
\end{itshape}
  %
%
%
%
%
%
\subsection{Overlapping pseudographs.}\label{pseudographs}
We present the main objects, overlapping pseudographs, 
and study some basic properties. The relevance of semi-concave
functions to this kind of problems was noticed by Albert Fathi.

%
%
%
%
%
%
\subsubsection{}
Given a Lipschitz function $u:M\lto \Rm$
and a smooth form $\eta$ on $M$, we define
the subset $\mG_{\eta,u}$ of $T^*M$ by
$$
\mG_{\eta,u}
=\big\{
(x,\eta_x+du_x), x\in M \text{ such that }  du_x \text{ exists} \big\}.
$$
We call the subset $\mG\subset T^*M$ an 
\textit{overlapping pseudograph} if there exists a smooth form
$\eta$ and a semi-concave function $u$ such that 
$\mG=\mG_{\eta,u}$. See Appendix \ref{semiconcave}
for the definition of semi-concave functions.
We shall denote by $\Pm$ the set of overlapping pseudographs.
Given a pseudograph $\mG$,
and a subset $U \subset M$, we will denote by $\mG_{|U}$ the set 
$
\mG_{|U}:=\mG\cap \pi^{-1}(U)
$.

\subsubsection{}\label{norm}
It is not hard to see that if 
an overlapping  pseudograph $\mG$ is represented in two ways as
$\mG_{\eta,u}$ and $\mG_{\mu,v}$, then the closed forms
$\eta$ and $\mu$ have the same cohomology in $H^1(M,\Rm)$.
As a consequence, it is possible to associate to each pseudograph
$\mG$ a cohomology $c(\mG)$, in such a way that
$$
c(\mG_{\eta,u})=[\eta].
$$
We will denote by $\Pm_c$ the set of overlapping pseudographs 
of cohomology $c$.
If $\mG$ is an overlapping pseudograph of
cohomology $c$, then $\mG$
can be represented in the form
$\mG=\mG_{c,u}$,
where $c$ is the standard form defined in \ref{etac}.
The function $u$ is then uniquely defined up to an additive constant.
As a consequence, denoting by $\Sm$ the set of semi-concave functions
on $M$, and by $\Pm$ the set of overlapping pseudographs, 
we have the identification
$$
\Pm=H^1(M,\Rm)\times \Sm/\Rm.
$$
This identification endows $\Pm$ with the structure
of a real vector space.
Let us endow the factor $\Sm/\Rm$ with the norm
$|u|=(\max u-\min u)/2$, which is the
norm induced from the uniform norm on $\Sm$.
More precisely, we have 
$|u|=\min _v \|v\|_{\infty}$, where the minimum is taken on functions
$v$ which represent the class $u$.
We put on $\Pm$ the norm 
$$
\|\mG_{c,u}\|=|c|+(\max u-\min u)/2\leq|c| + \|u\|_{\infty}.
$$
The set $\Pm$ is now a normed vector space.
It is also useful to define, for each subset $U\subset M$, 
the number
$$
\|\mG_{c,u}\|_U:=|c|+(\sup_U u-\inf_U u)/2.
$$
We define in the same way the set $\breve \Pm$ of 
\textit{anti-overlapping} pseudographs $\breve \mG$, which are the 
sets $\mG_{\eta,-u}$, with $\eta$ a smooth form
and $u\in \Sm$. This set is similarly a normed vector space.

\subsubsection{}\label{intersection}
\textsc{Lemma. }
\begin{itshape}
Let $\mG$ be an overlapping pseudograph, 
and $\breve \mG$ be an anti-overlapping pseudograph.
If $\mG$ and $\breve \mG$ have the same cohomology,
then they have nonempty intersection.
\end{itshape}
\vs\\
\proof
Let us write $\mG=\mG_{\eta,u}$
and $\breve \mG = \mG_{\eta,-v}$.
Let $x\in M$ be a point minimizing the continuous function  $u+v$.
Since they are semi-concave,
both $u$ and $v$ are differentiable at $x$, 
and $du_x=-dv_x$.
It follows that the point $(x,\eta_x+du_x)=(x, \eta_x-dv_x)$
belongs both to $\mG$ and to $\breve \mG$.
\qed
It is natural to introduce the following definition.
\vs\\
\textsc{Definition. }
\begin{itshape}
Let $\mG$ be an overlapping pseudograph, 
and $\breve \mG$ be an anti-overlapping pseudograph.
If $\mG$ and $\breve \mG$ have the same cohomology $c$,
write them $\mG=\mG_{c,u}$ and 
$\breve \mG=\mG_{c,\breve u}$.
We denote by 
$$
\mG \wedge \breve \mG \subset M
$$
the set of points of minimum of the difference
$u-\breve u$, and by 
$
\mG\tilde \wedge \breve \mG \subset \mG\cap \breve \mG
$
the set 
$$
\mG\tilde \wedge \breve \mG :=
\mG\cap \pi^{-1}(\mG \wedge \breve \mG )=
\breve \mG\cap \pi^{-1}(\mG \wedge \breve \mG )=
\mG\cap\breve \mG\cap \pi^{-1}(\mG \wedge \breve \mG )
\subset T^*M.
$$
This  set 
is compact,  not empty, and it is a Lipschitz graph
over its projection 
$
\mG  \wedge \breve \mG.
$
\end{itshape}
\vs\\
\proof
We have proved already that the set $\mG\wedge \breve \mG$
is not empty.
It follows from Appendix \ref{sum}
that both $u$ and $\breve u$ are differentiable on 
$\mG\wedge \breve \mG$, and that the map 
$x\lmto du_x=-d\breve u_x$
is Lipschitz on this set.
This makes the definition meaningful.
The set 
$
\mG \tilde  \wedge \breve \mG 
$
is compact because it is the image of the compact set 
$
\mG  \wedge \breve \mG 
$
by a Lipschitz map.
\qed

\subsubsection{}
Let us fix a closed form $\eta$.
We define the associated Lax-Oleinik mapping  on 
$C^0(M,\Rm)$ by the expression
$$
T_{\eta}u(x)=
\min _{q\in M}\big(
u(q)+A_{\eta}(0,q;1,x)\big)
$$
Let us recall some important properties
of the Lax-Oleinik mapping,
which are all direct consequences of 
the properties of the function $A$ given in \ref{surdiff}.
For each form $\eta$, 
The functions $T_{\eta}u, u\in C(M,\Rm)$
are equi-semi-concave, see
Appendix \ref{semiconcave}.
The mapping $T_{\eta}$ is a contraction:
$$
\|T_{\eta}u-T_{\eta}v\|_{\infty}
\leq \|u-v\|_{\infty}.
$$
To finish, the mapping $T_{\eta}$ is non-decreasing,
and it satisfies
$T_{\eta}(a+u)=a+T_{\eta}(u)
$
for all real $a$.

\subsubsection{}\label{Phi}
There exists a unique mapping  $\Phi:\Pm\lto \Pm$ 
such that 
$$
\Phi(\mG_{\eta,u})=\mG_{\eta,T_{\eta}u}
$$
for all smooth form $\eta$ and all semi-concave function $u$.
We have 
$$
c(\Phi(\mG))=c(\mG).
$$
The mapping $\Phi$ is continuous
(see \ref{phiu} for the proof of a more general result).
For each cohomology $c$, the image $\Phi(\Pm_c)$ is a
 relatively  compact subset
of $\Pm_c$, as follows directly from the properties of the 
Lax-Oleinik transformation  recalled above.
Note that this implies the existence of a fixed point
of $\Phi$ in each $\Pm_c$,
this is how Fathi proved the existence of fixed points.
See \ref{liminf} for another proof.
We call $\Vm_c$ the set of these fixed points, and 
$\Vm=\cup_c\Vm_c$.
We also define the sets
$$
\Om :=\bigcap _{n\in \Nm} \overline{\Phi^{n}(\Pm)}
=\bigcap _{n\in \Nm} \Phi^{n}(\Pm)
$$
and  $ \Om _c:=\Om \cap \Pm_c$.
Note that $\Om_c$ is compact and invariant under $\Phi$,
and that $\Vm\subset \Om$.
A pseudograph $\mG\in \Pm_c$ belongs to $\Om$ if and only if
there exists a sequence  $\mG_n\in \Pm, n\in \Zm$ of pseudographs
such that $\Phi^{m-n}(\mG_n)=\mG_m$ for all $m\geq n$, and such that 
$\mG_0=\mG$.
Note that we then have $\mG_n\in\Om_c$ for each $n\in \Zm$.

\subsubsection{}\label{inclusion}
The mapping $\Phi$ satisfies
$$
\overline {\Phi(\mG)}\subset \phi(\mG).
$$
This  inclusion is a consequence of the following Proposition, 
which will be central througrought the paper.

\subsubsection{}\label{initial}
\textsc{Proposition }
\begin{itshape}
Let us fix a
pseudograph $\mG_{\eta,u}\in \Pm$, an open set $U\subset M$
and two times $s<t$.
Let us set 
$$
v(x)=\min _{q\in \bar U}u(q)+A_{\eta}(s,q;t,x),
$$
where $\bar U$ is the closure of $U$.
Let  $V\subset  M$ be an open set and let $N\subset M$ be the
set of points where the minimum is reached in the definition
of $v(x)$ for some $x\in V$. If $\bar N\subset U$, 
then 
$$
\overline{\mG_{\eta,v|V}}\subset\phi_s^t\big(\mG_{\eta,u|\bar N}\big)
$$
and $\mG_{\eta,u|\bar N}$ is a Lipschitz graph 
above $\bar N$.
In other words, the function $u$ is differentiable at each point
of $\bar N$, and the mapping $q\lmto du_q$ is Lipschitz on $\bar N$.
\end{itshape}
\vs\\
\textsc{Addendum. }
\begin{itshape}
In addition, the Hamiltonian 
trajectories  $(x(\sigma),p(\sigma)):[s,t]\lto T^*M$
 which terminate
in $\overline{\mG_{\eta,v|V}}$, \textit{i. e. } such that
$(x(t),p(t))\in \overline{\mG_{\eta,v|V}}$
satisfy
$$
v(x(t))=u(x(s))+
\int_s^t L(\sigma,x(\sigma),\dot x(\sigma))-
\eta_{x(\sigma)}(\dot x(\sigma))d\sigma
$$
$$
=
u(x(s))+A_{\eta}(s,x(s);t,x(t)))
=
\min_{x\in U}
u(x)+
A_{\eta}(s,x;t,x(t))).
$$
\end{itshape}
\vs\\
\proof
Let us fix a point  $x\in V$,
and  consider a point $q\in N$ minimizing in the expression 
of $v(x)$.
Since $q$ is a point of local minimum of the function
$u+A_{\eta}(s,.,t,x)$, the semi-concave functions
 $u$ and $A_{\eta}(s,.,t,x)$
are differentiable at $q$ and satisfy $du_q+\partial_q A_{\eta}(s,q,t,x)=0$.
In view of  \ref{surdiff},
we have 
$\partial_q A_{\eta}(s,q,t,x)=\eta_q-p(s)$,
where 
$$
(x(\sigma), p(\sigma))=
(x(\sigma),\partial_v L(\sigma,x(\sigma),\dot x(\sigma))):[s,t]\lto T^*M
$$
is the Hamiltonian 
trajectory associated to the unique minimizing curve
$x(\sigma)\in \Sigma_m(s,q,t,x)$.
Uniqueness follows from the differentiability 
of $A_{\eta}(s,.,t,x)$ at $q$, see
\ref{surdiff}.
We have
$$
(x(s),p(s))=(q,du_q+\eta_q)\subset \mG_{\eta,u}
$$
and therefore 
$$(x(t),p(t))=\phi_s^t(q,du_q+\eta_q)\in \phi_s^t
\big(\mG_{\eta,u|N}\big).
$$
Since the functions $A_{\eta}(s,.,t,x),x\in M$
are equi-semi-concave, they are all $K$-semi-concave for some $K$.
It follows that the function $u$ has a $K$-sub-differential
at each point of $N$, and therefore at each point of $\bar N$.
We conclude using \ref{regularity} that the function
$u$ is differentiable on $\bar N$, and that
the map $q\lmto du_q$ is Lipschitz on $\bar N$.
As a consequence, we have 
$$ 
\mG_{\eta,u|\bar N}=\overline{\mG_{\eta,u| N}},
$$
and this set is a Lipshitz graph over $\bar N$.
Still exploiting \ref{surdiff},
we get that the function 
$A_{\eta}(s,q,t,.)$ is  differentiable 
at $x$, and satisfies
$\partial_x A_{\eta}(s,q,t,x)=p(t)-\eta_x$.
Noticing that  the function
$v-A_{\eta}(s,q,t,.)$ has a local maximum at $x$,
we conclude that $dv_x=p(t)-\eta_x$
if $v$ is differentiable at $x$, and therefore that 
$$
(x,\eta_x+dv_x)=(x(t),p(t))\in \phi_s^t
\big(\mG_{\eta,u|N}\big)
$$
for each point of differentiability $x\in V$ of $v$.
In other words, we have the inclusion 
$\mG_{\eta,v|V}
\subset \phi_s^t\big(\mG_{\eta,u|N}\big)$, hence 
$$
\overline{\mG_{\eta,v|V}}
\subset \phi_s^t\big(\overline{\mG_{\eta,u|N}}\big)
=\phi_s^t\big({\mG_{\eta,u|\bar N}}\big).
$$
\qed

\subsubsection{}\label{calibrated}
Let $\mG=\mG_{c,u}$ be a fixed point of $\Phi$.
And let $n< m$ be two relative integers.
Following Fathi, we say that a curve $x(t):[n,m]\lto M$
is \textit{calibrated} by $\mG$ if
$$
u(x(n))+
\int_n^m L(t,x(t),\dot x(t))-
c_{x(t)}(\dot x(t))dt=T^{n-m}_{c}u(x(m)),
$$
note, since $u$ is a fixed point of $\Phi$,
that $T^{n-m}_{c}u=u+(m-n)\alpha(c)$,
where $\alpha(c)$ is a  constant (which, as we will see below
depends on $c$ but not on $u$). 
A consequence of the addendum in \ref{initial}
is that 
the curve $x(t)$ is calibrated by $\mG$ if the 
curve 
$(x(t),p(t)=\partial_vL(t,x(t),\dot x(t)))$ 
is a Hamiltonian trajectory satisfying $(x(m),p(m))\in \bar \mG$.
Conversely, if $x(t)$ is calibrated by $\mG$, then  
 $(x(k),p(k))\in  \mG$
 for each integer $k\in [n,m[$.

\subsubsection{}
The following Corollary is the reason why we have called the elements
of $\Pm$ \textit{overlapping}.
\vs\\
\textsc{Corollary }
\begin{itshape}
All   overlapping pseudographs $\mG\in \Pm$ satisfy 
$
\pi\circ \phi(\mG)=M.
$
\end{itshape}
\vs\\
\proof
We have $\overline{\Phi(\mG)}\subset \phi(\mG)$,
and $\pi(\Phi(\mG))$ is dense in $M$, so that 
$\pi\big(\overline{\Phi(\mG)}\big)=M.$
\qed

\subsubsection{}
It is useful, still following Fathi,  to define ''dual'' concepts.
We define the dual  Lax-Oleinik operator associated to a closed form $\eta$
by the expression
$$
\breve T_{\eta}u(x)=
\max _{q\in M}\big(
u(q)-A_{\eta}(0,x;1,q)\big),\,\,\, u\in C(M,\Rm)
$$
and we associate to this operator
a mapping $\breve \Phi:\breve \Pm\lto \breve \Pm$
by the expression 
$
\breve \Phi(\mG_{c,-u})=\mG_{c,\breve T_c(-u)}\in \breve \Pm
$.
We have 
$$
\overline{\breve \Phi(\breve \mG)}
\subset \phi^{-1}(\breve \mG)
$$
if $\breve \mG\in \breve \Pm$. 
We denote by $\breve \Vm$ the set of fixed points of 
$\breve \Phi$.
Let $\breve \mG=\mG_{c,-u}$ be a fixed point of $\breve \Phi$,
and let $n< m$ be two relative integers.
We say that a curve $x(t):[n,m]\lto M$
is calibrated by $\breve \mG$ if
$$
u(x(m))-
\int_n^m L(t,x(t),\dot x(t))-
c_{x(t)}(\dot x(t))dt
=\breve T ^{m-n}_c u(x(n)).
$$

\subsubsection{}\label{inequal}
It is useful to collect a few remarks.
We have 
$$\breve T_cT_c u\leq u
$$
for all continuous functions $u$.
Writing $\breve T_c^{n+1}T_c^{n+1} u=
\breve T_c^n  \breve T_cT_c  T_c^nu$, we observe that  
$ \breve T_c^nT_c^n u$ is a  non-increasing sequence of functions.
Conversely, 
$T_c^n\breve T_c^n u
$
is non-decreasing for each continuous function $u$.
 %
%
%

%
%
%
%
%
%
\subsection{Aubry-Mather sets}\label{Aubrymather}
We use the overlapping pseudographs to recover various 
invariant sets introduced by Mather, and to study their
major properties. We also establish the equivalence 
between the different definitions of the same
sets given in the literature.
Most of the section follows Fathi \cite{Fathibook},
with some minor variations. 

\subsubsection{}\label{alpha}
\textsc{Proposition}
\begin{itshape}
There exists a function $\alpha:H^1(M,\Rm)\lto \Rm$
such that, for each continuous function $u$ and each form
$\eta$ of cohomology $c$, the sequence
$
T^n_{\eta}u(x) +n\alpha(c),n\geq 1
$
of continuous functions is equi-bounded and equi-Lipschitz.
The function $c\lmto \alpha(c)$ is convex and 
super-linear.
More precisely, there exists a constant $K(c)$, which does not depend
on the continuous function $u$, 
such that
$$
\min u -K(c)\leq
T^n_cu(x)+n\alpha(c)\leq \max u +K(c).
$$
for each $n\in \Nm$ and $x\in M$.
\end{itshape}
\vs\\
\proof
Let us fix a cohomology class $c$, and define
the sequences 
$$
M_n(c):=\max _{x\in M} T_c^n(0)(x)
\text{ and }
m_n(c):=\min _{x\in M}  T_c^n(0)(x),
$$
where $0$ is the zero function on $M$.
Since the functions  $T_c^n(0), n\geq 1,$ are  
equi-semi-concave, see Appendix \ref{semiconcave}, 
there exists a constant $K$ such that 
$$
0\leq M_n(c)-m_n(c)\leq K
$$
for $n\geq 1$.
We claim that  $M_{n+m}(c)\leq  M_n(c)+M_m(c)$.
This follows from the inequalities
$$
 T^{m+n}_c(0)(x)=T_c^m(T_c^n(0))(x)\leq
T_c^m(M_n(c))(x)\leq M_n(c)+ T_c^m(0)(x).
$$
Hence by a classical result on subadditive sequences,
we have 
$
\lim M_n(c)/n=\inf M_n(c)/n.
$
We denote by $-\alpha(c)$ this limit.
In the same way, the sequence 
$-m_n(c)$ is subadditive, hence 
$m_n(c)/n\lto \sup m_n(c)/n$.
This limit is also $-\alpha(c)$ because 
$
0\leq M_n(c)-m_n(c)\leq K.
$
Note that $m_1(c)\leq -\alpha(c)\leq M_1(c)$,
so that $\alpha(c)$ is indeed a finite number.
We have, for all $n\geq 1$,
$$
-K-n\alpha(c)\leq m_n(c)\leq -n\alpha(c)\leq M_n(c)\leq K-n\alpha(c).
$$
Now far all $u\in C(M,\Rm)$, $n\in \Nm$ and $x\in M$, we have
$$
\min _{ M}u-K
\leq 
\min _{ M}u+m_n(c)+n\alpha(c)
\leq
T^n_cu(x)+n\alpha(c)
\leq 
\max _{ M}u+M_n(c)+n\alpha(c)
\leq 
\max _{ M}u+K,
$$
and we obtain the first conclusion of the Proposition.
The explicit definition of the value $m_n(c)$
is
$$
m_n(c)=\min_{\gamma\in C^1([0,n],M)} 
\int _0^n L(s, \gamma(s),\dot \gamma(s))-
c_{\gamma(s)}(\dot \gamma(s))ds.
$$
As a consequence, the function $c\lmto m_n(c)$
is concave, as a minimum of linear functions.
Hence each of the functions $c\lmto m_n(c)/n$ is concave,
so that the limit $-\alpha(c)$ is concave, and
the function $\alpha(c)$ is convex.
Since $\alpha(c)\geq K-m_1(c)$, it is enough to prove
that $-m_1$ is super-linear as a function of $c$ in order to prove that
$\alpha$ is.
For each homology class $h\in H_1(M,\Zm)$, let 
$\gamma_h:[0,1]\lto M$ be a closed curve representing this homology class.
We have 
$$-m_1(c)\geq c(h)-\int_0^1 L(s, \gamma_h(s),\dot \gamma_h(s))ds.
$$
This implies that $-m_1$, hence $\alpha$, is super-linear.
Indeed, in order that a function $f:\Rm^n\lto \Rm$
is super-linear, it is enough that there exists,
for each $y\in \Zm^n$, a value $a_y$ such that 
$f(x)\geq y\cdot x-a_y$ for each $x$.
\qed
\subsubsection{}\label{liminf}
\textsc{Proposition. }
\begin{itshape}
Let us fix a closed form $\eta$ and a continuous function
$u$.
Let us set 
$$
v:=\liminf _{n\lto\infty}\big( T_{\eta}^n(u)+n\alpha([\eta])\big),
$$
then $v$ is a fixed point of $T_{\eta}+\alpha$ hence
$\mG_{\eta,v}$ is a fixed point of $\Phi$. 
\vs
\end{itshape}

\proof
The one-form $\eta$ is fixed once and for all in this 
proof,  we omit the subscript $\eta$, and denote by $\alpha$
the number $\alpha([\eta])$.
Let us first prove that
$Tv+\alpha\leq v$.
In order to do so, we 
fix $x\in M$ and consider 
 an increasing sequence $n_k$
of integers such that $T^{n_k}u(x)+n_k\alpha\lto v(x)$.
Let $q_k$ be a point such that 
$T^{n_k}u(x)=T^{n_k-1}u(q_k)+A(0,q_k;1,x)$
or equivalently,
$T^{n_k}u(x)+n_k\alpha =T^{n_k-1}u(q_k)+(n_k-1)\alpha+\alpha +
A(0,q_k;1,x)$.
We can suppose that the sequence $q_k$ has a limit $q$.
Taking the $\liminf$ in the equality above gives 
$$
v(x)\geq v(q)+A(0,q,1,x)+\alpha \geq Tv(x)+\alpha
$$
where we have used the equi-continuity of the functions 
$T^nu,n\in\Nm$.

In order to prove that $Tv+\alpha \geq v$,
just notice that 
$
T^nu(x)\leq T^{n-1}u(q)+A(0,q;1,x)
$
for each $q$ and $x$,
or equivalently that 
$
T^nu(x)+n\alpha \leq T^{n-1}u(q)+(n-1)\alpha+A(0,q;1,x)+\alpha
$
and take the liminf.
\qed

\subsubsection{}\label{minimum}
\textsc{Lemma}
\begin{itshape}
Let us fix a closed form $\eta$ of cohomology $c$.
Let $\Mm\subset C(M,\Rm)$ be a family of  fixed points
of  the Lax-Oleinik operator $T_{\eta}+\alpha(c)$.
Assume that the infimum $v(x)=\inf_{u\in \Mm}u(x)$
is finite  for one (and then each)  $x\in M$.
Then the function $v$ is a fixed point of  $T_{\eta}+\alpha(c)$.
\vs\\
\end{itshape}
\proof
For all functions $u\in \Mm$ and all points $x$ and $y$ in $M$,
we have $u(x)\leq  u(y)+A_{\eta}(0,y;1,x)+\alpha(c)$.
It follows that
$
v(x)\leq v(y)+A_{\eta}(0,y;1,x)+\alpha(c),
$
so that 
$$
v(x)\leq \inf_y \big( v(y)+A_{\eta}(0,y;1,x)\big)+\alpha(c).
$$
In order to prove the other inequality, let us
fix $\epsilon>0$, take a 
function $u\in \Mm$ such that $u(x)\leq v(x)+ \epsilon$,
 and consider a 
point $y\in M$ such that 
 $u(x)= u(y)+A_{\eta}(0,y;1,x)+\alpha(c)$.
We obtain 
$$
v(x)\geq u(x)-\epsilon
\geq  u(y)+A_{\eta}(0,y;1,x)+\alpha(c)-\epsilon
\geq v(y)+A_{\eta}(0,y;1,x)+\alpha(c)-\epsilon.
$$
As a consequence, we have $v(x)\geq T_{\eta}v(y)±\alpha(c)-\epsilon$,
and, since this holds for all $\epsilon>0$, the desired inequality follows.
\qed

\subsubsection{}
Fixed points of the Lax-Oleinik operator $T_c+\alpha(c)$
will be called weak KAM solutions, following Fathi.
We denote by $\Vm\subset \Pm$ the set of fixed points of $\Phi$,
and $\Vm_c\subset \Pm_c$ the set of fixed points of $\Phi$ of
cohomology $c$.
Sometimes, we will also denote by $\Vm_C$ the set of fixed points
of $\Phi$ whose cohomology belongs to the subset 
$C\subset H^1(M,\Rm)$.
The set $\Vm_c$ is non-empty for each $c$.
If $\mG\in \Vm$, then it follows from \ref{inclusion} that
$$
\overline{\mG}\subset \phi(\mG).
$$
It is then natural to define the set 
$$
\tilde \mI(\mG)=\bigcap_{n\in \Nm} \phi^{-n}(\bar \mG),
$$
which is a non-empty compact $\phi$-invariant subset of $T^*M$.
We also define 
$$
\mI(\mG)=\pi(\tilde \mI(\mG))\subset M.
$$
More generally, for each $\mG\in \Pm$,
we define the set 
$$
\tilde \mI(\mG):=
\bigcap _{n\in \Nm}
\phi^{-n}\big(\overline {\Phi^n(\mG)}\big).
$$
Since 
$
\phi^{-n}(\overline{\Phi^{n}(\mG)}),
$
is a non-increasing
 sequence of compact sets, the set $\tilde \mI(\mG)$
is compact and not empty for each $\mG\in \Pm$.

\subsubsection{}\label{ensembles}
For each $\mG\in \Vm$, we  define
the set $\tilde \mM(\mG)$ as  the union of the supports 
of invariant measures
of $\phi_{|\tilde \mI(\mG)}$.
If $\mG\in \Vm$ and $\mG'\in \Vm$ have the same cohomology $c$,
then it is known that
$$
\tilde \mM(\mG)\subset \tilde \mI(\mG')
$$
hence $\tilde \mM(\mG)=\tilde \mM(\mG')$.
As a consequence, the set $\tilde \mM$, usually called
the Mather set, depends only on the cohomology $c$.
It will be denoted by  
$$\tilde \mM(c),$$
and as usual, we will denote by  $\mM(c)$ 
the projection  $\pi(\tilde\mM(c))$.
We also define the Aubry set in a usual way
by
$$
\tilde \mA(c)=\bigcap_{\mG\in \Vm_c}\tilde \mI(\mG)
$$
and $\mA(c)=\pi(\tilde \mA(c))$.
The Ma\~n\'e set is defined  
by
$$
\tilde \mN(c)=\bigcup_{\mG\in \Vm_c}\tilde \mI(\mG)
$$
and $\mN(c)=\pi(\tilde \mN(c))$.
A bigger set will be useful in some occasions,
defined by
$$
\tilde \mB(c)=\bigcup_{\mG\in \Om_c}\tilde \mI(\mG),
$$
where $\Om_c$ is as defined in  \ref{Phi}.
As an intersection of Lipschitz graphs,  the Aubry set $\tilde \mA(c)$
is a Lipschitz graph over $\mA(c)$.
Note however that the Ma\~n\'e set is not a Graph in general.
The sets 
$$
\tilde \mM(c)\subset\tilde \mA(c)\subset \tilde \mN(c)\subset \tilde \mB(c)
$$
are compact and invariant under $\phi$.
The compactness of  $\tilde \mN(c)$ and  $\tilde \mB(c)$
is mentioned here for completeness, it will be proved later
in this section, in  \ref{maneset} and \ref{Bset} below.
These Lemma also prove that the Ma\~n\'e set is indeed
the set of orbits called $c$-minimizing by Mather
and semi-static by Ma\~n\'e, and that the set $\tilde \mB$
is the set of minimizing orbits, called $\tilde \mG$
in \cite{Fourier}.

\subsubsection{}\label{wedge}
It is  possible to associate to each dual fixed point 
$\breve \mG \in \breve \Vm$ the  invariant set 
$$
\tilde \mI(\breve \mG)=\bigcap_{n\in \Nm}\phi^n
\Big(\overline{ \breve \mG}\Big)
$$ 
and its projection $\ \mI(\breve \mG)$ on $M$.
The following is  due to Fathi, \cite{Fathi2}.
\vs\\
\textsc{Proposition. }
\begin{itshape}
Let us fix a cohomology $c$, and 
consider pseudographs $\mG\in \Vm_c$ and
 $\breve \mG\in \breve \Vm_c$.
The set $\mG \tilde \wedge \breve \mG$ is non-empty, compact
and invariant by $\phi$. In addition, this set intersects the Aubry set 
$\tilde \mA(c)$, and satisfies
$$
\mG \tilde \wedge \breve \mG
\subset \tilde \mI(\mG)\cap  \tilde \mI(\breve \mG)
$$
so that
$$
\mG  \wedge \breve \mG
\subset  \mI(\mG)\cap  \mI(\breve \mG).
$$
Furthermore, for each pseudograph $\mG\in \Vm_c$, there 
exists a pseudograph $\breve \mG\in \breve \Vm_c$
such that 
$$
\mG \wedge \breve \mG
= \mI(\mG)= \mI(\breve \mG).
$$
In a symmetric way, for each pseudograph  $\breve \mG\in \breve \Vm_c$,
there exists a pseudograph  $\mG\in \Vm_c$ such that this relation holds.
As a consequence, we have 
$$
\tilde \mA(c)=\bigcap_{\mG\in \Vm_c}\tilde \mI(\mG)
=\bigcap_{\breve \mG\in \breve \Vm_c}\tilde \mI(\breve \mG)
$$
and
$$
\tilde \mN(c)=\bigcup_{\mG\in \Vm_c}\tilde \mI(\mG)
=\bigcup_{\breve \mG\in \breve \Vm_c}\tilde \mI(\breve \mG).
$$
\end{itshape}
\proof
We have already proved that the set
 $\mG\tilde \wedge \breve \mG$
is compact and not empty, see \ref{intersection}.
Let us prove that it is invariant.
In order to do so, we
consider a weak KAM solution $u$ and a dual weak
KAM solution $\breve u$ such that 
$\mG=\mG_{c,u}$ and $\breve \mG=\mG_{c,\breve u}$.
Let 
$(x(t),p(t)):\Rm\lto T^*M$ be an orbit  of the Hamiltonian flow, such that 
$(x(0),p(0))\in 
\mG\tilde \wedge \breve \mG$.
Clearly, both $u$ and $\breve u$ are differentiable at 
$x(0)$, and $p(0)=c_{x(0)}+du_{x(0)}$.
For each $m\leq n$ in $\Nm$, we have 
$$
 u(x(n))=
\min _{x\in M}  u(x)-A_c(m,x,n,x(n))+(n-m)\alpha(c)
$$
$$
\leq u(x(m))+A_c(m,x(m),n,x(n))+(n-m)\alpha(c).
$$
On the other hand, we have $(x(0),p(0))\in  \breve \mG$
hence, in view of \ref{initial},
$$
\breve u(x(n))=\breve u(x(m))+A_c(m,x(m),n,x(n))+(n-m)\alpha(c).
$$
As a consequence, the sequence
$n\lmto (u-\breve u) (x(n))$
is non-increasing on $\Nm$.
Since its initial value $(u-\breve u)(x(0))$
has been chosen to be a minimum of the function
$u-\breve u$, the sequence must be constant,
so that $x(n)$ is a point of $\mG\wedge \breve \mG$
for each $n\geq 0$.
A symmetric argument shows that this is also true for $n\leq 0$.
In addition, we obtain that the inequality 
$u(x(n))\leq u(x(m))+A_c(m,x(m),n,x(n))+(n-m)\alpha(c)$
is in fact an equality for $0\leq m \leq n$.
Since this formula is true in view of \ref{initial}
for $m\leq n\leq 0$ in $\Zm$,
we obtain that, for all $m\leq n$ in $\Zm$,
$$
u(x(n))= u(x(m))+A_c(m,x(m),n,x(n))+(n-m)\alpha(c).
$$
In other words, the curve $x(t)$ is calibrated by $\mG$ and by $\breve \mG$,
see \ref{calibrated}.
This implies that $(x(n),p(n))\in \mG\cap \breve \mG$
for each $n\in \Zm$, and, since $x(n)\in \mG\wedge \breve \mG$,
we get $(x(n),p(n))\in \mG\tilde \wedge \breve \mG$.
This proves that $\mG \tilde \wedge  \breve\mG$
is invariant by $\phi$ and contained in $\mI(\mG)$
and in $\mI(\breve \mG)$.

Every compact invariant set of $\tilde \mI(\mG)$
carries an invariant measure.
As a consequence, every compact invariant set 
of $\tilde \mI(\mG)$ intersects the Mather set $\tilde \mM(c)$,
see \ref{ensembles}.
Since $\tilde \mM(c)\subset \tilde \mA(c)$, the set 
 $\mG \tilde \wedge  \breve\mG$, which is a compact and
invariant subset of $\tilde \mI(\mG)$, intersects 
$\tilde\mA(c)$.
Let us now fix the Pseudograph $\mG_{c,u}\in \Vm_c$,
and prove the existence of a pseudograph $\breve \mG\in \breve \Vm_c$
such that  
$
\mG \wedge \breve \mG
= \mI(\mG)= \mI(\breve \mG).
$
In order to do so, we set
$$
\breve u:= \lim_{n\lto \infty} \breve T_c ^nu-n\alpha(c)
=\lim _{n\lto \infty}\breve T_c ^n T_c ^n u.
$$ 
It follows from \ref{inequal} and \ref{alpha} that the limit exists and that
$\breve u\leq u$.
Let $(x(t),p(t)):\Rm\lto T^*M$ 
be a Hamiltonian orbit satisfying  $(x(0),p(0))\in \tilde \mI(\mG)$.
The orbit $x(t)$ is then calibrated by $\mG$, see \ref{calibrated},
so that the relation
$$
u(x(n))-u(x(m))=A_c(m,x(m);n,x(n))+(n-m)\alpha(c)
$$
holds for all $m\leq n$ in $\Zm$.
It is clear from this relation that,
for each $n\in \Nm$, 
$$
\breve T_c ^nu(x(0))
\geq u(x(n))-A_c(0,x(0);n,x(n))=u(x(0))+n\alpha(c),
$$
So that $\breve T_c ^nu(x(0))-n\alpha(c)=u(x(0))$,
hence $\breve u=u$ on $\mI(\mG)$.
As a consequence, the set of points minimizing 
$u-\breve u$ contains  $\mI(\mG)$.
Since we have already proved that this set is contained
in  $\mI(\mG)$, we can conclude, as desired, that
$$
\mG_{c,u}\wedge  \mG_{c, \breve u}
=\mI(\mG).
$$
Setting 
$u'=\lim T_c^n \breve u +n\alpha(c)=\lim T_c^n\breve T_c^n\breve u$,
the same proof gives that 
$$
\mG_{c,u'}\wedge  \mG_{c, \breve u}
=\mI(\mG_{c, \breve u}).
$$
We claim that $u'=u$, so that we have proved 
$$
\mG_{c,u}\wedge  \mG_{c, \breve u}
=\mI(\mG_{c,u})=\mI(\mG_{c, \breve u}).
$$
In order to prove that $u'= u$, 
we first recall that $\breve u\leq u$, so that 
$T_c^n \breve u+n\alpha(c)\leq T_c^n u+n\alpha(c)=u$,
and $u'\leq u$.
On the other hand, for each $\epsilon>0$, there exists $N\in \Nm$ such that 
$\breve T^N_cu-N\alpha(c)\leq \breve u+\epsilon$,
hence 
$$
u'\geq \lim_{n\lto \infty}T^n_c \breve T^N_cu+(n-N)\alpha(c)-\epsilon
\geq \lim_{n\lto \infty}T^{n-N}_c u+(n-N)\alpha(c)-\epsilon
=u-\epsilon.
$$
We have proved that $u'=u$.
\qed\vs\\
The pairs $(u,\breve u)$
of fixed points of $T_c+\alpha(c)$
and $\breve T_c-\alpha(c)$ which satisfy 
$$
\breve u=\lim \breve T_c^nu-n\alpha(c)\, \,;\,\,\,
 u=\lim T_c^n \breve u+n\alpha(c)
$$
are \textit{conjugate} in the sense of Fathi.

\subsubsection{}\label{semicontinu}
\textsc{Proposition. }
\begin{itshape}
The restriction to $\Vm$
of the function $c:\Pm\lto H^1(M,\Rm)$
is continuous and proper.
\end{itshape}
\vs\\
\proof
Let us consider a compact subset  $C$ of $H^1(M,\Rm)$.
The Family of Hamiltonians $H(t,x,c_x+p),c\in C$,
is a uniform family of Hamiltonians, see Appendix \ref{uniform}.
As a consequence, the associated functions
$A_c(0,.;1,.), c\in C$ form an equi-semi-concave family
of functions on $M\times M$.
As a consequence, the functions 
$A_c(0,x;1,.),c\in C, x\in M$ form an equi-semi-concave family of 
functions on $M$, see Appendix \ref{semiconcave}.
It follows that 
the functions $u(x)+A_c(0,x;1,.),c\in C, x\in M$ also form
an equi-semi-concave family,
hence that the functions 
$\min _x u(x)+A_c(0,x;1,.),c\in C$ 
 form
an equi-semi-concave family.
As a consequence, the set $\Phi(\Pm_C)$
is relatively compact.
Since the set $\Vm_C$
is obviously closed,
and contained in  $\Phi(\Pm_C)$, it is compact.
\qed
We have proved the following Lemma, which is interesting in itself:
\vs\\
\textsc{Lemma. }
\begin{itshape}
If $C$ is a compact subset of $H^1(M,\Rm)$,
the set $\Phi(\Pm_C)$ is equi-semi-concave.
\end{itshape}

\subsubsection{}\label{barrier}
Following Mather, we will use the function
$$
h_c(x,y):=\liminf_{n\lto \infty} A^n_c(x,y)+n\alpha(c).
$$
In view of \ref{liminf}, the function $h_c(x,.)$
is a fixed point of $T_c+\alpha(c)$.
Similarly, the function 
$-h_c(.,y)$ is a fixed point of $\breve T_c-\alpha(c)$.
Let us recall here some basic properties of the function $h_c$.
\begin{itemize}
\item
For each $x,y,z\in M$ and $c\in H^1(M,\Rm)$,
we have the triangle inequality
$
h_c(x,y)+h_c(y,z)\geq h_c(x,z).
$
\item
For each $x,y\in M$ and $c\in H^1(M,\Rm)$, we have $h_c(x,y)+h_c(y,x)\geq h_c(x,x)\geq 0$,
\item
For each compact set $C\subset H^1(M,\Rm)$, the set of functions
$h_c:M\times M\lto \Rm,c\in C,$ is equi-semi-concave.
\end{itemize}

\subsubsection{}\label{uh}
\textsc{Proposition. }
\begin{itshape}
If the pseudograph $\mG_{c,u}$ is a fixed point of $\Phi$, then we have 
$$u(y)-u(x)\leq h_c(x,y)$$
for each $x$ and $y$. In addition, 
$$
u(x)=\min _{y\in M}u(y)+h_c(y,x)=
\min _{a\in \mA(c)} u(a)+h_c(a,x).
$$
\end{itshape}

\proof
We have, for each $n$,
$u= T_c^nu+n\alpha(c)$.
As a consequence, for each $n$,
$$
u(x)=\min _{y\in M}\big(u(y)+A_c(0,y;n,x)+n\alpha(c)\big).
$$
We obtain the inequality
$u(x)\leq u(y)+A_c(0,y;n,x)+n\alpha(c)$
and, by taking the liminf, 
$u(x)\leq u(y)+h_c(y,x)$.
In order to obtain the first equality,
we 
consider a point $y_n\in M$ such that
$$
u(x)=u(y_n)+A_c(0,y_n;n,x)+n\alpha(c).
$$
We  consider an increasing sequence $n_k$ of integers
such that the subsequence  $y_{n_k}$ has a limit $y$, 
and refine this subsequence in such a way  that the subsequence
$A_c(0,y,n_k,x)$ has a limit $l$.
We have 
$$
u(x)= 
u(y)+l\geq u(y)+h_c(y,x).
$$
Cumulated with the previously shown inequality,
 this proves the first equality
in the statement.
In order to prove the second equality, notice that the set of points 
$y$ which minimize the function $u(.)+h_c(.,x)$ is precisely
the set $\mG\wedge \mG_{c,-h_c(.,x)}$, and that 
$\mG_{c,-h_c(.,x)}\in \breve \Vm_c$.
By \ref{wedge}, the set  $\mG\wedge \mG_{c,-h_c(.,x)}$
intersects $\mA(c)$.
\qed

\subsubsection{}
\textsc{Corollary. }
\begin{itshape}
For each $x$ and $y$ in $M$ and $c\in H^1(M,\Rm)$, we have 
\end{itshape}
$$
h_c(x,y)=
\min _{z\in M}
h_c(x,z)+h_c(z,y)=
\min _{a\in \mA(c)}
h_c(x,a)+h_c(a,y).
$$
The following result connects our definition of
the Aubry set to the one of Mather.

\subsubsection{}
\textsc{Proposition. }
\begin{itshape}
The Aubry set $\mA(c)$ is the set of points $x$ such that
$h_c(x,x)=0$.
\end{itshape}
\vs\\
\proof
Let us consider a Hamiltonian trajectory 
$(x(t),p(t)):\Rm\lto T^*M$ such that  $(x(0),p(0))\in \tilde \mA(c)$.
This trajectory is calibrated by each fixed point
of $T_c+\alpha(c)$, so  in particular by  $h_c(x(0),.)$.
Consequently, we have
$$
h_c(x(0), x(n))
-h_c(x(0), x(0))
=A_c(0,x(0);n,x(n))+n\alpha(c).
$$
Taking a subsequence such that $x(n)$ has a limit
$x$, and then a subsequence such that 
$A_c(0,x(0);n,x)+n\alpha(c)$ is converging to a limit 
$l\geq h_c(x(0),x)$, we get, at the limit,
$$
h_c(x(0), x)
-h_c(x(0), x(0))
\geq h_c(x(0),x)
$$
thus  $h_c(x(0),x(0))\leq 0$
and then  $h_c(x(0),x(0))= 0$.
We have proved that the function $h_c(x,x)$ vanishes on $\mA(c)$.

Conversely, let us assume that $h_c(x,x)=0$.
Then there exists an increasing sequence $n_k$ 
of integers and a sequence of trajectories
$(x_k(t),p_k(t)):[0,n_k]\lto T^*M$
such that 
$x(0)=x(n_k)=x$ and 
$$
\int _0^{n_k}L(t,x_k(t),\dot x_k(t))
-c_{x_k(t)}(\dot x_k(t))+\alpha(c) \,dt
=A_c(0,x;n_k,x)+n_k\alpha(c)
\lto 0.
$$
Let $y_k:[-n_k,n_k]\lto M$
be the curve such that 
$y_k(t)=x_k(t+n_k)$ for $-n_k\leq t\leq 0$ 
and 
$y_k(t)=x_k(t)$ for $t\geq 0$. 
The sequence $y_k$ is $C^2$-bounded hence,
by taking a subsequence, we can suppose that 
$y_k$ is converging with its derivative, uniformly on compact sets,
to a limit $y(t):\Rm\lto M$.
We claim that this limit $y$ is calibrated by each 
fixed points of $T_c+\alpha(c)$.
Let $u$ be a such a fixed point.
We have, for each $n\in \Nm$ and $k$ large enough, 
$$
0\geq u(y_k(n))-u(y_k(-n))-A_c(-n,y_k(-n);n,y_k(n))-2n\alpha(c)
$$
$$
=u(y_k(n_k))-u(y_k(-n_k))-A_c(-n_k,y_k(-n_k);n_k,y_k(n))-2n_k\alpha(c)
$$
$$
-\big(
u(y_k(n_k))-u(y_k(n))-A_c(n,y_k(n);n_k,y_k(n))-(n_k-n)\alpha(c)
\big)
$$
$$
-\big(
u(y_k(-n))-u(y_k(-n_k))-A_c(-n_k
,y_k(-n_k);-n,y_k(-n))-(n_k-n)\alpha(c)
\big)
$$
$$
\geq 
u(y_k(n_k))-u(y_k(-n_k))-A_c(-n_k,y_k(-n_k);n_k,y_k(n))-2n_k\alpha(c)
$$
$$
=-A_c(-n_k,x;n_k,x)-2n_k\alpha(c)=-2A_c(0,x;n_k,x)-2n_k\alpha(c).
$$
For fixed $n$, we now take the limit $k\lto \infty$, and get 
that 
$$
u(y(n))-u(y(-n))=A_c(-n,y(-n);n,y(n))+2n\alpha(c).
$$
Consequently, the curve $y(t)$ is calibrated by $u$.
Since this holds for each weak KAM solution $u$,
we have $x=y(0)\in \mA(c)$.
\qed
\vs\\

The following well-known result connects our definition of 
the Ma\~n\'e set with the usual one,
and implies its compactness.

\subsubsection{}\label{maneset}
\textsc{Lemma}
\begin{itshape}
The following properties 
are equivalent for a continuous  curve 
 $P(t)=(x(t),p(t)):\Rm\lto T^*M$.
\begin{itemize}
\item[($i$)]
The curve $P(t)$ is a Hamiltonian trajectory and 
$P(\Zm)\subset \tilde \mN(c)$.
\item[($ii$)]
The curve $P(t)$ satisfies 
$p(t)=\partial_vL(t,x(t),\dot x(t))$
and there exists  $\mG_{c,u}\in \Vm_c$ such that, 
for  each $m\geq n$ in $ \Zm$, we have 
$$
u(x(m))-u(x(n))=
\int_n^m L(t,x(t),\dot x(t))-
c_{x(t)}(\dot x(t))dt +(m-n)\alpha(c).
$$
\item[($iii$)]
The curve $P(t)$ satisfies 
$p(t)=\partial_vL(t,x(t),\dot x(t))$
and 
for  each $m\geq n$ in $ \Zm$, we have 
$$
\int_n^m L(t,x(t),\dot x(t))-
c_{x(t)}(\dot x(t))dt+(m-n)\alpha(c)
=\min_{l \in \Nm, l>0}
A_c(0,x(n);l,x(m))+l\alpha(c).
$$
\end{itemize}
\end{itshape}
\proof
We shall prove that $(iii)\Rightarrow(ii)$.
The other implications are left to the reader.
Let $P(t)$ be a curve satisfying $(iii)$.
let $n_k$ be an increasing sequence of 
integers such that $x(-n_k)$ has a limit $\alpha$.
Then we have, for $m\geq n$, 
$$
\int_n^m L(t,x(t),\dot x(t))-
c_{x(t)}(\dot x(t))dt+(m-n)\alpha(c)
$$
$$
=\int_{-n_k}^m L(t,x(t),\dot x(t))-
c_{x(t)}(\dot x(t)) +\alpha(c)dt
-
\int_{-n_k}^n L(t,x(t),\dot x(t))-
c_{x(t)}(\dot x(t))+\alpha(c)\,dt
$$
$$
=A_c(-n_k,x(-n_k);m,x(m))+(m+n_k)\alpha(c)
-A_c(-n_k,x(-n_k);n,x(n))+(n+n_k)\alpha(c).
$$
By $(iii)$, we have 
$$
A_c(-n_k,x(-n_k);m,x(m))+(m+n_k)\alpha(c)
=\min_{l\in \Nm, l>0} 
A_c(0,x(-n_k);l,x(m))+l\alpha(c)
$$
$$\leq h_c(x(-n_k),x(m))
$$
which implies that 
$$
A_c(-n_k,x(-n_k);m,x(m))+(m+n_k)\alpha(c)\lto h_c(\alpha,x(m))
$$
as $k\lto \infty$.
Similarly, 
$$
A_c(-n_k,x(-n_k);n,x(n))+(n+n_k)\alpha(c)\lto h_c(\alpha,x(n)),
$$
so that 
$$
h_c(\alpha,x(m))- h_c(\alpha,x(n))
=\int_n^m L(t,x(t),\dot x(t))-
c_{x(t)}(\dot x(t))dt +(m-n)\alpha(c).
$$
We have proved $(ii)$ with $u=h_c(\alpha,.)$.
\qed
\vs\\
We now give equivalent definitions for the set $\tilde \mB$.
The following Lemma shows that the set $\tilde \mB$ 
is the set called $\tilde \mG$ in \cite{Fourier},
and implies its compactness.

\subsubsection{}\label{Bset}
\textsc{Lemma}
\begin{itshape}
The following properties 
are equivalent for a continuous  curve 
 $P(t)=(x(t),p(t)):\Rm\lto T^*M$.\begin{itemize}
\item[($i$)]
The curve $P(t)$ is a Hamiltonian trajectory and 
$P(\Zm)\subset \tilde \mB(c)$ .
\item[($ii$)]
The curve $P(t)$ satisfies 
$p(t)=\partial_vL(t,x(t),\dot x(t))$,
and there exists  a sequence 
$u_n$ of functions such that, 
for each  $m\geq n$, we have 
$T_c^{m-n}u_n=u_m$
and 
$$
u_m(x(m))-u_n(x(n))=
\int_n^m L(t,x(t),\dot x(t))-
c_{x(t)}(\dot x(t))dt.
$$
\item[($iii$)]
The curve $P(t)$ satisfies 
$p(t)=\partial_vL(t,x(t),\dot x(t))$,
and 
for each  $m\geq n$, we have 
$$
\int_n^m L(t,x(t),\dot x(t))-
c_{x(t)}(\dot x(t))dt
=A_c(n,x(n);m,x(m))
$$
\end{itemize}
\end{itshape}
\proof
$(ii)\Rightarrow (i)$.
Then for each pair $m\geq n$ of integers,
the curve $x(t):[n,m]\lto M$ is minimizing the action between its
 endpoints.
Hence the curve $P(t)$ is a Hamiltonian trajectory.
It follows from \ref{initial}
that, for each $n\geq 0$, we have 
$P(n)\in \overline{\mG_{c,u_n}}$
and since $P(n)=\phi^n(P(0))$, we have 
$$
P(0)\in \phi^{-n}(\overline{\Phi^n(\mG_{c,u_0}))}.
$$
This inclusion holds for all $n$, 
so that $P(0)\in \tilde \mI(\mG_{c,u_0})$.
Now $(i)$ follows from the fact
$\mG_{c,u_0}\in \Om$ and 
 that $\tilde \mI(\mG_{c,u_0})$ is invariant under $\phi$.

$(i)\Rightarrow (ii)$.
There exists a pseudograph 
$\mG_{c,u_0}\in \Om$ such that
$P(0)\in \tilde \mI(\mG_{c,u_0})$.
There exists a sequence $u_n, n\in \Zm$ of functions 
on $M$ such that $T_c^{m-n}u_n=u_m$ for $m\geq n$.
For each $m\geq 0$, since $P(m)\in  \overline{\mG_{c,u_m}}$,
we have 
$$
u_m(x(m))-u_0(x(0))=\int_0^m L(t,x(t),\dot x(t))-
c_{x(t)}(\dot x(t))dt.
$$
On the other hand, for each $n\leq 0$, we can find by minimization
a curve $y_n(t):[n,0]\lto M$ such that $y_n(0)=x(0)$ and 
$$
u_0(y_n(0))-u_n(y_n(n))=\int_n^0 L(t,y_n(t),\dot y_n(t))-
c_{y_n(t)}(\dot y_n(t))dt.
$$
There exists a subsequence $n_k$ such that the curves 
$y_{n_k}(t)$ converge, uniformly on compact sets, 
to a limit
$y(t):(-\infty,0]\lto M$.
By \ref{mini}, this curve satisfies, for all  $n\leq 0$,
$$
u_0(y(0))-u_n(y(n))=\int_n^0 L(t,y(t),\dot y(t))-
c_{y(t)}(\dot y(t))dt.
$$
Hence, for $n\leq 0\leq m$, we have 
$$
u_m(x(m))-u_n(y(n))=\int_n^0 L(t,y(t),\dot y(t))-
c_{y(t)}(\dot y(t))dt
+\int_0^m L(t,x(t),\dot x(t))-
c_{x(t)}(\dot x(t))dt.
$$
As a consequence, the curve obtained by gluing $y$ on $\Rm^-$
and $x$ on $\Rm^+$ is the projection of a Hamiltonian trajectory,
which, by Cauchy-Lipschitz uniqueness, has to be $P(t)$.
In other words, we have proved that $x(t)=y(t)$ on $\Rm^-$.
The relation of calibration is now established.

$(iii)\Rightarrow (ii)$.
Let $P(t)$ be a curve satisfying $(iii)$.
Then we have, for $m\geq n\geq k$, 
$$
\int_n^m L(t,x(t),\dot x(t))-
c_{x(t)}(\dot x(t))dt+(m-n)\alpha(c)
$$
$$
=A_c(k,x(k);m,x(m))+(m-k)\alpha(c)
-A_c(k,x(k);n,x(n))+(n-k)\alpha(c).
$$
Let us denote by $u_k^n$ the function
$$
u^n_k(x)=A_c(k,x(k);n,x)+(n-k)\alpha(c),
$$
we obviously have $T_c^{m-n}u_k^n+(m-n)\alpha(c)=u_k^m$
for $m\geq n\geq k$.
By diagonal extraction, 
we find  an increasing sequence of integers $n_k$
such that $u_{-n_k}^{n}$ has a limit $u_n$ for each $n$
as $k\lto \infty$. 
We then have $T_c^{m-n}u_n+(m-n)\alpha(c)=u_m$
for each $m\geq n$, so that $\mG_{c,u_n}\in \Om$.
In addition, we have 
$$
\int_n^m L(t,x(t),\dot x(t))-
c_{x(t)}(\dot x(t))dt+(m-n)\alpha(c)
=u_m(x(m))-u_n(x(n)).
$$

$(ii)\Rightarrow (iii)$ is clear.
\qed

\subsubsection{}
\textsc{Lemma}
\begin{itshape}
For each $P\in\tilde  \mB(c)$, the orbit $\phi^n(P)$
is $\alpha$-asymptotic and $\omega$-asymptotic to 
the Aubry set 
$\tilde \mA(c)$. As a consequence, the Mather set $\tilde \mM(c)$ 
is the closure of the union of the supports of the  invariant measures of 
the action of $\phi$ on $\tilde \mB(c)$
\end{itshape}
\vs\\
\proof
Let $P(t)=(x(t), p(t))$ be the Hamiltonian orbit of $P$.
Let $u_n, n\in \Zm$ be a sequence of continuous functions
such that, for $m\geq n$, we have 
 $u_m=T_c^{m-n}u_n+(m-n)\alpha(c)$
and 
$$
u_m(x(m))-u_n(x(n))=\int _n ^m 
L(t,x(t), \dot x(t))-c_{x(t)}(\dot x(t))
dt+(m-n)\alpha(c).
$$ 
The sequence $u_m, m\in \Zm$ is equi-semi-concave, hence equi-Lipschitz.
Together with \ref{alpha}, this implies 
that this sequence is equi-bounded.
If  $v$ is  a weak KAM solution, that is a fixed point of 
$T_c-\alpha(c)$, we have, for $m\geq n$,
$$
v(x(m))-v(x(n))\leq \int _n ^m 
L(t,x(t), \dot x(t))-c_{x(t)}(\dot x(t))dt
+(m-n)\alpha(c)
$$ 
It follows that the sequence $n\lmto u_n(x(n))-v(x(n))$
is non-decreasing and bounded, 
so that 
  it has a limit $l$ as $n\lto -\infty$.
Let us now consider an increasing sequence $n_k$ of integers
such that the sequence
$P(t-n_k)$ converges, uniformly on compact
sets, to a limit $Z(t)=(y(t),z(t))$ which is a Hamiltonian 
trajectory.
Extracting a subsequence in $n_k$, we can suppose that 
there exists a sequence $w_n$ of continuous functions on $M$
such that $u_{n-n_k}\lto w_n$ uniformly, for each $n$,
as $k\lto \infty$.
Then, the sequence $w_n$ satisfies $T_c^{m-n}w_n=w_m$
for $m\geq n$.
In addition, we have 
$
w_n(y(n))-v(y(n))=l
$
and, for $m\geq n$, 
$$
w_m(y(m))-w_n(y(n))=\int _n ^m 
L(t,y(t), \dot y(t))-c_{y(t)}(\dot y(t))
dt+(m-n)\alpha(c).
$$ 
It follows that, for $m\geq n$,  
$$
v(x(m))-v(x(n))=\int _n ^m 
L(t,y(t), \dot y(t))-c_{y(t)}(\dot y(t))
dt+(m-n)\alpha(c)
$$ 
which implies that $Z(\Zm)\in \tilde \mI(\mG_{c,v})$.
Since this holds for all weak KAM solutions $v$,
we obtain that $Z(\Zm)\in \tilde \mA(c)$.
We have proved that the trajectory
$P(t)$ is $\alpha$-asymptotic to $\tilde \mA(c)$.
Similarly, one can prove that it is also 
 $\omega$-asymptotic to $\tilde \mA(c)$.
This implies that the  invariant measures of the action 
of $\phi$ on $\tilde \mB(c)$ are supported
on $\tilde \mA(c)$.
\qed

%
%

\subsection{Static classes and heteroclinics}\label{sectionstatic}
In this section, we see that there is a natural partition
of the Aubry set in compact invariant subsets, which we 
call static classes, following the terminology of Ma\~n\'e.
This partition was first considered by Mather in \cite{MatherFourier}.
We also discuss the existence of heteroclinic orbits between these
static classes, extending to the non-autonomous case results
of Fathi, Contreras and Paternain, see \cite{Fathi3} and \cite{CP}.
This survey is also an occasion to introduce several technical lemmas
and notations to be used later.

\subsubsection{}
\textsc{Lemma. }
\begin{itshape}
Let $x$ and $y$ be two points in $M$.
The following properties are equivalent:
\begin{enumerate}
\item[(i)]
The points $x$ and $y$ are in $\mA(c)$ and 
 the function $z\lmto h_c(x,z)-h_c(y,z)$ is constant on $\mA(c)$.
\item[(ii)]
$h_c(x,y)+h_c(y,x)=0$.
\item[(iii)]
The points $x$ and $y$ are in $\mA(c)$ and, for each pair  $(\mG,\breve \mG)\in \Vm_c\times \breve \Vm_c$,
either the set $\mG\wedge \breve \mG$ contains $x$ and $y$
or it contains neither $x$ nor $y$.
\end{enumerate}
If $x$ and $y$ satisfy these properties, we have, for each $z\in M$, 
$
h_c(x,z)= h_c(x,y)+h_c(y,z).$
\end{itshape}
\vs\\
\proof
$i \Rightarrow ii.$
Assuming $i$, we evaluate the constant function at $x$ and $y$
and get $h_c(x,x)-h_c(y,x)=h_c(x,y)-h_c(y,y)$,
hence $h_c(x,y)+h_c(y,x)=0$.

$ii \Rightarrow iii$.
We have $h_c(x,x)\leq h_c(x,y)+h_c(y,x)=0$ hence 
$x\in \mA(c)$. Similarly, $y\in \mA(c)$.
 Let us consider 
$\mG=\mG_{c,u}\in \Vm_c$ and 
$\breve \mG=\mG_{c,\breve u}\in \breve \Vm_c$
such that $x\in \mG\wedge \breve \mG$
(such a pair exists because $x\in \mA(c) $).
We have to prove that $y \in \mG\wedge \breve \mG$.
We have the inequalities
$u(y)\leq u(x)+h_c(x,y)$ and
$\breve u(y)\geq \breve u(x)-h_c(y,x)$.
We obtain the inequality
$$
(u-\breve u)(y)\leq (u-\breve u)(x)+h_c(x,y)+h_c(y,x).
$$
As a consequence, if 
$h_c(x,y)+h_c(y,x)=0$, then $y$ is also a point of minimum
of $u-\breve u$, which is the desired result.

$iii\Rightarrow ii$.
The point $x$ is a point of minimum of the function
$h_c(x,.)+h_c(.,x)$.
As a consequence, the point $y$ is also a point of minimum 
for this function, so that 
$h_c(x,y)+h_c(y,x)=h_c(x,x)+h_c(x,x)=0$.

$ii\Rightarrow i$.
We have the inequalities
$$
h_c(x,z)\leq h_c(x,y)+h_c(y,z)
\text{ and }
h_c(y,z)\leq h_c(y,x)+h_c(x,z).
$$
If $h_c(x,y)+h_c(y,x)=0$, then these inequalities sum to an equality,
hence they are both equalities. 
\qed

\subsubsection{}
The equivalent properties of the Lemma define an equivalence relation on 
$\mA(c)$. We call \textit{static classes} the classes of equivalence
for this relation. In other words, we say that the points $x$ and $y$
of $\mA(c)$ belong to the same static class if they satisfy $(i)$, $(ii)$
or $(iii)$ of the lemma.
We usually denote by $\mS$ a static class, and by 
$\mS(x)$
the static class containing $x$.
If $\mS$ is a static class, we denote by $\tilde \mS$
the set 
of points of $\tilde \mA(c)$ whose projection on $M$
belong to $\mS$. We will also call static classes
the sets of the form
$\tilde \mS$.
The static classes $\mS$ are compact and partition $\mA(c)$,
the static classes $\tilde \mS$ are compact, invariant,
 and they partition 
$\tilde \mA(c)$.
The invariance is a direct consequence of the caracterisation $(iii)$
of the equivalence relation.
To each point $x$ of the Aubry set $\mA(c)$, we associate 
the weak KAM solution $h_c(x,.)$, and we denote by  $E_{c,x}$
the associated element of  $\Vm_c$.
The pseudographs  of this form are called elementary solutions 
of $\Vm_c$. 
Two points of a same static class
give rise to the same elementary solution,
we will denote by $E_{c,\mS}$ the elementary solution
induced by points of $\mS$.
There is a one to one correspondence between 
the set of static classes and the set of elementary solutions.
We will denote by $\Em_c\subset \Vm_c$ this set.
We endow it with the induced metric,
it is clearly a compact set for this metric.
We also denote by $\breve E_{c,\mS}$
the fixed point of $\breve \Phi$ associated to the dual
weak KAM solution $-h(.,x)$ for $x\in \mS$. 

%
%
%
%
\subsubsection{}\label{limit}
\textsc{Proposition. }
\begin{itshape}
Let $\mG\in \Vm_c$ be a fixed point, and let
$P$  be a point of $\bar \mG$.
The $\alpha$-limit of the orbit $\phi^n(P)$ is contained
in one static class $\tilde \mS \subset \tilde \mA(c)$.
We also have $P\in \bar E_{c,\mS}$.
In a similar way, if 
$P\in \breve \mG\in \breve \Vm_c,
$ 
then the  $\omega$-limit 
of the orbit of $P$ is contained in one  static class
of $\tilde \mA(c)$.  
\end{itshape}
\vs\\
\proof
Let $\alpha\subset M$ be the projection of the 
 $\alpha$-limit of the orbit of 
$P\in \mG_{c,u}$.
Note that this $\alpha$-limit is contained in $\tilde \mI(\mG)$,
so that it is a Lipschitz graph above $\alpha$.
We claim that, for each weak KAM solution 
or backward weak KAM solution 
$v$,
the function $u-v$ is constant on $\alpha$.
Clearly, this implies that $\alpha$ is contained in
a static class.
In order to prove the claim,
we consider the projection  $x(t)$ on $M$ of the orbit of $P$.
The curve $x(t)$ is calibrated by $u$
on $(-\infty,0]$, hence the equality 
$$
u(x(-m))-u(x(-n))=A_c(-n,x(-n),-m, x(-m))
$$
holds for all $n,m\in \Nm$ such that $m\leq n$.
On the other hand, 
if $v$ is a weak KAM solution
or a backward weak KAM solution, 
we have the inequality
$$
v(x(-m))-v(x(-n))\leq A_c(-n,x(-n),-m, x(-m))
$$
for all $n,m\in \Nm$ such that $m\leq n$.
We deduce that the sequence
$(u-v)(x(-n)):\Nm \lto \Rm$ is non-increasing.
Now let  $y=\lim _{k\rightarrow \infty} x(-n_k)$
and  $z=\lim _{k\rightarrow \infty} x(-m_k)$
be  two points of $\alpha$. We can suppose that
$n_k\leq m_k\leq n_{k+1}$ by extracting subsequences.
We obtain 
$(u-v)(x(-n_k))\geq (u-v)(x(-{m_k}))\geq (u-v)(x(-n_{k+1}))$,
and at the limit
$(u-v)(y)\geq (u-v)(z)\geq (u-v)(y)$.
Hence the function $u-v$ is constant on $\alpha$.
This proves the first statement of the proposition.

Taking $v=h_c(\alpha,.)$, we obtain that 
 the sequence
$u(x(-n))-h_c(\alpha,x(-n)):\Nm \lto \Rm$ is non-increasing.
On the other hand, we have 
$u(x(0))-h_c(\alpha,x(0))\leq u(\alpha).$
It follows that the sequence is in fact constant,
so that the curve $x(t)$ is calibrated by 
$E_{c,\mS(\alpha)}$ on $(-\infty,0]$, and, by \ref{calibrated},
$((x(0),p(0))\in \bar E_{c,\mS(\alpha)}$.
\qed
\vs\\
\textsc{Corollary. }
\begin{itshape}
Let $P\in \tilde \mN(c)$ be a point whose $\alpha$-limit
is contained in $\tilde \mS$ and whose $\omega$-limit is contained
in $\tilde\mS'$. We have
$$
P\in E_{c,\mS}\tilde \wedge \breve E_{c,\mS'}.
$$
\end{itshape}

\proof
Let $(x(t),p(t))$ be the orbit of $P$.
Let $\alpha$ be an $\alpha$-limit point
of the curve $x(t)$, and let $\omega$ be an $\omega$-limit point.
It follows from the proposition that $(x(m),p(m))\in  \bar E_{c,\mS}$
for each $m\in \Zm$.
Applying the discussions in the proof
of the proposition with 
the point $P=(x(m),p(m))$
and the functions $u=h_c(\alpha,.)$ 
and $v=-h_c(.,\omega)$,
we get that the sequence
$h_c(\alpha,x(n))+h_c(x(n),\omega)$
is non-decreasing on $n\leq m$. Since we can take any $m\in \Zm$,
this sequence is non-decreasing on $\Zm$.
Taking a sequence $m_k\lto \infty$ such that $x(m_k)\lto \omega$
we obtain the inequality
$$
h_c(\alpha,\omega)
\leq  h_c(\alpha,x(n))+h_c(x(n),\omega)
\leq \lim _k (h_c(\alpha,x(m_k))+h_c(x(m_k),\omega))
=h_c(\alpha,\omega).
$$
It follows that, for each $n$, 
$$
 h_c(\alpha,x(n))+h_c(x(n),\omega)=
h_c(\alpha,\omega)
=\min _M h_c(\alpha,.)+h_c(.,\omega).
$$
This is precisely saying that 
$$
x(n)\in E_{c,\mS(\alpha)}\wedge \breve E_{c,\mS(\omega)}.
$$
Recalling that $(x(n),p(n))\in  \bar E_{c,\mS}$, we obtain
$
(x(n),p(n))\in E_{c,\mS(\alpha)}\tilde \wedge \breve E_{c,\mS(\omega)}.
$
\qed

%
%
%
\subsubsection{}\label{voisinage}
\textsc{Lemma. }
\begin{itshape}
If the static class $\mS$ 
is isolated in $\mA(c)$, then
there exists a neighborhood $V$ of $\mS$
in $M$ such that the $\alpha$-limit of every point
$P\in E_{c,\mS}$ satisfying  $\pi(P)\in V$
is contained in $\mS$.
\vs\\\end{itshape}
\proof
If the result did not hold true, we could find a sequence
$(x_n,p_n)\in E_{c,\mS}$ that has a limit 
$(x,p)\in \tilde \mS$ and a sequence $\alpha_n$ of $\alpha$-limit
points of $(x_n,p_n)$ that has a limit $\alpha$ in 
$\tilde \mA(c)-\tilde \mS$.
Note that the  orbit 
$
(x_n(t),p_n(t)):(-\infty,0]\lto T^*M
$
of $(x_n,p_n)$ is 
contained in $ E_{c,\mS}$.
Hence it is calibrated by the function 
$h_c(x,.)$, that is 
$$
h_c(x,x_n(0))=h_c(x,x_n(-k))+A_c(-k,x_n(-k);0,x_n(0))+k\alpha(c)
$$
for all $k\in \Nm$.
At the  liminf $k\lto \infty$, for fixed $n$, we obtain the inequality
$
h_c(x,x_n)\geq h_c(x,\alpha_n)+h_c(\alpha_n,x_n)
$
hence the equality
$
h_c(x,x_n)=h_c(x,\alpha_n)+h_c(\alpha_n,x_n).
$
At the limit $n\lto \infty$ we get
$
0=h_c(x,x)=h_c(x,\alpha)+h_c(\alpha,x).
$
This is in contradiction with the fact that $\alpha$ and $x$
do not belong to the same static class.
\qed

\subsubsection{}\label{chains}
Let $\tilde \mS$ and $\tilde \mS'$ be two different static classes
in $\tilde \mA(c)$.
The set $E_{c,\mS}\tilde \wedge \breve E_{c,\mS'}$ contains 
$\tilde \mS$ and $\tilde \mS'$ 
as well as other orbits of 
$\tilde \mN(c)$.
The following result is similar to Theorem A of \cite{CP}.
\vs\\
\textsc{Proposition. }
\begin{itshape}
The set 
$
E_{c,\mS}\wedge \breve E_{c,\mS'}-(\mS\cup \mS')
$
is not empty and contains points in every neighborhood of $\mS$,
as well as in every neighborhood of $\mS'$.
More precisely,
if $\mS$ and $\mS'$ are two, possibly equal, static classes,
and
if  $\tilde \mK\subset \tilde \mS$ 
and 
 $\tilde \mK'\subset \tilde \mS'$ 
are two disjoint  compact invariant sets,
then the set 
$
E_{c,\mS}\wedge \breve E_{c,\mS'}-(\mK\cup \mK')
$
contains points in every neighborhood of $\mK$,
as well as in every neighborhood of $\mK'$.
\end{itshape}
\vs\\
\proof
Let $V$ be an open neighborhood of $\mK$ in $M$ which 
does not intersect $\mS'$.
Let us fix a recurrent 
orbit $(y(t),z(t)):\Rm\lto T^*M$ such that
$(y(0),z(0))=(y,z)\in \tilde \mK$ and a 
recurrent orbit $(y'(t),z'(t)):\Rm\lto T^*M$ such that 
$(y'(0),z'(0))=(y',z')\in \tilde \mK'$.
Consider  a sequence $n_k\in \Nm$ of integers  and
a sequence $(x_k(t),p_k(t)):[0,n_k] \lto T^*M$
of Hamiltonian trajectories such that
$x_k(0)=y$ and $x_k(n_k)=y'$
and 
$$
\int_0^{n_k}L(t,x_k(t),\dot x_k(t)) 
-c_{x_k(t)}(\dot x_k(t))+\alpha(c)\;dt
= A_c(0,y;n_k,y')+n_k\alpha(c)
\lto h_c(y,y').
$$
We extend the curve $x_k:[0,n_k]\lto M$
to a curve $x_k:\Rm \lto \Rm$
by setting  $x_k(t)=y(t)$ for $t\leq 0$ and 
$x_k(t)=y'(t-n_k)$ for $t\geq n_k$. 
Let $a_k$ and $b_k$ be two increasing sequences of integers
such that 
$y(-a_k)\lto y$ and $y'(b_k)\lto y'$.
The existence of such sequences follows from the fact that the curves
$y(t)$ and $y'(t)$ are recurrent.
Since the curve $y(t)$ is calibrated by $h_c(y,.)$,
we have, as $k\lto \infty$,
$$
A_c(-a_k,y(-a_k);0,y)+a_k\alpha(c)=-h_c(y,y(-a_k))\lto 0
$$
and similarly
$$
A_c(0,y';b_k,y_k(b_k))+b_k\alpha(c)=h_c(y',y(b_k))\lto 0.
$$
As a consequence, we have, as $k\lto \infty$, 
$$
A_c(-a_k,x_k(-a_k);b_k+n_k,x_k(b_k+n_k))+(b_k+a_k+n_k)\alpha(c)
\lto h_c(y,y').
$$
For each $k$, let $T_k$ 
be the maximum of all times $i\in \Nm$ such that
$x_k(i)\in V$.
Note that $x_k(T_k+1)$ does not belong to $V$.
We can assume, taking a subsequence, that
the curve 
$x_k(t+T_k)$ is converging,
uniformly on compact sets
to a limit $x(t):\Rm \lto M$.
Let us now  fix $m\leq n$ in $\Zm$.
Summing the inequalities 
$$
\liminf _{k\lto \infty}
\Big(A_c(-a_k,x_k(-a_k);T_k+m,x_k(T_k+m))+(T_k+m+a_k)\alpha(c)\Big)
\geq  h_c(y,x(m)),
$$
$$
\liminf  _{k\lto \infty}
 \Big(A_c(T_k+m,x_k(T_k+m);T_k+n,x_k(T_k+n))\Big)
=
A_c(m,x(m);n,x(n))
$$
and
$$
\liminf _{k\lto \infty}
 \Big(A_c(T_k+n,x_k(T_k+n);b_k+n_k,x_k(b_k+n_k)))+(b_k+n_k-T_k-n)\alpha(c)\Big)
\geq  h_c(x(n),y'),
$$
we get 
$$
h(y,y')=
\liminf  _{k\lto \infty}
A_c(-a_k,x_k(-a_k);b_k+n_k,x_k(b_k+n_k))+(b_k+a_k+n_k)\alpha(c)
$$
$$
\geq  h_c(y,x(m))+ A_c(m,x(m);n,x(n))+
(n-m)\alpha(c)+h_c(x(n),y').
$$
Since the converse inequality obviously
holds, we obtain the equality
$$
h_c(y,y')= h_c(y,x(m))+ A_c(m,x(m);n,x(n))
+(m-n)\alpha(c)+h_c(x(n),y'),
$$
for all $m\leq n$.
It follows that all the inequalities above are in fact equalities, so that 
we also have
$$
h_c(y,x(m))+A_c(m,x(m);n,x(n))+(m-n)\alpha(c)
=h_c(y,x(n))
$$
so that the orbit $x(t)$ is calibrated by the weak KAM solution
$h_c(y,.)$ on $\Rm$.
Hence it is the projection of a Hamiltonian trajectory
$(x(t),p(t))$.
Moreover, we have the equality
$$
h_c(y,y')= h_c(y,x(n))+h_c(x(n),y'),
$$
so that  the point $x(n)$ is a point of minimum of the function
$h_c(y,.)+h_c(.,y')$. Hence
it belongs to 
 $E_{c,\mS}\wedge \breve E_{c,\mS'}$.
We have proved that the 
sequence $(x(n),p(n)),n \in \Zm$ is an orbit of $\phi$
which is contained in the invariant graph 
 $E_{c,\mS}\tilde \wedge \breve E_{c,\mS'}$.
Since the point $x(1)$ is not a point of $\mK$, this orbit
does not intersect the invariant set  $\tilde \mK$.
As a consequence,  the point $x(0)$ belongs to $\bar V-\mK$.
We have proved that the set 
$
E_{c,\mS}\wedge \breve E_{c,\mS'}-(\mK\cup \mK')
$
contains points in each neighborhood of $\mK$.
One can prove in a similar way that this set contains 
points in every neighborhood of $\mK'$.
\qed
%
%
%
%
\subsubsection{}\label{connected}
\textsc{Corollary. }
\begin{itshape}
A static class $\tilde \mS$ cannot be decomposed as the union of 
two disjoint invariant compact subsets.
\end{itshape}
\vs\\
\proof
Assume, by contradiction,  that there exists a static class 
$\tilde \mS=\tilde \mK_1 \cup \tilde \mK_2$,
with $\tilde \mK_i$ invariant, compact and disjoint.
In view of \ref{chains}, the set 
$$
E_{c,\mS}\wedge \breve E_{c,\mS}-
(\mK_1\cup \mK_2)
$$
is not empty.
This is a contradiction since
$
E_{c,\mS}\wedge \breve E_{c,\mS}=\mS
$
and 
$
\mK_1\cup \mK_2=\mS.
$
\qed
%
%
\subsubsection{}
Let $(x(t),p(t)):\Rm\lto T^*M$ be an orbit of the Ma\~n\'e set,
that is an orbit satisfying $(x(0),p(0))\in \tilde \mN(c)$.
This orbit is $\alpha$-asymptotic to a static class 
$\tilde \mS$, and $\omega$-asymptotic to a static class 
$\tilde \mS'$.
\vs\\
\textsc{Lemma. }
\begin{itshape}
The inclusion $(x(0),p(0))\in \tilde \mA(c)$
holds if and only if $\mS= \mS'$. In this case, we have
 $(x(0),p(0))\in \tilde \mS$
\end{itshape}
\vs\\
\proof
Let us first assume that $\mS=\mS'$.
In this case, we see from Corollary \ref{limit} that
$
(x(0),p(0))\in E_{c,\mS} \tilde \wedge \breve E_{c,\mS}.
$
But is is clear from the definition of static classes
that 
$E_{c,\mS} \tilde \wedge \breve E_{c,\mS}=\tilde \mS$.
Consequently, we have $(x(0),p(0))\in \tilde \mS\subset \tilde \mA(c)$.
Conversely, assume that $(x(0),p(0))\in  \tilde \mA(c)$.
Then this point is contained in one static class $\tilde \mS_0$.
Since this static class is compact and invariant, it contains
the $\alpha$ and the $\omega$-limits of the orbit 
$(x(t),p(t))$, so that we have
$\tilde \mS=\tilde\mS_0=\tilde\mS'$.
\qed
\\
\textsc{Corollary }
\begin{itshape}
We have the equality $\tilde \mA(c)= \tilde \mN(c)$
if and only if there is exactly one static class in $\tilde \mA(c)$.
\end{itshape}

\subsubsection{}
Let $\tilde \mH_c(\tilde \mS,\tilde \mS')$ 
be the set of orbits of  $\tilde \mN(c)$
which are heteroclinic orbits between 
the static classes $\tilde \mS$ and $\tilde \mS'$, we denote by 
$ \mH_c( \mS,\mS')$  its projection on $M$.
We have
$$
\tilde \mN (c)=\tilde \mA(c)\cup \bigcup_{\mS,\mS'}
\tilde \mH_c(\mS,\mS'),
$$
where the union is taken on all pairs $(\mS,\mS')$
of different static classes.
Recall, from Corollary \ref{limit}, that
$$
\tilde \mH_c(\tilde \mS,\tilde \mS') 
\subset 
E_{c,\mS}\tilde \wedge \breve E_{c,\mS'}.
$$
The following result is from 
\cite{Fathi3} and  \cite{CP}.

%
%
%
%
%
%
\subsubsection{}\label{heteroclinics}
\textsc{proposition. }
\begin{itshape}
If the static class $\tilde \mS$ 
is properly contained and isolated in $\tilde \mA(c)$, then
there exists an orbit of $\phi$ in 
$\tilde \mN(c)-\tilde \mA(c)$ which is $\alpha$-asymptotic
to $\tilde \mS$. 
This orbit is then $\omega$-asymptotic to another static class 
$\tilde \mS'$.
\end{itshape}
\vs\\
\proof
Let us chose, according to \ref{voisinage}, a neighborhood $V$
of $\mS$ such that every orbit of $E_{c,\mS}$ starting above $V$
has its $\alpha$-limit contained in $\mS$.
Now let us choose any static class $\mS''$ different from $\mS$.
In view of \ref{chains}, the set 
$E_{c,\mS} \wedge \breve E_{c,\mS''}$
intersects $V-\mS$.
Let $P(t)=(x(t),p(t)):\Rm\lto T^*M$
be an orbit such that 
$P(0)\in E_{c,\mS}\tilde \wedge \breve E_{c,\mS''}$
and 
  $x(0)\in V-\mS$.
 The $\alpha$-limit of the orbit $P(t)$ is contained
in $\tilde \mS$.
On the other hand, this orbits belongs 
to $\tilde \mN(c)$, hence its $\omega$-limit
is contained in some static class $\tilde \mS'$.
\qed
%
%
%
%
%
%
\subsubsection{}
We have treated so far the case where there exist several static classes.
We recall, however that the existence of a single static class
in $\mA(c)$, is, for $c$ fixed, a generic property of the Lagrangian,
see \cite{CP}.
We will explain in section \ref{covering} a device due
 to Fathi, as well as Contreras and Paternain, see 
\cite{Fathi3} and  \cite{CP}, which allows to treat this case.

%
%
%
%
%
%
%
%
%
%
%
%
%
%
%
%
%
%
%
\section{Abstract mechanisms}

%
\subsection{The relation and its dynamical consequences}\label{relation}
\label{dynamics}
We define the forcing relation $\ffleche$, and describe its dynamical consequences.
We prove Proposition \ref{orbites}.

\subsubsection{}
Let us introduce some useful notations.
Given two subsets $\mG$ and $\mG'$ of $T^*M$,
we define the relation 
$\mG\fleche_N \mG'$ as follows:
$$
\mG\fleche_N \mG'
\Longleftrightarrow
\bar \mG' \subset \bigcup _{n=1} ^{N} \phi^n(\mG),
$$
where as usual $\bar \mG$ is the closure of $\mG$.
We say that $\mG$ forces $\mG'$, and  write $\mG\fleche \mG'$ if there exists an integer 
$N$ such that $\mG \fleche_N \mG'$.
If $\mG$ is a subset of  $T^*M$
and if $c\in H^1(M,\Rm)$,
the relations 
$$
\mG \fleche c \text{ and } \mG \fleche_N c$$
mean that there exists 
an overlapping pseudograph $\mG'$ 
of cohomology $c$ and such that $\mG\fleche \mG'$
(resp.  $\mG\fleche_N \mG'$).
To finish, for $c$ and $c'$ two cohomology classes, the relation
$$
c\fleche_N c'
$$
means that, 
for each pseudograph $\mG\in \Pm_c$,
we have $\mG\fleche_N c'$.
As the reader may have guessed, we will then 
say that $c$ forces $c'$
($c\fleche c'$) if there exists an integer $N$ such that 
$
c\fleche_N c'.
$
The relation $\fleche$ (between subsets as well as between 
cohomology classes) is obviously transitive.
We will be concerned in this paper with understanding
the relation $\fleche$ between cohomology classes.
For this purpose, it is useful to introduce the symmetric relation
$$
c\ffleche c' \Longleftrightarrow 
c\fleche c' \text{ and } c'\fleche c.
$$
We say that $c$ and $c'$ force each other if $c\ffleche c'$.
\vs\\
\textsc{Proposition }
\begin{itshape}
The forcing relation 
$
\ffleche 
$
is an equivalence relation on 
$H^1(M,\Rm)$.\\
\end{itshape}
Note that we have $c\fleche_1 c$ for each $c$
since 
$\overline{ \Phi(\mG)}\subset \phi(\mG)$ for each $\mG\in \Pm_c$,
which can be written $\mG \fleche_1 \Phi(\mG)$.

\subsubsection{}\label{torus}
Let us   present a simple  (negative) result about this
relation.
If $\mG$ is a  the graph of a continuous section 
of $T^*M$, and is invariant under $\phi$, then
 $\mG\in \Vm\cap \breve \Vm$ is in fact a
Lipschitz graph,
and the relation $c(\mG)\fleche c$ holds if and only if
$c=c(\mG)$. 
Note that, if $C\subset H^1(M,\Rm)$ is bounded,
 it is possible to chose
a uniform constant $K$ such that all 
the invariant Lipschitz Graphs
$\mG$ 
whose  cohomology satisfies $c\in C$ are $K$-Lipschitz.
In other words, 
the elements of $\Vm_C\cap \breve \Vm$
are equi-Lipschitz graphs.
Of course, we would like to be able to prove that
the forcing relation $\ffleche$ has non-trivial classes.
We first restate and   prove Proposition \ref{orbites}.

\subsubsection{}
\textsc{Proposition  }
\begin{itshape}
\begin{itemize}
\item[$(i)$]
If the cohomology class $c$ forces the cohomology class 
$ c'$,
there exists a heteroclinic trajectory of the Hamiltonian flow
between $\tilde \mA(c)$ and $\tilde \mA(c')$.
For any closed forms $\eta$ of cohomology $c$ and $\eta'$
of cohomology $c'$, 
there exists a positive integer $N$ and a trajectory 
$(q(t),p(t)):[0,N]\lto T^*M$ of the Hamiltonian flow such that
$p(0)=\eta_{q(0)}$ and $p(N)=\eta'_{q(N)}$.
There exists a trajectory $(q(t),p(t)):[0,\infty)\lto T^*M$
which satisfies $p(0)=\eta_{q(0)}$ and
 is $\omega$-asymptotic to  $\tilde \mA(c')$. 
There exists  a trajectory $(q(t),p(t)):(-\infty,0]\lto T^*M$
which satisfies $p(0)=\eta'_{q(0)}$ and 
 is $\alpha$-asymptotic to  $\tilde \mA(c)$.
\item[$(ii)$]
Let $c_i,i\in \Zm,$ be a sequence of cohomology classes such
that $c_i$ forces $c_{i+1}$ for each $i\in \Zm$.
Fix for each $i$ a neighborhood $U_i$ of $\tilde \mM(c_i)$
in $T^*M$.
There exists a trajectory of the Hamiltonian flow which visits
in turn all the sets $U_i$. In addition, 
if the sequence stabilizes to $c-$ on the left, or to $c+$ on the right,
the trajectory can be assumed negatively asymptotic to   $\mA(c-)$
or positively asymptotic to  $\mA(c+)$.
\end{itemize}
\end{itshape}
\proof
Let us first assume that $c\fleche c'$.
Take a fixed point $\mG_c\in \Vm_c$.
There exists a graph $\mG\in \Pm_{c'}$
such that $\mG_c\fleche \mG$.
Now, consider a pseudograph  $\breve \mG_{c'}\in \breve\Vm_{c'}$.
It follows from Lemma \ref{intersection}
that $\mG$ intersects $\breve \mG_{c'}$.
In view of \ref{limit}, 
the points of intersection are $\alpha$-asymptotic
to $\tilde \mA(c)$ and $\omega$-asymptotic to $\tilde \mA(c')$.
In the same way, we can take for $\mG_c$ the graph of 
the closed form $\eta$, choose $\mG\in \Pm_{c'}$ such that 
 $\mG_c\fleche \mG$,
and take for $\breve \mG_{c'}$ the graph of $\eta'$.
The points of the intersection $\mG\cap  \breve \mG_{c'}$
have trajectories from $\mG_c$ to  $\breve \mG_{c'}$.
The other statements of  $(i)$ are proved similarly.

\subsubsection{}\label{proximite}
\textsc{Lemma }
\begin{itshape}
Let us fix a cohomology $c$.
\begin{itemize}
\item[$(i)$]
For each neighborhood $U$ of $\tilde \mB(c)$, there exists $N\in \Nm$
such that, for all  $l\geq N$ and  all $\mG \in \Pm_c$,  we have 
$$
\phi^{-l}\left(\overline{ \Phi^{2l}(\mG)}\right)\subset U.
$$
\item[$(ii)$]
If $V$ is an open neighborhood
of $\tilde \mM(c)$ in $T^*M$, 
there exists $N\in\Nm$ such that,
for each $\mG\in \Pm_c$ and 
each $P\in\overline{ \Phi^{N}( \mG)}$, one of the points 
$
\phi^{-i}(P), 1\leq i\leq N-1
$
belongs to $V$.
\end{itemize}
\end{itshape}
\proof
In order to prove $(i)$, it is sufficient 
to prove that, 
if $\mG_n\in \Pm_c$
is a sequence of pseudographs,
if $m_n$ is an increasing sequence of integers,
and if $(x_n(t), p_n(t)):[0, m_n]\lto T^*M$
is a Hamiltonian trajectory which satisfies 
$$
(x_n(m_n),p_n(m_n))\in \overline{\Phi^{2m_n}( \mG_n)}
$$
and which converges uniformly on compact sets to a limit
 $(x(t), p(t)):\Rm^+ \lto T^*M$,
then $(x(0),p(0))\in \tilde \mB(c)$.

Let us write the pseudographs $\Phi^{m_n}(\mG_n)$
on the form $\mG_{c,u_n}$.
For each $k,n\in \Nm$, we have 
$$
T^{k}_c u_n(x_n(k))
=
 u_n(x_n(0))+
\int_{0}^{k}L(t,x_n(t),\dot x_n(t)) 
+c_{x_n(t)}(\dot x_n(t))dt.
$$
Since the functions 
$
u_n
$
lie in the image of the operator $T^{m_n}_c$, they 
are  equi-Lipschitz, and there exists a real sequence $\lambda_n$
such that the sequence of functions 
$
\lambda_n+ u_n
$
has an  accumulation points in $C(M,\Rm)$.
As a consequence, we can assume, taking a subsequence
if necessary, that the functions 
$
\lambda_n+ u_n
$
converge uniformly to a limit $u$.
We have $\mG_{c,u}=\lim \Phi^{m_n}(\mG_n)\in \Om_c$.
For each fixed $k\in \Zm$, taking the limit as $n\lto \infty$, we get
$$
T^{k}_c u(x(k))
=
u(x(0))+
\int_{0}^{k}L(t,x(t),\dot x(t)) 
+c_{x(t)}(\dot x(t))dt.
$$
Hence we have 
$P(k-1)\in \Phi^{k-1}(\mG_{c,u})$, and therefore
$P(0)\in \phi^{1-k}(\Phi^{k-1}({\mG_{c,u})})$.
Since this holds for all $k\in \Nm$, 
we conclude 
that  
$$
P(0)\in \tilde \mI(\mG_{c,u})\subset \tilde \mB(c).
$$

In order to prove $(ii)$,
it is useful  to recall that $\tilde \mB(c)$ is a compact set,
invariant under the time-one flow $\phi$, and that the Mather set 
$\tilde \mM(c)$ is  the union of the supports
of the invariant measures of the action of $\phi$ on $\tilde \mB(c)$.
The claim below  follows from general  facts 
about  topological dynamics
on compact spaces:
For each neighborhood $W$ of $\tilde \mM(c)$ in  $\tilde \mB(c)$,
there exists an integer $k$ such that, for each point 
$P$ of  $\tilde \mB(c)$, one of the points 
$\phi^i(P), 1\leq i\leq k$, belongs to $W$.
As a consequence, if $V$ is a neighborhood of  $\tilde \mM(c)$
in $T^*M$, there exists a neighborhood $U$ of 
 $\tilde \mB(c)$ in $T^*M$ such that, for each $P\in U$,  one of the points 
$\phi^i(P), 1\leq i\leq k$, belongs to $V$.
Now let us take $l\geq k$ such that $(i)$ holds for this
neighborhood $U$, and set $N=2l$.
For each $\mG\in \Pm_c$ and each $P\in \overline{\Phi^N(\mG)}$,
we have $\phi^{-l}(P)\in U$.
Hence one of the points $\phi^{i-l}(P), 1\leq i\leq k$ is in $V$,
which proves $(ii)$.
\qed

\subsubsection{}
Let us now prove $(ii)$ of the proposition.
Let $M_i\in \Nm, i\in \Zm$
be a sequence of integers such that $c_i\fleche_{M_i}c_{i+1}$,
and let $W_i\subset  V_i$ be compact  neighborhoods of $\tilde \mM(c_i)$.
In view of lemma \ref{proximite},
there exists a sequence $N_i$ of integers such that,
for each 
$\mG\in \Pm_{c_i}$ and each
$$
P\in \phi^{-N_i}(\overline{\Phi^{N_i}(\mG)}),
$$
one of the points $\phi^l(P)$, $0\leq l \leq N_i$ 
belongs to $W_i$.

Let us first fix an integer $k\in \Nm$,
and choose a pseudograph $\mG^{k}_{-k}\in \Pm_{c_{-k}}$.
Since $c_{-k}\fleche _{M_{-k}} c_{1-k}$,
there exists a pseudograph $\mG^{k}_{1-k}\in \Pm_{c_{1-k}}$
such that 
$\Phi^{N_{-k}}(\mG^{k}_{-k})\fleche _{M_{-k}}\mG^{k}_{1-k}$.
We build, by induction, a sequence 
$\mG^{k}_i\in \Pm_{c_i}, i\geq -k$, of pseudographs such that
$$
\Phi^{N_i}(\mG^{k}_i)\fleche _{M_i}\mG^{k}_{i+1}
$$
for each $i\geq -k$.

Let us now take a point $P^k_k\in  \mG^{k}_k$.
There exists a positive integer $l^k_{k-1}\leq M_{k-1}$
such that 
$\phi^{-l^k_{k-1}}(P^k_k)\in \Phi^{N_{k-1}}(\mG^k_{k-1})$.
We then set $P^k_{k-1}=\phi^{-(l^k_{k-1}+N_{k-1})}(P^k_k)$,
this point belongs to $\mG_{k-1}^k$.
We can  build a sequence $P^k_i,-k\leq i\leq k$
of points of $\mG^k_i$
 and a sequence $l^k_i, -k\leq i\leq k-1$
of  integers satisfying 
$0\leq l^k_i\leq M_i$
such that 
$$
\phi^{N_i+l^k_i}(P^k_i)=P^k_{i+1}
$$
for each $i$.
In addition, one of the points $\phi^j(P^k_i), 0\leq j\leq N_i$
belongs to $W_i$.

There exists an increasing sequence $k_n$ of integers
such that each of the sequences $n\lmto l^{k_n}_i$, 
for fixed $i$, is the constant $l_i$ after a certain rank,
and each of the sequences $n\lmto P^{k_n}_i$, for fixed $i$,
is converging to $P_i$.
Clearly, we have 
$\phi^{l_i+N_i}(P_i)=P_{i+1}$ for each $i\in \Zm$,
and one of the points $\phi^j(P_i),0\leq j\leq  N_i$
belongs to $W_i$.
This proves the main part of the statement. 

If the sequence $c_i$ stabilizes to $c-$ on the right,
then it is possible to build a sequence $\mG_i\in \Pm_{c_i}$
as above 
which stabilizes to $\mG-\in \Vm_{c-}$ on the right,
and we obtain by the above method an orbit which is  
$\alpha$-asymptotic to $\tilde \mA(c-)$ and then visits
in turn all the sets $W_i$.
 If the sequence $c_i$ stabilizes to $c+$ on the right,
say for $i\geq I$, then it is possible to
impose that $P_I\in \breve \mG+ \in \breve \Vm_{c+}$
in the construction above, and we then obtain an orbit
which is $\omega$-asymptotic to $\tilde \mA(c+)$.
\qed

\subsection{Evolution operators}\label{operators}
We define operators on $\Pm$
that generalize the Lax-Oleinik operator $\Phi$.
These operators will play a central role in the proof of our main results.

\subsubsection{}\label{AN}
Given two integers $N'\geq N\geq 1$, and a cohomology $c$,
we define the function 
$A_c^{N,N'}:M\times M\lto \Rm$ 
by the expression 
$$A_c^{N,N'}(x,y)=\min_{k\in \Nm, N\leq k\leq N'}
A_c(0,x;k,y)+k\alpha(c).
$$
Since each of the mappings 
$$
H^1(M,\Rm)\lto C(M\times M,\Rm)
$$
$$
c\lmto A_c(0,.;k,.)
$$
is continuous (see appendix \ref{uniform}),
 it is easy to see that, for fixed $N'\geq N$, 
the mapping
$$
H^1(M,\Rm)\lto C(M\times M,\Rm)
$$
$$
c\lmto A^{N,N'}_c
$$
is continuous. 
\vs\\
\textsc{Proposition}
\begin{itshape}
Let $c$ be a fixed cohomology class.
For each $\epsilon>0$, there exist integers $N'\geq N\geq 1$
such that 
$$\|A^{N,N'}_c-h_c\|_{\infty} \leq \epsilon.$$
More precisely, there exists an integer $N_0$ and a function 
$\Lambda:\Nm\lto \Nm$ such that the conclusion holds
if $N\geq N_0$ and $N'\geq \Lambda(N)$.
\end{itshape}

\proof
The functions $A_c^{N,N'}$, $N\leq N'\in \Nm$, and the function
$h_c$ have a common modulus of continuity and a common uniform
bound.
As a consequence, the pointwise limit 
$$h_c(x,y)=
\lim_{N\lto \infty}
\lim_{N'\lto \infty}A_c^{N,N'}(x,y)
$$
is uniform.
\qed

\subsubsection{}\label{phiu}
It is useful to generalize the operators $\Phi^N:\Pm\lto\Pm$.
Given two integers $N'\geq N\geq 1$ and an open set
$U \subset M$, we define the operator
$$
\Phi^{N,N'}_U:\Pm \lto \Pm
$$
by the relation $\Phi^{N,N'}_U(\mG_{c,u})=\mG_{c,T_{c,U}^{N,N'}u}$
where 
$$
T_{c,U}^{N,N'}u(x)
=\min _{y\in \bar U,N\leq k\leq N'} T_c^ku(y)+k\alpha(c)
=\min _{y\in \bar U} u(y)+A_c^{N,N'}(y,x).
$$
For simplicity we will denote by $\Phi^{N,N'}$ the operator
$\Phi^{N,N'}_M$.
For each $\mG\in \Pm$, we  have
$$
\mG \fleche_{N'}\Phi^{N,N'}(\mG).
$$ 
\textsc{Lemma. }
\begin{itshape}
For each integers  $1\leq N\leq N'$ and each open set $U\subset M$,
the operator 
$
\Phi^{N,N'}_U:\Pm \lto \Pm
$
is continuous, when the source is endowed with the seminorm
$\|.\|_U$ and the image with the norm $\|.\|$, see \ref{norm}.
\end{itshape}
\vs\\
\proof
Let $\mG=\mG_{c,u}$ and $\mG_1=\mG_{c_1,u_1}$ be two pseudographs.
We have
$$
\|
\Phi_U^{N,N'}(\mG_1)-\Phi_U^{N,N'}(\mG)
\| \leq
|c_1-c|+\|T^{N,N'}_{c_1,U}u_1-T^{N,N'}_{c,U}u\|.
$$
In order to estimate the term 
$\|T^{N,N'}_{c_1,U}u_1-T^{N,N'}_{c,U}u\|$, 
let us write 
$$
T^{N,N'}_{c,U}u=u(y)+A^{N,N'}_c(y,x)
$$
with $y\in \bar U$.
Then, we have 
$$
T^{N,N'}_{c_1,U}u_1-T^{N,N'}_{c,U}u\leq 
u_1(y)+A^{N,N'}_{c_1}(y,x)-
u(y)-A^{N,N'}_c(y,x)
$$
and by symmetry
$$
\|T^{N,N'}_{c_1,U}u_1-T^{N,N'}_{c,U}u\|
\leq \sup_{y\in U}
|u_1(y)-u(y)|+\|A^{N,N'}_{c_1}-A^{N,N'}_c\|
$$
The conclusion  follows from the continuity 
of the mapping
$c\lmto A_c^{N,N'}$, see \ref{AN}
\qed

%
%
%
%
%
\subsubsection{}\label{convergence}
Similarly, we define the operator
$$
\Phi^{\infty}_U:\Pm \lto \Vm
$$
by the relation 
$\Phi^{\infty}_U(\mG_{c,u})=\mG_{c,T_{c,U}^{\infty}u}$
where 
$$
T_{c,U}^{\infty}u(x)=\min _{y\in \bar U} u(y)+h_c(y,x).
$$
In the autonomous case Fathi proved that the sequence $\Phi^n(\mG)$
is converging to a fixed point of $\Phi$ for each $\mG\in \Pm$.
Such a result would be very useful to us, but does not hold
in our non-autonomous setting. It is replaced by the following one.\vs\\
\textsc{Proposition}
\begin{itshape}
Let $c$ be a fixed cohomology class.
For each $\epsilon>0$ there exist integers $N\leq N'$ such that,
for each pseudograph $\mG=\mG_{c,u}\in \Pm_c$ 
and each open set $U \subset M$,
we have
$$
\| \Phi^{N,N'}_U(\mG)-\Phi^{\infty}_{U}(\mG)\|
\leq \epsilon.
$$
More precisely, there exists an integer $N_0$ and a function 
$\Lambda:\Nm\lto \Nm$ such that the conclusion holds
if $N\geq N_0$ and $N'\geq \Lambda(N)$.
\end{itshape}
\vs\\
\proof
It is not hard to see that, for each continuous function $u$,
$$
\|
T^{N,N'}_{c,U}u
-
T^{\infty}_{c,U}u\|\leq 
\|A^{N,N'}_c-h_c\|.
$$
The proposition follows from \ref{AN}.
\qed

\subsubsection{}\label{phiufleche}
\textsc{Proposition}
\begin{itshape}
Let $\mG_0=\mG_{c_0,u_0}\in \Pm$ be a pseudograph, and let $\epsilon>0$
be fixed.
Assume that there exists an  open set $U\subset M$ 
 and two  compact sets $\mK\subset U$  and $\mK_1\subset M$ 
such that, for each 
$x\in \mK_1$, the minimum in the expression
$T^{\infty}_{c_0,U}u_0(x)=\min_{y\in \bar U}u_0(y)+h_{c_0}(y,x)$
is never reached outside of  $\mK$.
Then there exists integers $N\leq N'$,  a positive number $\delta $ 
and an open neighborhood $U_1$ of $\mK_1$
such that,
for each pseudograph 
$\mG\in \Pm$ satisfying $\|\mG-\mG_0\|_{U}\leq \delta$ 
we have 
$$
\mG_{|U}\fleche _{N'} 
\Phi^{N,N'}_U(\mG)_{|U_1}.
$$
More precisely, there exists an integer $N_0$ and a function 
$\Lambda:\Nm\lto \Nm$ such that the conclusion holds
if $N\geq N_0$ and $N'\geq \Lambda(N)$.
\end{itshape}
\vs\\
\proof
Let us denote by $\partial U$ the boundary of $U$.
There exists a positive number $\epsilon$ 
and a neighborhood $U_1$ of $\mK_1$ such that,
for each $x\in \bar U_1$,
$$
\min_{y\in \partial  U}u_0(y)+h_{c_0}(y,x)\geq 
\min_{y\in  U}u_0(y)+h_{c_0}(y,x)+7\epsilon.
$$
In view of \ref{AN}, there 
exist integers $N$ and $N'$ such that 
$$
\|A_{c_0}^{N,N'}-h_{c_0}\|\leq \epsilon.
$$
For fixed $N$ and $N'$, the function 
$A^{N,N'}_c$ depends continuously on $c\in H^1(M,\Rm)$, 
see \ref{AN}.
As a consequence, if $c$ is sufficiently close to $c_0$, 
we have 
$$
\|A_{c_0}^{N,N'}-A_{c}^{N,N'}\|\leq \epsilon.
$$
For these values of $N$ and $N'$, 
if 
 $u\in C(M,\Rm)$ is such that 
$\sup_U |u-u_0|\leq \epsilon$, 
we have, for each $x\in M$ and $y\in \bar U$,
$$
|u_0(y)+h_{c_0}(y,x)-u(y)-A^{N,N'}_c(y,x)|\leq 3\epsilon.
$$
Hence we have the inequality
$$
\min_{y\in \partial  U}u(y)+A^{N,N'}_c(y,x)\geq 
\min_{y\in  U}u(y)+A^{N,N'}_c(y,x)+\epsilon.
$$
As a consequence,
if $\|\mG_{c,u}-\mG_0\|_U$
is sufficiently small, then  there exists a compact set $\mK'\subset U$
such that the minimum in the expression 
$$
T^{N,N'}_{c,U}u(x)=\min_{y\in \bar U} u(y)+A^{N,N'}_c(y,x)
$$
is reached in $\mK'$ for all $x\in \bar U_1$.
Now let us set 
$v=T^{N,N'}_{c,U}u$ and consider a point
$$
(x,p)\in \overline{\mG_{c,v|U_1}}.
$$
The point $(x,p)$ is the limit of a sequence 
$(x_n,p_n)\in \mG_{c,v|U_1}$.
In other words, the points $x_n\in U_1$ are points of differentiability of
$v$, and we have $dv_{x_n}+c_{x_n}=p_n$.
Let $y_n\in \mK'$ and $k_n\in \Nm, N \leq k_n \leq  N'$, 
satisfy
$$
v(x_n)=u(y_n)+A_c(0,y_n;k_n,x_n)+k_n\alpha(c).
$$ 
By extracting a subsequence, we can suppose that the
sequence $k_n$ is a constant $k$.
By arguments similar to those of \ref{initial},
recalling that the function $u$ is semi-concave,
we conclude that the function $u$ is differentiable
at $y_n$, and, setting $z_n=c_{y_n}+du_{y_n}$,
that $\phi^k(y_n,z_n)=(x_n,p_n)$.
By extracting another subsequence, we can suppose that the sequence $y_n$
has a limit $y\in \mK'$.
We then have 
$$
v(x)=u(y)+A_c(0,y;k,x)+k\alpha(c),
$$
so that the function 
$u$ is differentiable at
$y$. 
Since the function $u$ is semi-concave, we then have
$du_y=\lim du_{y_n}$, see Appendix \ref{closed}.
At the limit in  $\phi^k(y_n,z_n)=(x_n,p_n)$,
we get  $\phi^k(y,z)=(x,p)$, where $z:=dv_y+c_y$.
We have proved that 
$$
\overline{\mG_{c,v|U_1}}
\subset \bigcup_{k=N}^{N'}\phi^k(\mG_{c,u|U}).
$$
\qed

\subsubsection{}\label{reduction}
\textsc{Proposition }
\begin{itshape}
Let $c$ and $c'$ be two cohomology classes.
Assume that, for each weak KAM solution 
$\mG_0 \in \Vm_c$,
there exists a positive number $\epsilon>0 $ and
an integer $N$  with the following property:
For each pseudograph $\mG \in \Pm_c$ such that 
$\|\mG-\mG_0\|\leq \epsilon$, there exists 
 a pseudograph $\mG'\in \Pm_{c'}$
such that $\mG\fleche_N \mG'$.
Then $c\fleche c'$.
\end{itshape}

\proof
By compactness of $\Vm_c$, there exists a neighborhood $\Um$ of
$\Vm_c$ in $\Pm_c$ and an integer $N$ such that, for all
$\mG\in \Um$, we have $\mG \fleche_N c'$.
In view of Proposition \ref{convergence}, there exist
integers $k\leq k'$ such that 
$\Phi^{k,k'}(\mG)\in \Um$ for each $\mG\in \Pm_c$.
We obtain, for each $\mG\in \Pm_c$, the existence of a $\mG'\in \Pm_{c'}$
such that 
$$
\mG\fleche_{k'}
\Phi^{k,k'}(\mG)
\fleche _{N} \mG'
$$ 
so that  $\mG \fleche_{k'+N}\mG'$.
\qed

\subsection{Coverings}\label{covering}
As was noticed by  Fathi, as well as Contreras and Paternain, see 
\cite{Fathi3} and  \cite{CP}, it is useful to study the effect of 
taking finite Galois coverings.

\subsubsection{}
Let $P:M_0\lto M$ be a finite connected  covering, and 
$P^*:H^1(M,\Rm)\lto H^1(M_0,\Rm)$ the induced mapping.
Let us also denote by $L\circ TP:\Rm\times TM_0\lto \Rm$ 
the lifted Lagrangian 
$$
L\circ TP(t,x,v)=L(t,P(x),dP_x(v)),
$$
and  by $T^*P:T^*M_0\lto T^*M$ the covering 
$$(x,p)\lmto (P(x), p\circ dP_x^{-1}).$$
The lifted Hamiltonian $H\circ T^*P$ is in natural duality
with the Lagrangian $L\circ TP$.
As a consequence, the Hamiltonian flow associated to 
the Lagrangian $L\circ TP$ is the 
Hamiltonian flow of $H\circ T^*P$, which is the lifting of
the Hamiltonian flow of $H$.
Each overlapping pseudograph $\mG=\mG_{c,u}$ on $M$
lifts to a pseudograph
$$
P^*\mG:=T^*P^{-1}(\mG)=\mG_{P^*c,u\circ P}
$$
on $M_0$. Note that $c(P^*\mG)=P^*(c(\mG))$.
It is not hard to see that the
Aubry set $\tilde \mA_{L\circ TP}(P^*(c))$ associated to $L\circ TP$ 
on $M_0$ is precisely
$$
\tilde \mA_{L\circ TP}(P^*(c))=
T^*P^{-1}(\tilde \mA_L(c)),
$$
while we only have the inclusion
$$
\tilde \mN_{L\circ TP}(P^*(c))\supset
T^*P^{-1}(\tilde \mN_L(c)).
$$
Finally, if $ \mS_{L\circ TP}$ is a static class of 
 $\mA_{L\circ TP}(P^*(c))$ , then 
 $P( \mS_{L\circ TP})$ is a static class of $\mA_L(c)$.
Note however that the lifting 
$P^{-1}(\mS_L)$ of a static class of $\mA_L(c)$
can contain several static classes of  $\mA_{L\circ TP}(P^*(c))$.
This is illustrated by the following result which,
in conjunction with \ref{heteroclinics},
allows to prove the existence of heteroclinic
orbits in the case where there is only one static class,
 see 
\cite{Fathi3} and  \cite{CP}.
We need first another definition.
If $\tilde \mX\subset T^*M$ is an invariant set of
the time-one flow $\phi$, then we denote
by $s\tilde \mX\subset  T^*M\times\Tm$
the set $\cup _{t\in \Rm, x\in \tilde \mX}(\phi^t_0(x),t)$
and by $s\mX$ its projection on $M\times \Tm$.

\subsubsection{}\label{coveringexists}
\textsc{Proposition. }
\begin{itshape}
Assume that
the set $\mA(c)$ contains finitely many static classes, 
and that 
 there exists an open neighborhood $U\subset M\times\Tm$
of the compact set 
$
s\mA(c)
$ 
such that the
mapping $h:H_1(U,\Zm)\lto H_1(M,\Zm)$
is not surjective,
where $h$ is the composition  of  the mappings
$$
H_1(U,\Zm)\overset{i_*}{\lto} H_1(M\times \Tm,\Zm)
\overset{p_*}\lto H_1(M,\Zm)
$$
induced from the inclusion 
and the projection.
Then there  exists a finite connected Galois covering
$P:M_0\lto M$ with $k$ sheets, $k\geq 2$, such that,
for each static class $\tilde \mS$ of $\tilde \mA(c)$,
the lifting $T^*P^{-1}(\tilde \mS)$ is the union of exactly $k$ 
different static classes of $\tilde \mA_{L\circ TP}(P^*(c))$.
\end{itshape}
\vs\\
\proof
Let $N$ be the number of static classes in $\mA(c)$.
First, we claim that for each static class $\mS$, the set $s\mS$
is connected. This follows easily from \ref{connected}.
As a consequence, we can suppose that the neighborhood
$U$ is a union of finitely many connected open sets $U_i,1\leq i\leq N$, 
each of which 
contains exactly one of the sets $s\mS$.
Since the group $H^1(M,\Zm)$ is Abelian and  of finite type, 
and since the  mapping  $h:H_1(U,\Zm)\lto H_1(M,\Zm)$ is not surjective,
there exists an integer $k\geq 2$ and a  surjective morphism
$g:H^1(M,\Zm)\lto \Zm/k\Zm$
whose kernel contains the subgroup  $h(H_1(U,\Zm))$.
There is a connected Galois covering $P:M_0\lto M$ with 
$k$ sheets 
associated to this morphism.
This means that, if $\chi:\pi_1(M)\lto H_1(M,\Zm)$
is the Hurewitz map, 
then the image  $P_*(\pi_1(M_0))$ in $\pi_1(M)$
is precisely the kernel of $g\circ \chi$.
The following diagram commutes.
$$
\begin{CD}
@.  \pi_1(M_0\times\Tm)@>{p_*}>>\pi_1(M_0)\\
@. @V{(P\times Id)_*}VV @V{P_*}VV \\
\pi_1(U) @>{i_*}>>\pi_1(M\times\Tm)@>p_*>>\pi_1(M)\\
@V{\chi}VV @V{\chi}VV @V{\chi}VV\\
H_1(U,\Rm)
@>{i_*}>>
H_1(M\times\Tm,\Rm)
@>{p_*}>>
H_1(M,\Rm)
@>g>>\Zm/k\Zm
\end{CD}
$$
We claim that 
$$
i_*(\pi_1(U))\subset (P\times Id )_*(\pi_1(M_0\times\Tm)),
$$
which implies that the covering $P\times Id$ is trivial above $U$.
In order to prove the claim,
let us first notice that 
$g\circ \chi \circ p_* \circ i_*=g\circ  p_* \circ i_*\circ \chi:
\pi_1(U)\lto \Zm/k\Zm$
is the zero map.
This implies that the image of the map $i_*:\pi_1(U)\lto \pi_1(M\times \Tm)$
is contained in the kernel of $g\circ \chi \circ p_*$.
In order to finish the proof of the claim, we check that
the image of $(P\times Id)_*$ is precisely the kernel
of  $g\circ \chi \circ p_*$.
This follows from the fact $\ker(g\circ \chi)=\text{im}( P_*)$,
and from the fact that the upper square of the diagram can be 
identified with 
$$
\begin{CD}
\pi_1(M_0)\times \Zm @>>>\pi_1(M_0)\\
 @V{P_*\times Id}VV @V{P_*}VV \\
\pi_1(M)\times\Zm @>>>\pi_1(M)
\end{CD}
$$
where the horizontal arrows are projections on the first factor.
The claim is proved,
so that the covering $P\times Id$ is trivial above $U$.
It follows that each connected component $U_i$ of 
$U$ has $k$ disjoint connected  preimages $V_i^j\subset M_0\times \Tm$.
Now it is not hard to see that 
the static classes of  $\mA_{L\circ TP}(P^*(c))$
are precisely the intersections
$$ 
TP^{-1}(\mA(c))\cap  V_i^j=
\mA_{L\circ TP}(P^*(c))\cap V_i^j, 
1\leq i\leq N, 
1\leq j\leq k.
$$
\qed

\subsubsection{}\label{coveringrelation}
\textsc{Proposition. }
\begin{itshape}
Let  $P:M_0\lto M$ be a finite Galois covering.
Let $c$ and $c'$ be two cohomology classes in $H^1(M,\Rm)$.
If $P^*(c)\fleche _N P^*(c')$ for the forcing relation 
associated to the Lagrangian $L\circ TP$ on $M_0$, then 
$c\fleche_N c'$.
\vs\\
\end{itshape}
\proof
Let us consider a pseudograph $\mG \subset \Pm_c$.
 If $P^*(c)\fleche_N P^*(c')$ then there exists a pseudograph
$\mG'$ on $M_0$ of cohomology $P^*(c')$ and such that 
$P^*\mG\fleche_N \mG'$. 
Let $\Dm$ be the group of deck transformations of the covering $P$.
The elements of $\Dm$ are the diffeomorphisms $D$ of $M_0$
such that $P\circ D=P$. 
To each element $D$ of $\Dm$ we associate
the fibered diffeomorphism $T^*D$ of $T^*M$ defined by 
$$
T^*D(x,p)=(D(x), p\circ dD_x^{-1}).
$$
This diffeomorphism is a Deck transformation of the covering 
$T^*P$.
Let us prove that there exists a pseudograph $\mG''$ on $M_0$
which is invariant by deck transformations, 
which has  cohomology $P^*(c')$, and such that 
$P^*\mG\fleche_N \mG''$. 
Let $\eta$ be a form on $M$ with cohomology $c'$, and 
let $P^*\eta$ be its lifting to $M_0$.
We write  $\mG'$ on the form $\mG_{P^*\eta,u}$.
Since the flow of $H\circ T^*P$ commutes 
with Deck transformations, and since the pseudograph 
$P^*\mG$ is invariant by deck transformations, 
we have 
$$
P^*\mG=T^*D(P^*\mG)\fleche_N T^*D(\mG')
$$
for each deck transformation $D$. It is easy to check that
$T^*D(\mG')=\mG_{P^*\eta, u\circ D^{-1}}$.
Setting
$$v:= \min _{D\in \Dm}  u\circ D^{-1},$$
and $\mG''=\mG_{\eta,v}$, we have 
$\mG''\subset \cup _{D\in \Dm} T^*D(\mG')$
hence $P^*\mG\fleche_N \mG'' $, and we have the desired 
properties for $\mG''$.
Since $P$ is a Galois covering, functions 
on $M_0$ which are invariant by deck transformations 
are liftings of functions on $M$.
As a consequence, there exists a continuous function $w:M\lto\Rm$
such that $v=w\circ P$. 
Hence the pseudograph $\mG''$
is the lifting of the  pseudograph $\mG_{\eta,w}$ on $M$.
Since $P^*\mG\fleche_N \mG''=P^*\mG_{\eta,w}$, we have
$\mG\fleche_N \mG_{\eta,w}$.
We have asssociated, to each pseudograph
$\mG\in \Pm_c$, a pseudograph $\mG_{\eta,w}\in \Pm_{c'}$
such that $\mG\fleche_N \mG_{\eta,w}$. This proves that 
$c\fleche_N c'$
\qed
%
%
%
%
%
%
%
%
%
%
%
%
%
%
%
\subsection{Mather's mechanism}\label{mathersection}
We comment and prove Theorem \ref{mather}.
Let us first  discuss some properties  of the subspace $R(c)$
as defined in \ref{mather}.

\subsubsection{}
It is useful to consider the \v Cech cohomology $\check H(.,\Rm)$
with real coefficients.
Recall that $H(.,\Rm)$ is the De Rham cohomology.
We identify the De Rham cohomology $H^1(M,\Rm)$
with the \v Cech cohomology $ \check H^1(M,\Rm)$.
If $K$ is any subset of $M$, we denote by
$\check \imath ^*_{K}$ the mapping
$\check H^1(M,\Rm)\lto \check H^1(K,\Rm)$ induced from the inclusion
$i_K:K\lto M$.
Recalling that the subspace $R(\mG)$ was defined in \ref{mather},
we have:
\vs\\
\noindent
\textsc{Lemma. }
\begin{itshape}
For each $\mG\in \Vm_c$,
we have 
$
R(\mG)= \ker (\check \imath^*_{\mI(\mG)}).
$
\end{itshape}
\vs\\
\proof
Consider  an open neighborhood 
$V_0$ of $\mI(\mG)$ which is such that,
for each  open neighborhood $V\subset V_0$ of $\mI(\mG)$ in $M$,
 $R(\mG)$ is the set of 
cohomology classes of smooth closed one-forms which vanish on $V$.
Let us fix an open neighborhood $V\subset V_0$ of $\mI(\mG)$ in $M$.
Clearly, the cohomology classes of  smooth 
closed one-forms which vanish on $V_0$ belong to 
$\ker( i^* _{V})$
(where $i^* _{V}$ is the mapping associated to the inclusion 
in De Rham cohomology).
 We have proved the inclusion
$
R(\mG)\subset \ker( i^* _{V})$.
Conversely, let $\omega$ be a smooth closed one-form on $M$ whose De Rham
cohomology belongs to 
$\ker(i^* _{V})$.
This means that the restricted form
$\omega_{|V}$ is exact. Hence there exists a smooth function $f$ on 
$V$ such that $\omega_{|V}=df$.
There exists a smooth function $\tilde f:M\lto \Rm$ which is equal
to $f$ in a neighborhood $\tilde V\subset V$ of $\mI(\mG)$. 
 The form $\omega-d\tilde f$ is a smooth closed form on $M$,
cohomologous to $\omega$, and vanishing on $\tilde V$.
As a consequence, we have
$[\omega]=[\omega-d\tilde f]\in R(\mG)$.
We have proved that 
$$
 \ker( i^* _{V})= R(\mG)
$$
for all open neighborhood $V\subset V_0$ of  $\mI(\mG)$ in $M$.
The Lemma follows because
$\ker (\check \imath^*_{\mI(\mG)})$
is equal to 
$\ker( i^* _{V}) $ when $V$ is a sufficiently small open neighborhood 
 of  $\mI(\mG)$ in $M$.\qed

\subsubsection{}
In order to avoid confusion, we shall denote by 
$j_{\tilde \mK}:\tilde \mK\lto T^*M$ the inclusion of a subset $\tilde \mK$
of $T^*M$ into $T^*M$, and by 
$\check \jmath_{\tilde \mK}^*:\check H^1(T^*M,\Rm)\lto \check H^1(\tilde \mK,\Rm)$
the associated mapping in \v Cech cohomology.
We identifying the \v Cech cohomology $\check H^1(T^*M,\Rm)$
with the De Rham cohomology $ H^1(T^*M,\Rm)$, and the mapping
  $\check \pi^*: \check H^1(M,\Rm)\lto \check H^1(T^*M,\Rm)$
with 
$\pi^*: H^1(M,\Rm)\lto H^1(T^*M,\Rm)$.
\vs\\
\noindent
\textsc{Lemma. }
\begin{itshape}
We have 
$(\pi^*)^{-1}\big(
\ker(\check  \jmath_{\tilde \mN(c)}^*) 
\big)
\subset  R(c)$.
\end{itshape}
\vs\\
\proof
It is enough to prove that, for each $\mG\in \Vm_c$, we have
$
\ker(\check  \jmath_{\tilde \mN(c)}^*) 
\subset \pi^*( R(\mG))$.
For each $\mG\in \Vm_c$, we have $\tilde \mI(\mG)\subset \tilde \mN(c)$,
hence 
$
\ker(\check  \jmath_{\tilde \mN(c)}^*) 
\subset 
\ker(\check  \jmath_{\tilde \mI(\mG)}^*).
$
So it is enough to prove that the inclusion
$\ker(\check  \jmath_{\tilde \mI(\mG)}^*)\subset \pi^*( R(\mG))$
holds for each $\mG\in \Vm_c$.
Let us consider such a pseudograph $\mG$.
Since the projection $\pi:TM\lto M$ induces a homeomorphism
$\pi_{|\tilde \mI(\mG)}:\tilde \mI(\mG)\lto \mI(\mG)$,
the commutative diagram
$$
\begin{CD}
TM @<{j_{\tilde \mI(\mG)}}<<
\tilde \mI(\mG)\\
@V{\pi}VV  @V{\pi_{|\tilde \mI(\mG)}}VV\\
M @<{i_{ \mI(\mG)}}<< \mI(\mG)
\end{CD}
$$
gives rise in \v Cech cohomology to the diagram
$$
\begin{CD}
H^1(TM,\Rm)@= \check H^1(TM,\Rm)@>
{\check \jmath^*_{\tilde \mI(\mG)}}>> 
\check H^1(\tilde \mI(\mG),\Rm)\\
@A{\pi^*}AA @A{\check \pi^*}AA  @A{\check \pi^*_{|\tilde \mI(\mG)}}AA\\
  H^1(M,\Rm)@=\check  H^1(M,\Rm) @>{\check \imath^*_{\mI(\mG)}}>> 
\check H^1( \mI(\mG),\Rm)
\end{CD}
$$
where the vertical arrows are isomorphisms.
We conclude that 
$$
\ker (\check \jmath^*_{\tilde \mI(\mG)})=
\check \pi^*(\ker (\check \imath^*_{\mI(\mG)}))=
 \pi^*(R(\mG)).
$$
\qed

\subsubsection{}
\textsc{Lemma. }
\begin{itshape}
The space $R(c)\subset H^1(M,\Rm)$
depends semi-continuously on $c$ in the following sense:
For each $c_0\in H^1(M,\Rm)$, there exists a neighborhood
$V$ of $c_0$ in $H^1(M,\Rm)$ such that, for each $c\in V$,
we have $R(c_0)\subset R(c)$.
\end{itshape}
\vs\\
\proof
Let us fix a cohomology class $c_0$.
We claim that, for each $\mG_0\in \Vm_{c_0}$, there
exists $\epsilon_0>0$ such that each fixed pseudograph
$\mG$ of $\Vm$ which satisfies $\|\mG-\mG_0\|\leq \epsilon_0$
has to satisfy
$$
 R(\mG_0)\subset R(\mG).
$$
This claim follows from the existence of a neighborhood
$U$ of $\mI(\mG_0)$ such that 
$R(\mG_0)$ is the set of cohomology classes 
of smooth closed one-forms which vanish on $U$,
 and from the fact that the inclusion
$\mI(\mG)\subset U$ holds when $\mG$ is sufficiently close
to $\mG_0$.
Let $B_{\Pm}(\mG,\epsilon)$ denote the open ball of center 
$\mG$ and radius $\epsilon$ in $\Pm$.
The compact set $\Vm_{c_0}$ is covered by a finite family of balls 
$B_{\Pm}(\mG_i,\epsilon_i)$ 
such that $\mG_i\in \Vm_{c_0}$ and such that
$R(\mG_i)\subset R(\mG)$ for each 
$\mG \in B_{\Pm}(\mG_i,\epsilon_i)$.
Since the map $c$ is continuous and proper
on $\Vm$, see \ref{semicontinu}, there exists a $\delta>0$
such that, for $|c-c_0|\leq \delta$, we have
$\Vm_c\subset \cup_i B_{\Pm}(\mG_i,\epsilon_i)$.
As a consequence, if $\mG$ belongs to some
$\Vm_c$ with  $|c-c_0|\leq \delta$,
then there exists $i$ such that 
$R(c_0)\subset R(\mG_i)\subset R(\mG)$.
We conclude that $R(c_0)\subset R(c)$ for  $|c-c_0|\leq \delta$.
\qed
\vs\\
The following proposition is the main step in the proof
of Theorem  \ref{mather}.
We denote by $B_E(r)$ the open  ball of radius $r$ centered
at the origin in the normed vector space $E$.
\subsubsection{}\label{matherprop}
\textsc{Proposition. }
\begin{itshape}
For  each   $\mG_0\in \Vm_{c_0}$,
there exists  a positive number $\epsilon_0$, and  an integer $N$
such that the following holds:
For each pseudograph $\mG\in \Pm$ 
satisfying  $\|\mG-\mG_0\|<\epsilon_0$ and 
$c(\mG)-c_0\in R(c_0)$,
for each cohomology class  $c$
satisfying
$c-c_0\in B_{R(c_0)}(\epsilon_0)\subset R(c_0) $,
we have 
$$
\mG\fleche_N c.
$$
\end{itshape}
Note in this statement that
the cohomology class of the pseudograph $\mG$ is denoted
by $c(\mG)$, and that the symbol $c$ alone  denotes another 
cohomology class.

\subsubsection{}
\textsc{Proof of Theorem \ref{mather}. }
We assume the proposition.
For Each $\mG_0\in \Vm_{c_0}$, we consider the number $\epsilon_0$
given by the proposition, and the open
ball $B_{\Pm}(\mG_0,\epsilon_0)$
of center $\mG_0$ and radius $\epsilon_0$ in $\Pm$.
Since $\Vm_{c_0}$ is compact, it can be covered by a finite number 
of these balls, we denote $\mG_i$ and $\epsilon_i$
the associated centers and radii.
Since the function $c$ restricted to $\Vm$ is proper,  \ref{semicontinu},
there exists a positive number $\delta$ such that
$\Vm_c\subset \cup_i B_{\Pm}(\mG_i,\epsilon_i)$
when $|c-c_0|\leq \delta$.
Consider two cohomology classes  $c$ and $c'$ 
in $c_0+B_{R(c_0)}(\epsilon)$, with  
$\epsilon=\min \{\delta, \epsilon_i\}$.
It follows from \ref{reduction}
that 
$c$ forces $c'$.
The theorem clearly follows.
\mbox{}\hfill {\small \fbox{}} 
$ _{\textsc{theorem}} $
 {\small \fbox{}\\}

\subsubsection{}\label{local}
\textsc{Proof of the Proposition. }
Let us fix 
a $\mG_0\in \Vm_{c_0}$ and 
choose a neighborhood $U$
of $\mI(\mG_0)$ in such a way that
$R(\mG_0)$ is the set of cohomology classes of 
smooth closed one-forms vanishing on $U$.
\vs\\
\textsc{Lemma}
\begin{itshape}
There exist $\delta>0$ and $N'\leq N$  in $\Nm$ such that,
for all overlapping pseudographs $\mG$
satisfying $\|\mG-\mG_0\|\leq\delta$,
we have 
$$
\mG_{|U}\fleche_N
\Phi_U^{N',N}(\mG)
$$
\end{itshape}
\proof
Let us write the pseudograph 
$\mG_0$ on the form $\mG_{c_0,u_0}$.
We have seen  in \ref{uh} that 
$$
u_0(x)=\min _{y\in M}u_0(y)+h_{c_0}(y,x)=
\min _{y\in \mA(c_0)}u_0(y)+h_{c_0}(y,x).
$$
As a consequence, we have 
$T^{\infty}_U u_0=T^{\infty}_M u_0= u_0$,
and the minimum in the definition of $T^{\infty}_U u_0(x)$
is not reached outside of $\mI(\mG)$, which is a compact set contained
in $U$.
The lemma now follows from proposition \ref{phiufleche}.
\qed

\subsubsection{}\label{deform}
\textsc{Lemma. }
\begin{itshape}
Let us fix a $\delta>0$.
There exists $\epsilon_0>0$ such that, if we take :\\
One one hand 
a cohomology class $c$
satisfying
$c-c_0\in R(\mG_0)$ and $\|c-c_0\|\leq \epsilon_0$;\\
On the other hand 
a pseudograph $\mG\in \Pm$
satisfying $\|\mG-\mG_0\|\leq \epsilon_0$
and $c(\mG)\in c_0+R(\mG_0)$;\\
Then there exists a pseudograph $\mG'\in\Pm_{c}$
with the following properties:
 $\|\mG'-\mG_0\|\leq \delta$ and $\mG_{|U}=\mG'_{|U}$.\vs
\end{itshape}

\proof
Let us write $\mG_0=\mG_{\eta_0,u_0}$.
In view of the definition of $R(\mG_0)$,
it is possible to associate to each cohomology
class $d\in R(\mG_0)$ a closed one-form 
$\mu_d$ which is null on $U$. In addition, we can impose
that the correspondence 
$d\lmto \mu_d$
is linear.
Given a pseudograph 
$\mG\in \Pm$
and a cohomology $c$ satisfying the hypotheses of the Lemma,
we consider  the  pseudograph 
$$\mG'
=\mG+\mG_{\mu_{(c-c(\mG))},0}\in \Pm_c .
$$ 
It is clear that $\mG'_{|U}=\mG_{|U}$,
that $c(\mG')=c$,  and that 
 $\|\mG'-\mG_0\|\leq \delta$
if $\epsilon_0$ is small enough.
\qed

\subsubsection{}
We are now in a position to end the proof of the proposition.
Let us consider $\delta$ given by Lemma \ref{local},
and the associated $\epsilon_0$ as given by Lemma \ref{deform}.
If $\mG$ and $c$ satisfy the hypotheses of the Proposition with this
value of $\epsilon_0$, then, by Lemma \ref{deform},
there exists a pseudograph $\mG'$ such that 
$c(\mG')=c$ and 
 $\mG'_{|U}=\mG_{|U}$
and  $\|\mG'-\mG_0\|\leq \delta$.
In view of Lemma \ref{local}, 
we have
 $\mG_{|U}\fleche_N \Phi^{N',N}_U(\mG') $, so that
 $\mG\fleche_N \Phi^{N',N}_U(\mG')$.
\mbox{}\hfill {\small \fbox{}} 
$ _{\textsc{proposition}} $
 {\small \fbox{}\\}

%
%
%
%
%
%
%
\subsection{Arnold's mechanism for 
systems with finitely many static classes}\label{finitesection}
%
%
We prove and generalize Theorem \ref{weakfinite}.

\subsubsection{}
Let $\tilde \mH_c(\tilde \mS,\tilde \mS')$ 
be the set of orbits of  the Ma\~n\'e  set $\tilde \mN(c)$
which are heteroclinic orbits between 
the static classes $\tilde \mS$ and $\tilde \mS'$, we denote by 
$ \mH_c( \mS,\mS')$  its projection on $M$.
We have, from section \ref{sectionstatic},
$$
\tilde \mN (c)=\tilde \mA(c)\cup \bigcup_{\mS,\mS'}
\tilde \mH_c(\mS,\mS'),
$$
where the union is taken on all pairs $(\mS,\mS')$
of different static classes.
In addition, it is useful to recall that 
$$
\tilde \mH_c(\tilde \mS,\tilde \mS') 
\subset 
E_{c,\mS}\tilde \wedge \breve E_{c,\mS'}.
$$
 We say that the set  $\tilde \mH_c(\mS,\mS')$
is \textit{neat} if it admits a compact subset $\tilde \mK$ which
contains one and only one  point in each orbit of 
$\phi_{|\tilde \mH_c(\mS,\mS')}$ and whose
projection $\mK$ on $M$ is acyclic.
This means that $\mK$ has a neighborhood $U$  whose inclusion $i$
into $M$ induces the null map $i_*:H_1(U,\Rm)\lto H_1(M,\Rm)$.

\subsubsection{}\label{finite}
\textsc{Theorem }
\begin{itshape}
Let $c_0$ be a cohomology class such that
the number of static classes in $\mA(c_0)$ is finite and greater than one.
Assume in addition that all the sets  $\tilde \mH_{c_0}(\mS,\mS')$ are neat.
Then the class $c_0$ is in the interior of its class of 
\ffleche-equivalence.
\end{itshape}
\vs\\
Let us  gather some preliminary consequences of the hypotheses.
\subsubsection{}\label{isolation}
\textsc{Lemma}
\begin{itshape}
We assume the hypotheses of the theorem.
Let $\mS_0$ be a static class and $V_0$ be a neighborhood of $\mS_0$.
\begin{itemize}
\item[($i$)] There exists an open neighborhood $V$ 
of $\mS_0$, contained in $V_0$, such that 
the boundary of $V$ does not intersect $\mI(E_{c_0,\mS_0}).$
\item[($ii$)] There exists an acyclic open set $U\subset V_0-\mS_0$
and a static class $\mS_1$
such that the intersection $U\cap \mI(E_{c_0,\mS_0})$
is not empty, compact,  and contained in $\mH(\mS_0, \mS_1)$.
\end{itemize}
\end{itshape}
\proof
Let $V_0$ be a neighborhood of $\mS_0$,
sufficiently small for lemma
\ref{voisinage} to apply, so that we have
$$
V_0\cap \mI(E_{c_0,\mS_0})
=
\mS_0\cup \bigcup_{\mS \in \Em(c_0)-\mS_0} \big( \mH(\mS_0,\mS)\cap V_0
\big),
$$
where the union is taken on all static classes $\mS\neq \mS_0$.
We shall also assume that $\bar V_0 \cap \mA(c_0)=\mS_0$.

For each static class $\mS$, let us consider an acyclic
compact set 
 $ \tilde \mK(\mS_0, \mS)$ 
which contains one and only one point in every orbit of 
$ \tilde \mH(\mS_0, \mS)$, and denote  by 
$\mK(\mS_0, \mS)$ its projection on the base.
Clearly, the sets 
 $\tilde \mK(\mS_0, \mS),\mS \in \Em(c_0)-\mS_0, $ are pairwise disjoint
and they all belong to the  Lipschitz graph
$\tilde \mI(E_{c_0,\mS_0})$, so that their projections
 $\mK(\mS_0, \mS)$ on the base are also pairwise disjoint.
Let us consider a static class $\mS \neq \mS_0$.
For $n$ large enough, we have 
$\pi\circ \phi^{-n}(\tilde \mK(\mS_0,\mS))\subset V_0$.
In addition, since $\mK(\mS_0,\mS)$ is acyclic in $M$,
the compact  $\tilde \mK(\mS_0,\mS)$ is acyclic in $T^*M$.
As a consequence, the compact set $\phi^{-n}(\tilde \mK(\mS_0,\mS))$
is acyclic in $T^*M$ and contained in the Lipschitz graph
$\tilde \mI(E_{c_0,\mS_0})$, so that 
$\pi\circ \phi^{-n}(\tilde \mK(\mS_0,\mS))$
is acyclic in $M$.
Consequently, 
recalling that the number of static classes is finite,
there is no loss of generality in supposing that the
sets $\mK(\mS_0, \mS),\mS \in \Em(c_0)-\mS_0,$ are all contained in  $V_0$.

Let us prove that each of the sets $\tilde \mK(\mS_0,\mS)$ is  isolated
in $\tilde \mI(E_{c_0,\mS_0})$.
Let $F$ be a compact neighborhood of $\mS_0$ which does not intersect 
any of the finitely many compact 
sets $\mK(\mS_0,\mS),\mS \in \Em(c_0)-\mS_0 $.
Since the points of $\tilde \mK(\mS_0,\mS)$ are $\alpha$-asymptotic to 
$\tilde \mS_0$ and $\omega$-asymptotic to $\tilde \mS$,
and since there are finitely many static classes,
there exists an integer $N$ such that all the sets 
$\pi\circ \phi^{n}(\tilde \mK(\mS_0,\mS)), n\in \Nm,\mS \in \Em(c_0)-\mS_0,$
are contained in $F$ for $n\leq -N$, and 
do not intersect $\bar V_0$ for $n\geq N$.
The set $(V_0-F)\cap \mI(E_{c_0,\mS_0})$
is thus covered by finitely many pairwise disjoint compact sets
of the form 
$\pi\circ \phi^{n}(\tilde \mK(\mS_0,\mS)), n\in \Zm,\mS \in \Em(c_0)-\mS_0$.
As a consequence, each of the sets $\mK(\mS_0,\mS)$
is isolated in  $(V_0-F)\cap \mI(E_{c_0,\mS_0})$, 
and then also in 
 $\mI(E_{c_0,\mS_0})$.
Let us fix a static class $\mS_1\neq \mS_0$ such that 
 $\mK(\mS_0,\mS_1)$ is not empty. Such a static class 
exists by \ref{heteroclinics}.
Then, we can find an open neighborhood  $U\subset V_0$ of 
$\mK(\mS_0,\mS_1)$  such that 
$U$ is acyclic and such that 
 $U \cap \mI(E_{c_0,\mS_0})=\mK(\mS_0,\mS_1)$
is a non-empty compact set contained in 
$\mH(\mS_0,\mS_1)$.
We have proved $(ii)$. 

Let us consider again the finite family of pairwise disjoint compact
sets 
$\pi\circ \phi^{n}(\tilde \mK(\mS_0,\mS)), n\in \Zm,
|n|\leq N,\mS \in \Em(c_0)-\mS_0$.
There exists a finite family of pairwise disjoint compact sets
$\mK'_n(\mS_0,\mS), n\in \Nm,
|n|\leq N,\mS \in \Em(c_0)-\mS_0$
such that $\mK'_n(\mS_0,\mS)$  is a neighborhood of 
$\pi\circ \phi^{n}(\tilde \mK(\mS_0,\mS))$.
We can clearly assume in addition that the sets $\mK'_n(\mS_0,\mS)$
do not intersect $\mS_0$.
The set 
$$
V=V_0-\bigcup_{n\in \Nm,|n|\leq N, \mS \in \Em(c_0)-\mS_0}
\mK'_n(\mS_0,\mS)
$$
is an open neighborhood of $\mS_0$ which is contained in $V_0$,
and its boundary does not intersect $\mI(E_{c_0,\mS_0})$.
We have proved $(i)$.
\qed
\vs\\
The following proposition is the main step in the proof of the theorem.

\subsubsection{}\label{sousfinite}
\textsc{Proposition}
\begin{itshape}
Let $c_0$ satisfy the hypotheses of Theorem \ref{finite}.
For each weak KAM solution $\mG_0 \in \Vm_{c_0}$, there exists a number $\epsilon>0$
and an integer $N$ 
such that, if $\mG\in \Pm$ and $c\in H^1(M,\Rm)$ satisfy
 $$\|\mG-\mG_0\|\leq \epsilon
\text{ and }
|c-c_0|\leq \epsilon$$
then
$
\mG\fleche_N c.
$
\end{itshape}

\subsubsection{}
\textsc{Proof of the Theorem. }
We assume the Proposition.
Let us cover the compact set $\Vm_{c_0}$ by a finite number of balls
$B(\mG_i,\epsilon_i)$, where $\epsilon_i$ is given by
the Proposition applied to $\mG_i$.
Since the function $c$ restricted to $\Vm$ is proper,
the union of these finite balls covers the sets $\Vm_c$ for $c$
sufficiently close to $c_0$.
The Theorem holds by Proposition \ref{reduction}.
\mbox{}\hfill {\small \fbox{}} 
$ _{\textsc{theorem}} $
 {\small \fbox{}\\} 
\vs\\
We now prove the Proposition in three steps.

\subsubsection{}\label{S0}
\textsc{Step 1. }
\begin{itshape}
Let $\mG \in \Vm_{c_0}$ be a fixed point.
If there exist only finitely many static classes in $\mA(c_0)$,
then there exists an elementary solution $E_0$ and a neighborhood
$U_0$ of the corresponding  static class 
$\mS_0$ such that 
$
\mG_{|U_0}=E_{0|U_0}.
$
\end{itshape}

\proof
Let us fix the solution $\mG=\mG_{c_0,u}$.
We define a partial  order on the set of static classes by 
saying that $\mS\leq \mS'$ 
if, for each $x\in \mS$ and $x'\in \mS'$, we have
$h_{c_0}(x,x')=u(x')-u(x)$.
It is easy to check that this relation 
satisfies the following three axioms of order relations:
\begin{itemize}
\item
$\mS\leq \mS,$
\item 
$\mS\leq \mS' \text{ and } \mS'\leq \mS''
\Longrightarrow
\mS\leq \mS'',$
\item
$\mS\leq \mS' \text{ and } \mS'\leq \mS 
\Longrightarrow \mS=\mS'.$
\end{itemize}
Since the number of static classes is finite,
 there exists an initial element $\mS_0$,
that is an element which is not greater than any other element.
Let us write 
$$
u(x)=\min_{a\in \mA} u(a)+h_{c_0}(a,x),
$$
and consider, for each point $x$, the set
$\mA(c_0)\cap (\mG\wedge \breve E_{c_0,x})$
 of points $a$ where 
the minimum is reached.
Let us first assume that $x\in \mS_0$.
In this case, $a$ is a point of minimum if and only if 
the static class $\mS_a$ of $a$ satisfies $\mS_a \leq \mS_0$.
Since the class $\mS_0$ is initial, this implies
that $\mS_a=\mS_0$, or equivalently, that $a\in \mS_0$.
In other words, for $x\in \mS_0$, the compact set
$\mA(c_0)\cap (\mG\wedge \breve E_{c_0,x})$
does not intersect other static classes than $\mS_0$.
This implies that, for $x$ sufficiently close to $\mS_0$,
the set $\mA(c_0)\cap (\mG\wedge \breve E_{c_0,x})$
does not intersect other static classes than $\mS_0$.
Since, for each $x$, the set 
 $\mA(c_0)\cap (\mG\wedge \breve E_{c_0,x})$
contains a static class, we conclude that, for 
$x$ sufficiently close to $\mS_0$, we have 
$$
 \mS_0= \mA(c_0)\cap (\mG\wedge \breve E_{c_0,x}).
$$
As a consequence, we have, if $x$ is sufficiently close to 
$\mS_0$, 
$$
u(x)= u(a)+h_{c_0}(a,x),
$$
for each $a\in \mS_0$.
In other words, the difference $x\lmto h_{c_0}(a,x)-u(x)$ is the constant
$u(a)$ in a neighborhood of $\mS_0$.
\mbox{}\hfill {\small \fbox{}} 
$ _{\textsc{step 1}} $
 {\small \fbox{}\\}

\subsubsection{}\label{S1}
\textsc{step 2. }
\begin{itshape}
Let $\mS_0$ be a static class of $\mA(c_0)$
 and let $U_0$
be a neighborhood of $\mS_0$
satisfying  $(ii)$ of \ref{isolation}.
There exists a static class
$\mS_1$, an open  neighborhood 
$U_1$ of $\mS_1$  and,
for each $\delta>0$, a number  $\epsilon>0$ and an integer $N$
with the following property :
If $\mG\in \Pm$ satisfies $\|\mG-E_{c_0,\mS_0}\|_{U_0}\leq \epsilon$
and $c\in H^1(M,\Rm)$ satisfies $|c-c_0|\leq \epsilon$,
then there exists a pseudograph $\mG'\in \Pm_{c}$
such that $\|\mG'-E_{c_0,\mS_1}\|_{U_1}\leq \delta$
and 
$$
\mG_{|U_0}\fleche_N \mG'_{|U_1}
$$ 
\end{itshape}

\proof
There exists a static class $\mS_1$ and
an acyclic open set $U\subset  U_0-\mA(c_0)$
such that
$$
\mI(E_{c_0,\mS_0})\cap \bar U=\mI(E_{c_0,\mS_0})\cap U
$$ 
is a compact set $\mK\subset \mH(\mS, \mS_1)$.
Let us fix a point $x_0\in \mS_0$,
and denote by $u_0$ the function $h_{c_0}(x_0,.)$.
\vs\\
%
%
%
%
\subsubsection{}\label{step2lem1}
\textsc{Lemma. }
\begin{itshape}
There exists a neighborhood $U_1$ of $\mS_1$ such that
the equality 
$$
T^{\infty}_{c_0,U}u_0(y)=h_{c_0}(x_0,x_1)+h_{c_0}(x_1,y)
= h_{c_0}(x_0,x)+h_{c_0}(x,x_1) +h_{c_0}(x_1,y)
$$
holds for all 
$x\in \mK$, $y\in U_1$, and $x_1\in \mS_1$. 
As a consequence, we have 
$$
\Phi^{\infty}_U\big(E_{c_0,\mS_0}
\big)_{|U_1}
=E_{c_0,\mS_1|U_1}, $$
and the minimum in the definition of 
$T^{\infty}_{c_0,U}u_0(y)$ is not reached outside of $\mK$ when $y\in U_1$.
\end{itshape}
\vs\\
\proof
Let us set $v=T^{\infty}_{c_0,U}u_0$ for simplicity.
Recall, from \ref{uh}, that  
all  weak KAM solutions $v\in \Vm_{c_0}$ satisfy 
$v(y)=\min_{a\in \mA(c_0)}v(a)+h_{c_0}(a,y)$. 
Here, we obtain
\begin{equation}
\tag{$\diamond$}
v(y)=
\min_{x\in \bar U, a\in \mA(c_0)}
 h_{c_0}(x_0,x)+h_{c_0}(x,a) +h_{c_0}(a,y).
\end{equation}
We claim that, for $y\in \mS_1$, the set of minimizing pairs
$(x,a)$ is $\mK\times \mS_1$.
Indeed, if $(x,a)\in \mK\times \mS_1$, 
then $x\in E_{c_0,x_0}\wedge \breve E_{c_0,a}$,
so that 
$ h_{c_0}(x_0,x)+h_{c_0}(x,a) = h_{c_0}(x_0,a)$,
and
$$
 h_{c_0}(x_0,x)+h_{c_0}(x,a) +h_{c_0}(a,y)
=h_{c_0}(x_0,y)
=\min_{(z,z')\in M\times M}
 h_{c_0}(x_0,z)+h_{c_0}(z,z') +h_{c_0}(z',y).
$$
Hence we have
$$
 h_{c_0}(x_0,x)+h_{c_0}(x,a) +h_{c_0}(a,y)
=
\min_{(z,z')\in \bar U\times \mA(c)}
 h_{c_0}(x_0,z)+h_{c_0}(z,z') +h_{c_0}(z',y).
$$
We have proved that the pairs of  $\mK \times \mS_1$
are minimizing in the  equation ($\diamond$)
for $y\in \mS_1$.

Let us now prove that they are the only minimizing pairs.
A pair $(x,a)$ is minimizing if and only if
$
h_{c_0}(x_0,a)+h_{c_0}(a,y)=h_{c_0}(x_0,y)
$
and 
$
 h_{c_0}(x_0,x)+h_{c_0}(x,a)=h_{c_0}(x_0,a).
$
The second equality implies 
$$
x\in E_{c_0,\mS_0}\wedge \breve E_{c_0,\mS(a)}
\subset
\mI(E_{c_0,\mS_0}).
$$
Since $\mI(E_{c_0,\mS_0})\cap \bar U=\mK$, this implies $x\in \mK$.
If $x\in \mK$ and $a\in \mA(c_0)$, then
the equality
$
 h_{c_0}(x,a) = h_{c_0}(x,y)+h_{c_0}(y,a)
$
holds for all $y\in \mS_1$.
Indeed, let $x(n)=\pi \circ \phi^n(x,-\partial_1h(x,a))$ be the projection of the  orbit of the only point of 
$\breve E_{c_0,\mS(a)}$ above $x$.
 We have, for each $n\in \Nm$,
the equality of calibration by $-h_{c_0}(.,a)$:
$$
 A_{c_0}(0,x,n,x(n))+n\alpha(c_0) = h_{c_0}(x,a)-h_{c_0}(x(n),a).
$$
Let $n_k$ be an increasing sequence of integers such that 
the subsequence $x(n_k)$ has a limit $\omega\in \mS_1$.
Taking the liminf as $k\lto \infty$, we get 
$
h_{c_0}(x,\omega)\leq h_{c_0}(x,a)-h_{c_0}(\omega,a),
$
which implies the desired equality for $\omega$, and then
for all points of $\mS_1$.

Since  $(x,a)$ is a minimizing pair for  $v(y)$, we get,
by decomposing $h_{c_0}(x,a)$ in the expression of $v$,
$$
v(y)=h_{c_0}(x_0,x)+h_{c_0}(x,y)+
h_{c_0}(y,a)+h_{c_0}(a,y)
$$
and, since $v(y) \leq h_{c_0}(x_0,x)+h_{c_0}(x,y)$,
we finally obtain that 
$h_{c_0}(y,a)+h_{c_0}(a,y)\leq 0$
so that  $a\in \mS_1$.
We have proved the claim. In addition, we have proved,
for $x_1\in \mS_1$ and $x\in \mK$, the equality
$$
v(x_1)=
 h_{c_0}(x_0,x)+h_{c_0}(x,x_1)
=h_{c_0}(x_0,x_1).
$$
As a  consequence, for $y\in \mS_1$,
each point $a\in \mA(c)$
which is minimizing in the equation 
$$
v(y)=\min_{a\in \mA(c_0)}v(a)+h_{c_0}(a,y)
$$
belongs to $\mS_1$.
Since $\mS_1$ is isolated in $\mA(c)$, the conclusion holds
also for $y$ sufficiently close to $\mS_1$. We then have the
equality 
$$
v(y)=v(x_1)+h_{c_0}(x_1,y)
=h_{c_0}(x_0,x)+
h_{c_0}(x,x_1)+h_{c_0}(x_1,y)
$$
for all  $x_1\in \mS_1$ and $x\in \mK$
(and no other $x$ in $\bar U$).
\qed

\subsubsection{}
Applying \ref{phiufleche}, we get the existence of a positive $\epsilon'$
and of integers $N\leq N'$ such that 
each $\mG\in \Pm$ which satisfies
$\|\mG-E_{c_0,\mS_0}\|_U\leq \epsilon'$
also satisfies
$$
\mG_{|U}\fleche_{N'}
\Phi^{N,N'}_{U}(\mG)_{|U_1}.
$$
Since
$
\Phi^{\infty}_U\big(E_{c_0,\mS_0}
\big)_{|U_1}
=E_{c_0,\mS_1|U_1} $ in view of the lemma, and by \ref{convergence},
the integers $N$ and $N'$ can be chosen such that, in addition,
we have 
$$
\|\Phi^{N,N'}_U(E_{c_0,\mS_0}) -E_{c_0,\mS_1}\|_{U_1}\leq \delta/2.
$$
Reducing $\epsilon'$ if necessary, we can 
furthermore assume, by continuity of $\Phi_U^{N,N'}$, that
$$
\|\Phi^{N,N'}_{U}(\mG)-
\Phi^{N,N'}_{U}(E_{c_0,\mS_0})\|_{U_1}\leq \delta/2
$$
when $\|\mG-E_{c_0,\mS_0}\|_U\leq \epsilon'$, so that 
$$
\|\Phi^{N,N'}_{U}(\mG)-
 E_{c_0,\mS_1}   \|_{U_1}\leq \delta.
$$
Since $U$ is acyclic, for each cohomology $c$ and each pseudograph
$\mG$, there exists a pseudograph $\mG(c)$ which has cohomology $c$
and such that $\mG_{|U}=\mG(c)_{|U}$.
There exists a positive $\epsilon$ such that,
if $|c-c_0|\leq \epsilon$ and if 
$\|\mG-E_{c_0,\mS_0}\|_{U}\leq \epsilon$, then
we have 
$$
\|\mG(c)-E_{c_0,\mS_0}\|_{U}\leq \epsilon'.
$$
Note that this norm does not depend on the choice
of $\mG(c)$.
As a consequence, setting  $\mG'=\Phi^{N,N'}_{U}(\mG(c))$, we have 
$c(\mG')=c$, 
$$
\mG_{|U} =\mG(c)_{|U} \fleche_{N'}
\mG'_{|U_1}
$$
and 
$$
\|\mG'-
 E_{c_0,\mS_1}   \|_{U_1}\leq \delta.
$$
\mbox{}\hfill {\small \fbox{}} 
$ _{\textsc{step 2}} $
 {\small \fbox{}\\}

\subsubsection{}
\textsc{Step 3. }
\begin{itshape}
Let $\mS_1$ be a static  class in $\mA(c_0)$
satisfying  $(i)$ of \ref{isolation},
and let $U_1$ be a fixed neighborhood of $\mS_1$.
There exists a number $\delta>0$ and an integer $N$ such that,
if $\mG'\in \Pm$ satisfies 
$\|\mG'-E_{c_0,\mS_1}\|_{U_1}\leq \delta$,
then $\mG'_{|U_1}\fleche_N c(\mG')$
\end{itshape}
\vs\\
\proof
There exists 
an open  neighborhood $V_1\subset U_1$ of $\mS_1$ such that
$\mI(E_{c_0,\mS_1})\cap V_1=\mI(E_{c_0,\mS_1})\cap \bar V_1$.
(this is 
 $(i)$ of \ref{isolation}).
Let $x_1$ be a point of $\mS_1$ and set $u_1=h_{c_0}(x_1,.)$.
Recall that (by definition) 
$$
T^{\infty}_{c_0,V_1}u_1(x)=
\min_{y\in \bar V_1} h_{c_0}(x_1,y)+h_{c_0}(y,x).
$$
Taking $y=x_1$ in this expression, we obtain 
the inequality  $T^{\infty}_{c_0,V_1}u_1(x)\leq u_1(x)$.
On the other hand, we have the triangle inequality 
$u_1(x)\leq h_{c_0}(x_1,y)+h_{c_0}(y,x)$ for each $y$,
so that 
$
T^{\infty}_{c_0,V_1}u_1(x)=u_1(x),
$
and 
$$
\min_{y\in \bar V_1} h_{c_0}(x_1,y)+h_{c_0}(y,x)=
h_{c_0}(x_1,x)
=\min_{y\in  M} h_{c_0}(x_1,y)+h_{c_0}(y,x).
$$
By \ref{wedge} the points $y$ where this last minimum is reached 
belong to $\mI(E_{c_0,\mS_1})$.
As a consequence, for each $x\in M$,
the points where the minimum is
reached in the definition of 
$T^{\infty}_{c_0,V_1}u_1(x)$
belong to $\mI(E_{c_0,\mS_1})\cap V_1$,
which is a compact set contained in $V_1$.
In view of \ref{phiufleche},
there exist integers $N$ and $N'$ and a positive real number $\delta$
such that, if 
$\mG'\in \Pm$ satisfies 
$\|\mG'-E_{c_0,\mS_1}\|_{V_1}\leq \delta$,
then 
$$
\mG'_{|V_1}\fleche _N
\Phi^{N,N'}_{V_1}(\mG').
$$

\mbox{}\hfill {\small \fbox{}} 
$ _{\textsc{step 3}} $
 {\small \fbox{}\\} 
The proposition obviously follows from the three steps above.
\mbox{}\hfill {\small \fbox{}} 
$ _{\textsc{Proposition}} $
 {\small \fbox{}\\} 
%
%
%
%
%
%
%
%
%
%
%
%
%
\section{Applications}\label{application}
\subsection{Twist Maps}\label{TM}
The case where $M=\Tm$ is well known and have been studied
many times.
The resulting time-one flow is then a finite composition of
right twist maps of the biinfinite annulus $T^*\Tm$.
In view of \ref{orbites}, much of what is known on the existence
of orbits with prescribed behavior is summed up in the 
following discussion.

\subsubsection{}
Let $G\in H^1(\Tm, \Rm)$ be the set of cohomology classes 
of invariant curves which are Lipschitz graphs.
The set $G$ is closed,
and  every point $c\in G$ is alone in its class of $\ffleche$-equivalence,
as follows from \ref{torus}.
Conversely, if $c$ does not belong to $G$, then all the sets 
$\mI(\mG), \mG\in \Vm_c$ are properly contained in $\Tm$.
It follows that $R(\mG)=\Rm$ for each $\mG\in \Vm_c$, so that 
$R(c)=\Rm$ and, in view of \ref{mather}, $c$ is in the interior
of its class of equivalence.
The classes 
of $\ffleche$-equivalence are the points of $G$ and the connected components
of the complement of $G$.

\subsubsection{}
For completeness, we recall without proof some of the special properties
of Aubry sets in dimension one, see \cite{WE} for example.
The function $\alpha$ is differentiable, and its differential
$\alpha'(c)$ is the rotation number of every orbit of 
$\tilde\mN(c)$.
If $\alpha'(c)$ is irrational, then there is only one element 
in $\Vm_c$.
If $\alpha'(c)$ is rational, then the Mather set $\tilde \mM(c)$
is made of periodic orbits.
%
%
%
%
%
%
\subsection{Generalized Arnold Example}\label{GAE}
\subsubsection{}
In this application, we take
$$
M=\Tm\times N,
$$
where $N$ is a compact manifold of dimension $d-1$,
and denote by $q=(q_1,q_2)$ the points of $M$.
We assume that the homology group $H_1(N, \Zm)$
is not trivial. 
We denote the points of $TM$ by
$(q,v)=(q_1,q_2,v_1,v_2),$ where  $(q_1,v_1)\in T\Tm$
and $(q_2,v_2)\in TN$.
In the same way, we denote by
$(q,p)=(q_1,q_2,p_1,p_2)$ the points of $T^*M$.
We will consider the projection  $\pi_1:\Tm\times N\lto \Tm$ 
 and the induced mapping
$$
\pi_1^*:H^1(\Tm,\Rm)\lto 
H^1(\Tm\times N,\Rm).
$$

\subsubsection{}\label{hypinegalite}
Let us fix a point $0$ in $N$.
We will consider Lagrangian systems
which satisfy
$$
L(t,q_1,q_2,v_1,v_2)>L(t,q_1,0,v_1,0)
$$
for all $(q_2,v_2)\neq (0,0)$, all  $t\in \Rm$ and all
$(q_1,v_1)\in T\Tm$.
Let $\partial_v L:TM\lto \Tm^*M$ be the Legendre transform
associated to $L$.
We denote by $\Tm_1$
the submanifold $\Tm\times \{0\}$ of $M$,
 by  $T^*\Tm_1$  the submanifold
$\{q_2=0,p_2=0\}$ of $T^*M$, and 
$T\Tm_1$  the submanifold
$\{q_2=0,v_2=0\}$ of $TM$.
We have 
$$
\partial_vL(T\Tm_1)=T^*\Tm_1,
$$ 
and this manifold is invariant under the Hamiltonian flow.
Moreover, the restriction of the flow to $T^*\Tm_1$
is the Hamiltonian flow of the restriction 
$H_1(t,q_1,p_1):=H(t,q_1,p_1,0,0)$ of $H$. 
Setting $L_1(t,q_1,v_1)=L(t,q_1,0,v_1,0)$, we see that
$L_1$ is the Lagrangian associated to $H_1$.
We denote by  $\phi_1$  the restriction of $\phi$  to $T^*\Tm$.

\subsubsection{}
\textsc{Theorem }
\begin{itshape}
Under the non-degeneracy conditions \ref{internal}
and \ref{external} to be specified below, 
the image of $\pi_1^*$ is contained in one
class of $\ffleche$-equivalence.
\end{itshape}

\subsubsection{}\label{internal}
\textsc {Genericity property for $\phi_1$.}
We assume that every rotational invariant circle
of $\phi_1$ which contains a periodic orbit
is completely periodic (every orbit of this circle is periodic).
We could, more simply, require that the map $\phi_1$ does not have 
any invariant circle containing a periodic orbit. This
property is known to be generic in any reasonable sense of the term.
However, allowing periodic circles includes the important case where 
$\phi_1$ is integrable, as in the original Arnold's example.

\subsubsection{}\label{external}
\textsc{nondegeneracy of external homoclinics.}
We assume that, 
for each $c\in \pi^*_1(H^1(\Tm,\Rm))$,
there exists a finite Galois covering $P:M_0\lto M$ such that 
the set 
$$
\tilde \mN_{L\circ TP}(P^*(c))-T^*P^{-1}(T^*\Tm_1)
$$
is not empty and contains finitely many orbits.
Note that, since
$H^1(N,\Zm)$ is not zero,
it follows from \ref{coveringexists}, \ref{chains},
and   \ref{localisation} below 
that there exists a  finite Galois covering $P:M_0\lto M$
such that the set under consideration is not empty.
So the important point of our assumption is finiteness.
As the reader will see it in the proof, this assumption could 
be somewhat weakened.

\subsubsection{}\label{localisation}
\textsc{Lemma}
\begin{itshape}
For each cohomology $c= \pi_1^*(c_1)$, with $c_1\in H^1(\Tm,\Rm)$,
we have $\mN(c)\subset \Tm_1$.
As a consequence, the restriction to $\Tm_1$ gives a bijection
between the set $\Vm_c$  and the set $\Vm_{c_1}$
associated to the Lagrangian $L_1$ on $T\Tm_1$.
\end{itshape}
\vs\\
\proof
Let us fix a cohomology $c_1\in H^1(\Tm,\Rm)$
and its image $c:=\pi_1^*(c_1)$.
Let $\mu$ be a form on $\Tm$ which represents $c_1$,
and $\eta$ be its pull back on $M=\Tm\times N$.
Consider a pseudograph
$\mG\in \Vm_{c}$, and write it
$\mG=\mG_{\eta,u}$.
We want to prove that $\tilde \mI(\mG)\subset T^*\Tm_1$.
Let  $(q(t),p(t))$
be the trajectory of the Hamiltonian flow starting in $\tilde \mI(\mG)$.
We have, for $k<l$ in $\Zm$, 
$$
u(q(l))-u(q(k))
=
\int_k^l L(\sigma,q(\sigma),\dot q(\sigma)) 
-\mu_{q_1(\sigma)}(\dot q_1(\sigma)) +\alpha(c) \,d\sigma
$$
and
$$
u(q_1(t),0)-u(q_1(s),0)
\leq
\int_k^l L(\sigma,(q_1(\sigma),0,\dot q_1(\sigma),0) 
-\mu_{q_1(\sigma)}(\dot q_1(\sigma)) +\alpha(c) \,d\sigma.
$$
It follows that 
$$
\int_k^l L(\sigma,q(\sigma),\dot q(\sigma)) -
L(\sigma,(q_1(\sigma),0),(\dot q_1(\sigma),0)) 
\,d\sigma \leq 2(\max u-\min u)
$$
Let us denote by $\tilde L$
the function
$$\tilde L(t,q,v)=
L(t,q,v)-L(t,(q_1,0),(v_1,0))
$$
which is positive except on $T\Tm_1$.
Since the integral 
$\int_{\Rm} \tilde L(\sigma, q(\sigma),\dot q(\sigma))d\sigma$
 is finite, we have 
$$
\liminf _{|\sigma|\lto \infty}
\tilde L(\sigma,q(\sigma),\dot q(\sigma))=0,
$$
and consequently 
$
\liminf _{|\sigma|\lto \infty}
(q_2(t),v_2(t))=0.
$
We now return to the inequality
$$
\int_k ^l\tilde L(\sigma,q(\sigma),\dot q(\sigma))d\sigma
\leq u(q(t))-u(q_1(t),0)-u(q(s))+u(q_1(s),0),
$$
from which we get 
$$
\int_{-\infty} ^{\infty}\tilde L(\sigma,q(\sigma),\dot q(\sigma))d\sigma=0,
$$
which implies that 
$(q_2,v_2)\equiv 0$.
We have proved that $\tilde \mI(\mG)\subset T^*\Tm_1$.
\qed

\subsubsection{}
Let us fix cohomologies  $c=\pi_1^*(c_1), c_1\in H^1(\Tm,\Rm)$,
such that  there exists an invariant Lipschitz Graph $\mG$ in $\Vm_{c_1}$.
If the rotation number of $\phi_{1|\mG}$
is irrational, then $\Vm_{c_1}$ contains only one element.
As a consequence, $\Vm_c$ also contains only one   element,
 so that   $\tilde \mN(c)=\tilde \mA(c)=\mG$,
and there is only one static class in $\tilde \mA(c)$.
If the rotation number is rational, then in view of 
$\ref{internal}$ the graph $\mG$ is a union of periodic orbits,
so that $\mG=\tilde \mM(c)$.  
As a consequence, we have $\mA(c)=\Tm_1$, and there is only one static
class.

In view of \ref{external}, there exists a finite Galois covering
$P:M_0\lto M$ such that the Lagrangian $L\circ TP$ satisfies the
hypotheses of \ref{finite}.
As a consequence, the cohomology  $P^*(c)$ 
is in the interior of its forcing class for $L\circ TP$.
It follows from \ref{coveringrelation} that the 
cohomology $c$ is in the interior of its forcing  class 
 for $L$.

\subsubsection{}
Let  $c=\pi_1^*(c_1)$ be such that 
each set $\mI(\mG), \mG\in \Vm_c$ is properly contained in $\Tm_1$.
Applying  \ref{mather}, we observe that $R(c)=H^1(M,\Rm)$,
and $c$ is in the interior of its forcing  class.

\subsubsection{}
We have proved that each $c\in \pi_1^*(\Tm,\Rm)$ is in the interior
of its forcing class.
Since the subspace $ \pi_1^*(H^1(\Tm,\Rm))$ is obviously connected, 
it is contained in one forcing class.
\qed

%
%
%
%
%
%
%
%
%
%
%
\renewcommand{\thesubsection}{\Alph{subsection}}

\section{Appendix}
\setcounter{subsection}{0}

\subsection{Semi-concave functions}\label{semiconcave}
We recall some useful facts on semi-concave functions.
In all this section, $M$ is a compact manifold of dimension $d$.
It is useful for the sequel to fix once and for all 
 a finite atlas
$\Phi$ of $M$ composed of
charts  $\varphi:B_3\lto M$, where 
$B_r$ is the open  ball of radius $r$
centered at zero in $\Rm^d$.
We assume that the sets $\varphi(B_1),\varphi\in \Phi$ cover $M$.
A family $F$  of  $C^2$ functions
is said bounded if there exists a constant $C>0$
such that
$$
\|d^2(u\circ \varphi)_x\|\leq C
$$
for all $x\in B_1, \varphi\in \Phi, u\in F$.
Note that  a bounded family is not required
to be bounded in $C^0$ norm, but will automatically
be bounded in $C^1$ norm and thus equi-Lipschitz.
The notion of bounded family of functions
does not depend on the atlas $\Phi$.

\subsubsection{}
A function $u:M\lto \Rm$ is called semi-concave 
if there exists a bounded subset $F_u$ of the set
$C^2(M,\Rm)$  such that
$$
u=\inf_{f\in F_u}f.
$$
A family $\Um$ of functions is called equi-semi-concave
if there exists a bounded set $F$ of functions in
$C^2(M,\Rm)$, and, for each function $u\in \Um$,
a subset $F_u$ of $F$
such that
$$
u=\inf_{f\in F_u}f.
$$
Let us first collect some easy consequences of this definition.
We shall prove later that the infima could be replaced by 
minima in these definitions.

\subsubsection{}
\begin{itshape}
An equi-semi-concave set of functions is equi-Lipschitz.
\end{itshape}

\subsubsection{}
\begin{itshape}
If $\Um$ is an equi-semi-concave set of functions on $M$,
and if the infimum  $\inf _{u\in \Um}u(x_0)$ is finite
 for some  $x_0\in M$,
then the function  
$v(x)=\inf _{u\in \Um} u(x)$ is finite and semi-concave.
\end{itshape}

\subsubsection{}\begin{itshape}
Let $\Um$ be an equi-semi-concave 
set of functions on $N\times M$,
where $N$ is another compact manifold.
Then  the functions $u(x,.):M\lto \Rm,x\in N, u\in \Um$  form
an equi-semi-concave set.
\end{itshape}

\subsubsection{}
We say that the linear form $p\in T_xM$ is a
proximal super-differential
of the function $u$ at point $x$ if there exists a $C^2$
function $f$ 
such that $f-u$ has a minimum at $x$ and $df_x=p$.
The definition would not be changed by
requiring that the function $f$ is smooth and that the minimum is 
strict.
We say that a linear form $p\in T_xM$
is a $K$-super-differential of the function
$u$ at point $x$ if for each chart
$\varphi\in \Phi$ and each  point $y\in B_2$
satisfying
$\varphi(y)=x$, the inequality
$$
u\circ\varphi(z)-u\circ\varphi(y)
\leq p\circ d\varphi_y(z-y)+K\|z-y\|^2
$$
holds for each $z\in B_2$.
It is plain that $p$ is a proximal super-differential of $u$
if and only if there exists a $K>0$ such that
$p$ is a $K$-super-differential of $u$.

\subsubsection{}
A function $u$ on $M$ is called
$K$-semi-concave if it
has a $K$-super-differential at each point.
It is equivalent
to require that, for each $\varphi\in \Phi$,
the function
$$
u\circ \varphi(y)-K\|y\|^2
$$
is concave on $B_2$.
As a consequence, if $u$ is $K$-semi-concave
and if $p$ is a proximal super-differential of $u$
at $x$, then $p$ is a $K$-super-differential of
$u$ at $x$.

\subsubsection{}\label{regularity}
\begin{itshape}
Let $u$ be a continuous function on $M$, and let 
Let $A$ be the  compact subset of $M$
formed by points $x\in M$ 
at which both the functions $u$ and $-u$
have a $K$-super-differential.
Then the function $u$ is differentiable at each point of $A$,
and the mapping $x\lmto du(x)$ is Lipschitz on $A$, with a Lipschitz
constant that depends only on $K$.
\end{itshape}
\vs\\
This  follows from Proposition 4.5.3
in Fathi's book \cite{Fathibook}. We have the following useful 
corollary.

\subsubsection{}\label{sum}
\begin{itshape}
Let $\Um$ be an equi-semi-concave set of functions.
Then there exists a constant $K$ with the following property:
If $u$ and $v$ are two functions of $\Um$,
and if $A$ is the set of points minimizing the sum
$u+v$, then the functions $u$ and $v$ are differentiable 
at each point of $A$, and 
 the mapping $A\ni x\lmto du_x=-dv_x$
is $K$-Lipschitz.
\end{itshape}

\subsubsection{}\label{closed}
\begin{itshape}
If $u_n$ is a sequence of $K$-semi-concave functions
converging uniformly to a function $u$,
then the function $u$ is $K$-semi-concave.
In addition,
if $x_n$ is a sequence of points of differentiability of 
$u_n$, converging to a point of differentiability $x$ 
of $u$, then  $du_n(x_n)\lto du(x)$\vs.\\
\end{itshape}

\proof
By the Theorem below,
the functions $u_n$ form an equi-semi-concave, hence equi-Lipschitz
family of functions.
Let $x_n$ be a sequence converging to $x$,
and $p_n$ be a $K$-super-differential of $u_n$
at $x_n$.
The sequence $p_n$ is bounded, hence  we assume that $p_n\lto p$.
Let $y\in B_2$ and $\varphi\in \Phi$
be such that $\varphi(y)=x$.
For $n$ large enough, the point $x_n$ can be written
$\varphi(y_n)$ with $y_n\in B_2$.
We have the inequality
$$
u_n\circ \varphi(z)
\leq u_n\circ \varphi(y_n) +p_n\circ d\varphi_{y_n}(z-y_n)
+K\|z-y_n\|^2,
$$
for each $z\in B_2$,
and at the limit, we obtain 
$$
u\circ \varphi(z)
\leq u\circ \varphi(y) +p\circ d\varphi_{y}(z-y)
+K\|z-y\|^2.
$$
It means that $p$ is a $K$-super-differential of 
$u$ at $x$.
Under the assumptions of the statement, we have 
$p_n=du_{x_n}$, and $p=du_x$ is the only possible limit
of this bounded sequence, which is thus converging to $p$.
\qed

\subsubsection{}
\textsc{Theorem}
\begin{itshape}
A family $\Um $ of functions is equi-semi-concave if and only if
there exists a number  $K>0$ such that all the functions of $\Um$
are $K$-semi-concave.
In this case,  there  exists a bounded subset $F\subset C^2(M,\Rm)$
and, for each $u\in \Um$, a subset $F_u$ of $F$
which has the following properties:
$$
u=\min _{f\in F_u} f
$$
and, for each point $x\in M$ and each super-differential
$p$ of $u$ at $x$, there exists a function $f\in F_u$
such that $(f(x),df(x))=(u(x),p)$.\vs\\
\end{itshape}

In order to prove this result, we need a Lemma:\vs\\

\noindent
\textsc{Lemma }
\begin{itshape}
For each $K>0$, there exists a bounded subset $L_K$
of $T^*M$ which contains all the proximal super-differentials
of all $K$-semi-concave functions.
As a consequence, the $K$-semi-concave functions
are equi-Lipschitz\vs\\
\end{itshape}

\noindent
\textsc{proof of the Lemma}
Let us consider a chart $\varphi\in \Phi$,
a $K$-semi-concave function $u$,
a point $y_0\in B_1$, and  the point $x_0=\varphi(y_0)$.
Let $p_0$ be a proximal super-differential of $u$
at $x_0$, and let us set $l=p_0\circ d\varphi_{y_0}$.
\vs\\
\noindent
\textsc{claim}
\begin{itshape}
If $\|l\|\geq 11K$, then there exists a point $y\in B_2$
which is a point of differentiability of $u\circ \varphi$
and satisfies 
$$\|d(u\circ \varphi)_{y}\|\geq (\|l\|-11K)/3.$$
and
$$
u\circ \varphi(y)<\inf_{B_1} u\circ \varphi.
$$
\end{itshape}
\vs\\
\textsc{Proof of the Claim. }
Let us prove   first that the infimum of 
$u\circ \varphi$ in $B_2$ is not reached in $\bar B_1$.
Assume, by contradiction, that there exists a point
$m\in \bar B_1$ such that 
$u\circ \varphi(m)=\inf _{B_2} u\circ \varphi$.
Then clearly the function $u\circ \varphi$ is differentiable at 
$m$, its differential is zero,  and the inequality
$$
u\circ \varphi(y_0)\leq u\circ \varphi(m)+K\|y_0-m\|^2
$$
holds.
On the other hand, we have 
$$
u\circ \varphi(m)\leq u\circ \varphi(z) \leq u\circ \varphi(y_0)
+ l(z-y_0)+K\|z-y_0\|^2
$$
for all $z\in B_2$.
Combining these inequalities gives 
$$
l(y_0-z)\leq K\|z-y_0\|^2+K\|y_0-m\|^2
$$
for all $z\in \bar B_2$.
Hence $\|l\|\leq 5K$, which is in contradiction with the hypothesis.

Let us now
consider a vector $v\in \Rm^d$ of norm $1$ and such that 
$l(v)= -\|l\|$. 
We get 
$$
u\circ \varphi(y_0+v)-
u\circ \varphi(y_0)\leq l(v)+K\|v\|^2
= K-\|l\|.
$$
Hence the infimum of $u\circ \varphi$ on $B_2$ is not greater
than $u \circ \varphi(y_0)+K-\|l\|$.
It is then possible to choose a point $y$ in $B_2$ 
such that 
$$u\circ \varphi(y)<
\min \big ( \inf _{ B_1}u\circ \varphi,
 u \circ \varphi(y_0)+2K-\|l\|\big).
$$
In addition, since the function $u\circ \varphi$ 
is differentiable almost everywhere, we can 
assume that the function $u\circ \varphi$ is differentiable at $y$.
We have the inequality
$$
u\circ \varphi(y_0)\leq 
u\circ \varphi(y)+d(u\circ \varphi)_{y}(y-y_0)+K\|y-y_0\|^2
$$
from which follow 
$$
d(u\circ \varphi)_{y}(y_0-y)\leq 
u\circ \varphi(y)-u\circ \varphi(y_0)+
K\|y-y_0\|^2
\leq 11K-\|l\|.
$$
Hence $\|d(u\circ \varphi)_{y}\|\geq (\|l\|-11K)/3$.
This ends the proof of the claim. 

\mbox{}\hfill {\small \fbox{}} 
$ _{\textsc{claim}} $
 {\small \fbox{}\\}

In order to continue the proof of the Lemma,
we consider the point $y\in B_2$
given by the claim.
There exists a chart $\varphi_1\in \Phi$  and a point
$y_1\in B_1$ such that $\varphi_1(y_1)=\varphi(y)=:x_1$.
Note that $u\circ \varphi_1$ is differentiable at 
$x_1$, and define
$$
l_1:= d(u\circ \varphi_1)_{y_1}
= d(u\circ \varphi)_{y}\circ
d(\varphi^{-1}\circ \varphi_1)_{y_1}.
$$
There exists a constant $C>1$, which depends
only on the atlas $\Phi$, and such that 
$$
\|l_1\|\geq (\|l_0\|-11K)/C.
$$
If $l_0$ is large enough, then we have 
$\|l_1\|\geq 11K$, hence we can apply the lemma again,
and find a chart $\varphi_2$, a point $y_2\in B_2$ and a linear 
form $l_2$.
In addition, we have 
$$
u\circ \varphi_2(y_2)<\inf _{\varphi(B_1)\cup \varphi_1(B_1)}
u,
$$
so that the charts $\varphi$, $\varphi_1$ and $\varphi_2$ 
are different.
Now if $\|l_0\|$ is sufficiently large, the process
can be continued further and we can build
inductively, for $0\leq i\leq N$, a sequence
$x_i\in B_1$ of points, a sequence $\varphi_i\in \Phi$
of different charts, and a sequence $l_i$ of linear forms
such that 
$
\|l_{i+1}\|\geq (\|l_i\|-11K)/C.
$
The process can be continued as long as $\|l_i\|\geq 11K$.
Recall that the cardinal of $\Phi$ is finite, and denote it by
$|\Phi|$.
Since all the charts involved in the construction above
are different, at most $|\Phi|$ steps 
can be performed. Hence there exists an integer $N\leq |\Phi|$
such that 
$\|l_i\|\geq 11K$ for $i <N$, and $\|l_N\|\leq 11K$.
This gives a bound to $\|l\|$, hence to $\|p\|$.

\mbox{}\hfill {\small \fbox{}} 
$ _{\textsc{lemma}} $
 {\small \fbox{}\\}

We now finish the proof of the Theorem.
Let us consider a smooth function 
$g:\Rm^d\lto \Rm$ such that 
$0\leq g \leq 1$, and such that
$g=0$ outside of $B_2$ and $g=1$ inside $B_1$.
Since the $K$-semi-concave functions are equi-Lipschitz,
and since the manifold $M$ is compact, there 
exists a number $\Delta>0$ such that
$$
\max u-\min u\leq \Delta
$$
for each $K$-semi-concave function $u$.
Let us associate, to each chart $\varphi\in \Phi$,
and each point $(x,p)\in T_xM$
satisfying $x\in \varphi(B_1)$,
 the function $f_{x,p,\varphi}:M\lto \Rm$
defined by
$$
f_{x,p,\varphi}\circ \varphi(z)
:=g(z)\big(p\circ d\varphi_y(z-y)+K\|z-y\|^2\big)
+(1-g(z))\Delta
$$
for $z\in B_2$, where $y=\varphi^{-1}(x)$,
and 
$f_{x,p,\varphi}=\Delta$
outside of  $\varphi(B_2)$.
The functions $f_{x,p,\varphi}, (x,p)\in L_K, \varphi\in \Phi$
 form a bounded subset $F$ of $C^2(M,\Rm)$.
For each $K$-semi-concave function $u$,
let $F_u\subset F+\Rm$ be the set of functions 
$$
z\lmto f_{x,p,\varphi}(z)+u(x)
$$
where $p$ is a $K$-super-differential
of $u$ at $x$.
We claim that $u=\min_{f\in F_u} f$.
In order to prove this claim, observe that,
for each $y\in B_1$, $\varphi\in \Phi$ and 
$p$ a $K$-super-differential of $u$ at $x=\varphi(y)$,
we have   
$$f_{x,p,\varphi}-u(x)\geq u$$
with equality at $x$.
Indeed, we have the inequalities
$$
u\circ \varphi(z)-u(x)\leq p\circ d\varphi_x(z-y)+K\|z-y\|^2
$$
for $z\in B_2$
and 
$
u\leq u(x)+\Delta.
$
\mbox{}\hfill {\small \fbox{}} 
$ _{\textsc{theorem}} $
 {\small \fbox{}\\}

%
%
%
\subsection{Uniform families of Hamiltonians}\label{uniform}
Let us fix once and for all a Riemann metric
on the compact manifold $M$. We use this metric
to define a norm $|v|$ for tangent vectors,
and a norm $|p|$ for tangent covectors.

\subsubsection{}
A family of pairs $(H,L)$ of dual Hamiltonians and Lagrangians
satisfying the hypotheses \ref{HypothesesH} and \ref{HypothesesL}
is called uniform if:
\begin{itemize}
\item[$(i)$]
There exist two superlinear functions $h_0$ and $h_1:\Rm^+\lto \Rm$
such that each Hamiltonian $H$ of the family satisfies 
$h_0(|p|)\leq H(t,x,p)\leq h_1(|p|).$
\item[$(ii)$]
There exists an increasing function $K(k): \Rm^+\lto \Rm^+$ such 
that, if   $\phi $ 
is the flow of a Hamiltonian of the family and if 
the times $t$ and $s$ satisfy 
$t-1\leq s\leq t+1$,
then 
$$
\phi_t^s\big( \{|p|\leq k\}\big)
\subset \{|p|\leq K(k)\}\subset T^*M.
$$
\item[$(iii)$]
There exists a finite atlas $\Psi$ of 
$M$ such that, for each chart $\psi\in \Psi$ and each
Lagrangian $L$ of the family,
we have $\|d^2(L\circ T\psi)_{(t,x,v)}\|\leq K(k)$
for $|v|\leq k$.
\end{itemize}
Note that condition $(i)$ could have equivalently been 
replaced by the following:
\begin{itemize}
\item[$(i')$]
There exist two superlinear functions $l_0$ and $l_1:\Rm^+\lto \Rm$
such that each Lagrangian  $L$ of the family satisfies 
$l_0(|v|)\leq L(t,x,v)\leq l_1(|v|).$
\end{itemize}

\subsubsection{}
The uniform families of highest use are the following.
If $H$ is a Hamiltonian,
and if $\omega$ is a bounded finite-dimensional convex
family of closed one-forms on $M$,
then the Hamiltonians 
$H(t,x,p+\omega_x),\omega\in \Omega$ form
a uniform family.
Equivalently, the Lagrangians 
$L(t,x,v)-\omega_x(v)$ form a uniform family.

\subsubsection{}\label{proprete}
\begin{itshape}
In a uniform family, we have
$$
|\partial_p H(t,x,p)|\geq \frac{h_0(|p|)-h_1(0)}{|p|}
$$
and
$$
|\partial_v L(t,x,v)|\geq \frac{l_0(|v|)-l_1(0)}{|v|}.
$$
In other words, the Legendre transforms are uniformly proper.
\end{itshape}
\vs\\
\proof
In view of the convexity of $H$, we have
$$
|\partial_p H_{(t,x,p)}|\geq \frac {H(t,x,p)-H(t,x,0)}{|p|}.
$$
\qed

\subsubsection{}
Given a Lagrangian $L$ satisfying the hypotheses of \ref{HypothesesL},
we define the function $A_L(t,x;s,y):\Rm\times M\times \Rm\times M\lto \Rm$
 by
$$
A_L(t,x;s,y)=
\inf_{\gamma\in \Sigma( t,x;s,y)}
\int _t^s L(\sigma,\gamma(\sigma),\dot\gamma(\sigma))\, d\sigma,
$$
Where  $\Sigma( t,x;s,y)$ is the set of absolutely continuous curves
$\gamma:[s,t]\lto M$
satisfying  $\gamma(t)=x$ and $\gamma(s)=y$.
We denote by  $\Sigma_m^L( t,x;s,y)$ the set of curves of 
 $\Sigma( t,x;s,y)$ which realize the minimum.

\subsubsection{}\label{boundedspeed}
\begin{itshape}
For each uniform family of Lagrangians,
there exists a decreasing function $K_1(\epsilon):]0,\infty)\lto\Rm^+$
such that, if $L$ is a Lagrangian of the family
and  if $t$ and $s$
are two real times satisfying $t\geq s+\epsilon$,
then each curve $\gamma\in \Sigma_m^L( t,x;s,y)$  is $C^2$ and
satisfies  
$|\dot \gamma(\sigma)|\leq K(\epsilon)$ for each $\sigma\in [s,t]$.
\end{itshape}
\vs\\
\proof
Without loss of generality, we can assume that
$0<\epsilon<t-s<1$.
By comparing the action of $\gamma$ with that of a geodesic
with the same endpoints, we get 
$$
\int_s^t l_0(|\dot \gamma(\sigma)|)d\sigma
\leq \int_s^t L(\sigma, \gamma(\sigma),\dot \gamma(\sigma))d\sigma
\leq (t-s)l_1\left( \frac{\text{diam}(M)}{t-s}
\right)
$$
The right hand side is clearly bounded by a constant which
depends only of the parameters of the uniform family and of 
$\epsilon$.
We obtain
$$
(t-s)\min l_0(|\dot \gamma(\sigma)|)\leq C,
$$
from which follows, with another constant $C$, that 
$
\min |\dot \gamma(\sigma)|\leq C.
$
But then in view of \ref{proprete}, we have 
$$
\min_{\sigma \in [s,t]}
 |\partial_vL(\sigma,\gamma(\sigma),\dot \gamma(\sigma))|\leq C,
$$
then  in view of $(ii)$,
$$
\max_{\sigma \in [s,t]}
 |\partial_vL(\sigma,\gamma(\sigma),\dot \gamma(\sigma))|\leq C,
$$
so that finally, using   \ref{proprete} again, we get 
$\max |\dot \gamma|\leq C$.
We have used the symbol $C$ for different constants which depend only
of $\epsilon$ and of the parameters of the family.
\qed
Note that the proof does not use $(iii)$ in the definition of
uniform families.

\subsubsection{}
\begin{itshape}
For each times $s<t$, the mapping which, to a Lagrangian
$L$,
associates the function 
$$
(x,y)\lmto A_L(s,x,t,y)
$$
of $C(M \times M,\Rm)$, is continuous on each uniform
family of Lagrangians endowed with the topology of uniform
convergence on compact sets.
\vs\end{itshape}

\proof
Let $L_0$ and $L_1$ be two Lagrangians of the family.
Let $\gamma(\sigma):[s,t]\lto M$ be such that
$$
A_{L_0}(s,\gamma(s);t,\gamma(t))
=
\int_s ^t L_0(\sigma,\gamma(\sigma),\dot\gamma(\sigma))\, d\sigma.
$$
We have 
$$
A_{L_1}(s,\gamma(s);t,\gamma(t))
\leq
\int_s ^t L_1(\sigma,\gamma(\sigma),\dot\gamma(\sigma))\, d\sigma,
$$
so that 
$$
A_{L_1}(s,\gamma(s);t,\gamma(t))-
A_{L_0}(s,\gamma(s);t,\gamma(t))
\leq
(t-s)\max_{|v|\leq K_1(t-s)} L_1-L_0,
$$
where $K_1$ is defined in \ref{boundedspeed}.
By symmetry, we get  that
$$
\|A_{L_0}(s,.;t,.)-A_{L_1}(s,.;t,.)\|_{\infty}
\leq (t-s)\max_{|v|\leq K_1(t-s)} |L_1-L_0|.
$$
\qed

\subsubsection{}
\textsc{Theorem. }
\begin{itshape}
For each uniform family of Lagrangians and each
$\epsilon>0$, consider  the set $\Um_{\epsilon}$ of continuous 
 functions $M\times M\lto \Rm$ given by
$$
(x,y)\lmto A_L(s,x;t,y)
$$
where $t\geq s+\epsilon$ and $L$ is a Lagrangian 
of the family.
This set is equi-semi-concave, hence equi-Lipschitz
 on $M\times M$. 
In addition, for each curve $\gamma\in \Sigma^L_m(t,x;s,y)$,
the covector
$$
\big( -\partial_vL(s,\gamma(s),\dot \gamma(s)),
\partial_vL(t,\gamma(t),\dot \gamma(t))  \big)
$$
is a proximal super-differential of the function 
$A_L(s,.;t,.)$ at point $(x,y)$.
\end{itshape}
\vs\\
\proof
Let us consider a finite atlas $\Psi$ of $M$ formed by charts
$\psi:B_6^d\lto M$, where $B_r^d$ is the Euclidean ball of radius 
$r$ in $\Rm^d$. Assume that the open sets $\psi(B^d_{1/2})$, 
$\psi\in \Psi$, cover $M$.
Let $\Phi$ be the atlas of $M\times M$ composed of products 
$\psi \times \psi'$, with $\psi\in \Psi$ and $\psi'\in \Psi'$. 
The charts $\varphi$ of $\Phi$ are defined on $B^{2d}_3$,
and the images $\varphi(B^{2d}_1)$, $\varphi\in \Phi$, cover $M\times M$.
In order to prove that the set $\Um_{\epsilon}$ is equi-semi-concave,
we shall check that it is $K$-semi-concave for some $K$.
So from now on we shall work in a fixed chart 
$\varphi=\psi_0\times \psi_1$.

Let $(x_0,x_1)$ be a point in $\psi_0(B_2)\times \psi_1(B_2)$,
and let $y_0$ and $y_1$ be the preimages in $B_2$.
Let $\gamma(t): [s,t]\lto M$ be a curve in 
$\Sigma_m(s,x_0;t,x_1)$.
In view of \ref{boundedspeed}, we have 
$|\dot \gamma|\leq K_1(\epsilon)$.
As a consequence, there exists a constant $a>0$, which
depends only on the atlas, on the parameters of the family,
and  of $\epsilon$,  such that
the curve 
$\psi_0^{-1}\circ \gamma:[s,s+1/a]\lto B^d_4$ 
is well defined and $a$-Lipschitz, 
as well as the curve 
$\psi_1^{-1}\circ \gamma:[t-1/a,t]\lto B^d_4$.
Let us call $y_0(\sigma)$ and $y_1(\sigma)$ these curves,
note that $y_0(s)=y_0$ and $y_1(t)=y_1$.
Let us now define, for each points $z_0$ and $z_1$ in $B_4$,
the curves
$$
y_0(\sigma,z_0):=
y_0(\sigma)+(1+a(s-\sigma))(z_0-y_0)
$$
and
$$
y_1(\sigma,z_1):=
y_1(\sigma)+(1+a(\sigma-t))(z_1-y_1).
$$
For simplicity
we define the Lagrangians $L_0$ and $L_1$ on 
$\Rm\times B^d_4\times \Rm^d$ by
the expression
$L_i(\sigma, x,v)=L(\sigma, \psi_i(x), d\psi_{ix}(v))$,
shortly, $L_i=L\circ T\psi_i$.
We have 
$$
A(s, \psi_0(z_0);t,\psi_1(z_1))
$$
$$
\leq A(s,x_0;t,x_1)+
\int_s^{s+1/a}
L_0(\sigma, y_0(\sigma,z_0), \dot y_0(\sigma,z_0))
-L_0(\sigma,y_0(\sigma),\dot y_0(\sigma))\,d\sigma
$$
$$+
\int_{t-1/a}^t
L_1(\sigma, y_1(\sigma,z_1), \dot y_1(\sigma,z_1))
-L_1(\sigma,y_1(\sigma),\dot y_1(\sigma))\,d\sigma.
$$
There exists a constant $C>0$, which
depends only on the atlas, on the parameters of the family,
 of $\epsilon$, and of $a$,  such that,
for $(t,x,v)\in \Rm \times B^{d}_4\times B^{d}_a $ and 
$(y,w)\in \Rm\times  B^{d}_4\times B^{d}_a $, we have 
$$
L_i(\sigma, y,w)-L_i(\sigma, x,v)\leq
\partial _x L_{i(\sigma, x,v)}(y-x)+
\partial _v L_{i(\sigma, x,v)}(w-v)+
C(\|y-x\|^2 +\|w-v\|^2).
$$
We  get 
$$
A(s, \psi_0(z_0);t,\psi_1(z_1))
\leq  A(s,x_0;t,x_1)
$$
$$
+\int_s^{s+1/a}
\partial _x L_{0(\sigma, y_0(\sigma),\dot y_0(\sigma))}
(y_0(\sigma,z_0)-y_0(\sigma))+
\partial _v L_{0(\sigma, y_0(\sigma),\dot y_0(\sigma))}
(\dot y_0(\sigma,z_0)-\dot y_0(\sigma))\,d\sigma
$$
$$
+\int_{t-1/a}^t
\partial _x L_{0(\sigma, y_0(\sigma),\dot y_0(\sigma))}
(y_0(\sigma,z_0)-y_0(\sigma))+
\partial _v L_{0(\sigma, y_0(\sigma),\dot y_0(\sigma))}
(\dot y_0(\sigma,z_0)-\dot y_0(\sigma))\,d\sigma
$$
$$
+C
\int_s^{s+1/a}\| y_0(\sigma)-y_0(\sigma,z_0)\|^2
+\| \dot y_0(\sigma)-\dot y_0(\sigma,z_0)\|^2
\,d\sigma
$$
$$
+C
\int_{t-1/a}^t\| y_1(\sigma)-y_1(\sigma,z_1)\|^2
+\| \dot y_1(\sigma)-\dot y_1(\sigma,z_1)\|^2
\,d\sigma.
$$
Taking advantage of the Euler-Lagrange equations, this simplifies to
$$
A(s, \psi_0(z_0);t,\psi_1(z_1))
\leq  A(s,x_0;t,x_1)
-\partial_v L_{0(\sigma, y_0, \dot y_0(s))}(z_0-y_0)
+\partial_v L_{1(\sigma, y_1, \dot y_1(t))}(z_1-y_1)
$$
$$
+C\frac{1+a^2}{a}
(\|y_0-z_0\|^2+\|y_1-z_1\|^2) .
$$
\qed

\subsubsection*{Acknowledgements. }\hspace{,1cm}
I thank Marie-Claude Arnaud, Alain Chenciner and Albert Fathi
for their comments which helped me to improve the present text.

\begin{small}

\end{small}
\end{document}